\documentclass[a4paper,10pt, openany]{article}

\usepackage[cp1252]{inputenc}
\usepackage{amsthm}
\usepackage{amsmath}
\usepackage{amssymb}
\usepackage{amsfonts}
\newcommand{\stackunder}[2]{\underset{#1}{#2}}
\newcommand{\R}{\mathbb{R}}

\newcommand{\ve}{\varepsilon}

\newcommand{\g}{\mathbf{g}}
\newcommand{\E}{(E_{h,f,\mathbf{g}})}

\newtheorem{definition}{Definition}

\newtheorem*{theo}{Theorem}

\newtheorem{theorem}{Theorem}

\newtheorem{proposition}{Proposition}

\begin{document}
\title{Critical functions and elliptic PDE on compact riemannian manifolds}
\author{Stephane Collion\\ \sl \small Institut de Math{\'e}matiques, Universit{\'e} Paris VI, Equipe
  G{\'e}om{\'e}trie et Dynamique, \\ \sl \small 175 rue Chevaleret, 75013
  Paris.  \\ \sl \small email: collion@math.jussieu.fr}
\date{}
\maketitle
\begin{abstract}
We study in this work the existence of minimizing solutions to the critical-power type equation 
$\triangle _{\textbf{g}}u+h.u=f .u^{\frac{n+2}{n-2}} $ on a compact riemannian manifold in the 
limit case normally not solved by variational methods. For this purpose, we use a concept of "critical 
function" that was originally introduced by E. Hebey and M. Vaugon for the study of second best 
constant in the Sobolev embeddings. Along the way, we prove an important estimate concerning concentration 
phenomena's when $f$ is a non-constant function. We give here intuitive details.
\end{abstract}

\section{Introduction}
In the beginning was the Yamabe problem:\\

\textbf{Yamabe problem}\textit{: Given a compact riemannian manifold }$(M,\textbf{g})$\textit{\ of dimension }$%
n\geq 3$,\textit{\ does there exist a metric} $\textbf{g'}$ \textit{conformal to }$\textbf{g}$\textit{\ having constant scalar curvature?}\\

If we write $\textbf{g'}=u^\frac{4}{n-2}.\textbf{g}$ where $u>0$ is a smooth function on $M$, the scalar 
curvatures are linked by the partial differential equation :
$$\triangle _{\textbf{g}}u+\frac{n-2}{4(n-1)}S_{\textbf{g}}.u=\frac{n-2}{4(n-1)}S_{\textbf{g'}} .u^{\frac{n+2}{n-2}} $$
where $S_{\textbf{g}}$ is the scalar curvature of $\textbf{g}$ and where $\triangle _{\textbf{g}}=-\nabla^{i}\nabla_{i}$
 is the riemannian laplacian of $\g$.

To solve the Yamabe problem, one therefore has to prove the existence of a solution $u>0$ to this partial differential equation
when $S_{\textbf{g'}}$ is a constant. More generaly, the prescribed curvature problems, which consist in deciding, 
given a smooth function $f$ on $M$, if $f$ is the scalar curvature of a metric conformal to $\g$, come down
to prove the existence of a positive smooth solution $u$ to the above equation when $S_{\g'}$ is replaced by $f$.

These problems launched the study of elliptic PDE on compact riemannian manifolds of the form 
$$(E_{h,f,\mathbf{g}}): \triangle _{\textbf{g}}u+h.u=f .u^{\frac{n+2}{n-2}} $$
In all this paper $M$ will be a compact riemannian manifold of dimension $n\geqslant 3$, we will use 
the letter $\g$ or $\g'$ to denote a riemannian metric on $M$; $h$ and $f$ will always be smooth functions on 
$M$. We will always suppose the functions to be smooth, however in the definitions and in most of the 
theorems, continuity is in general sufficient. Beside, we will keep these notations, letter $\g$ for the metrics,
letter $h$ for the function on the left of equation $E_{h,f,\g}$, (defining the opperator $ \triangle _{\textbf{g}}+h$);
and letter $f$ for the function on the right of the equation; the unknown function will be designated by $u$.

One of the possible methods to study these equation is the use of variational methods, which have the 
advantage of giving minimizing solutions, or solution of minimal energy. If one multiply equation $\E$ 
by $u$ and integrate over $M$, one gets
$$\int_{M}\left| \nabla u\right| _{\mathbf{g}%
}^{2}dv_{\mathbf{g}}+\int_{M}h.u{{}^{2}}dv_{\mathbf{g}}= \int_{M}f\left| u\right| ^{\frac{2n}{n-2}%
}dv_{\mathbf{g}}$$
The variational methods therefore lead to consider the functional
$$I_{h,\mathbf{g}}(w)=\int_{M}\left| \nabla w\right| _{\mathbf{g}%
}^{2}dv_{\mathbf{g}}+\int_{M}h.w{{}^{2}}dv_{\mathbf{g}} $$
defined for $w\in H_{1}^{2}(M)$, the Sobolev space of $L^{2}$ functions whose gradient is also in $L^{2}$, 
and the minimum of this functional
\[
\lambda _{h,f,\mathbf{g}}=\stackunder{w\in \mathcal{H}_{f}}{\inf }I_{h,%
\mathbf{g}}(w) 
\]
on the set
\[
\mathcal{H}_{f}=\{w\in H_{1}^{2}(M)/\int_{M}f\left| w\right| ^{\frac{2n}{n-2}%
}dv_{\mathbf{g}}=1\}. 
\]
The Euler equation associated with the minimization problem of this functional by a function $u$ such that
\[
I_{h,\mathbf{g}}(u)=\stackunder{w\in \mathcal{H}_{f}}{\inf }I_{h,\mathbf{g}%
}(w) 
\]
is indeed exactly
\[
(E_{h,f,\mathbf{g}}):\triangle _{\mathbf{g}}u+hu=\lambda _{h,f,\mathbf{g}%
}.f.u^{\frac{n+2}{n-2}} 
\]
where $\lambda _{h,f,\mathbf{g}}\,$ appears as a normalizing constant due to the condition
\[
\int_{M}f\left| u\right| ^{\frac{2n}{n-2}}dv_{\mathbf{g}}=1. 
\]
It is sometimes usefull to consider the functional 
\[
J_{h,f,\mathbf{g}}(w)=\frac{\int_{M}\left| \nabla w\right| _{\mathbf{g}%
}^{2}dv_{\mathbf{g}}+\int_{M}h.w{{}^{2}}dv_{\mathbf{g}}}{\left(
\int_{M}f\left| w\right| ^{\frac{2n}{n-2}}dv_{\mathbf{g}}\right) ^{\frac{n-2%
}{n}}} 
\]
and the subset of $H_{1}^{2}(M)$ where it is defined
\[
\mathcal{H}_{f}^{+}=\{w\in H_{1}^{2}(M)/\int_{M}f\left| w\right| ^{\frac{2n}{%
n-2}}dv_{\mathbf{g}}>0\}. 
\]
One then consider the minimisation problem by a function $u$ such that
\[
J_{h,f,\mathbf{g}}(u)=\stackunder{w\in \mathcal{H}_{f}^{+}}{\inf }J_{h,f,\mathbf{g}}(w), 
\]
the Euler equation being identical but without the normalizing constant. This functional sometimes 
present the advantage of being homogeneous in the sense that $J_{h,f,\mathbf{g}}(c.w)=J_{h,f,%
\mathbf{g}}(w)$ for any constant $c.$ One therefore see that
\[
\stackunder{w\in \mathcal{H}_{f}}{\inf }I_{h,\mathbf{g}}(w)=\stackunder{w\in 
\mathcal{H}_{f}^{+}}{\inf }J_{h,f,\mathbf{g}}(w)=\lambda _{h,f,\mathbf{g}} 
\]
This functional $J$ also has the particularity, when $h=\frac{n-2}{4(n-1)}S_{\mathbf{g}},$ 
of being invariant by conformal changes of metrics; it is therefore especially usefull when studying 
problems of prescribed scalar curvatures. We shall mostly use  $I_{h,\mathbf{g}}$ and $\mathcal{H}_{f}$, 
but for some problems $J_{h,f,\mathbf{g}}$ will prove to be more convenient when we shall want
 to avoid the constraint $\int_{M}f\left| u\right| ^{\frac{2n}{n-2}}dv_{\mathbf{g}}=1.$

We will say that a function $u\in H_{1}^{2}(M)$ is a solution of minimal energy, or a minimizing 
solution, if either $I_{h,\mathbf{g}}(u)=\lambda _{h,f,\mathbf{g}}$ with $\int_{M}fu^\frac{2n}{n-2}=1$, 
or $J_{h,f,\mathbf{g}}(u)=\lambda _{h,f,\mathbf{g}} $.
Then, up to multiplying it by a constant, $u$ is stricly positive and smooth, and it is a solution of
$$(E_{h,f,\mathbf{g}}):\triangle _{\mathbf{g}}u+hu=\lambda _{h,f,\mathbf{g}}.f.u^{\frac{n+2}{n-2}} $$
with or without the normalizing constant which can always be supressed just by multipliying again $u$ 
by a constant. Please, note that we will use these notations $\E$ and $\lambda _{h,f,\mathbf{g}}$ 
throughout all this article.

 Th. Aubin discoverded a very important relation between equation $\E$ and the notion of best constant 
 in the Sobolev imbedding theorems. Remember that the inclusion of $H_{1}^{2}(M)$ in $L^p(M)$ is 
 compact for $p<\frac{2n}{n-2}$ and only continuous for $p=\frac{2n}{n-2}$ which is called the 
 critical exponant for the Sobolev imbeddings and will be noted $2^{*}=\frac{2n}{n-2}$. The continuous 
 imbedding $H_{1}^{2}(M) \subset L^{2^{*}}(M)$ is expressed by the existence of two positive constants 
 $A$ and $B$ such that :
\begin{equation}
\forall u\in H_{1}^{2}(M):\,\left( \int_{M}\left| u\right| ^{\frac{2n}{n-2}}dv_{\mathbf{g}}\right) ^{%
\frac{n-2}{n}}\leq A\int_{M}\left| \nabla u\right| _{\mathbf{g}}^{2}dv_{%
\mathbf{g}}+B\int_{M}u{{}^{2}}dv_{\mathbf{g}}
\end{equation}
The best first constant is the minimum $A$ that one can put in (1) such that there exist $B$ with (1) still true. 
It was proved by E. Hebey and M. Vaugon \cite{H-V 1} that this minimum is attained, and its value is known to 
be the same as for the sharp euclidean Sobolev inequality,
\[
A_{min}=K(n,2){{}^{2}}=\frac{4}{n(n-2)\omega _{n}^{\frac{2}{n}}}
\]
where $\omega _{n}$ is the volume of the unit sphere of dimension $n$. One then take
 $B_{0}(\mathbf{g})$ to be the minimum $B$ such that (1) remains true with $A_{min}$; it is proved that $B_{0}(%
\mathbf{g})<+\infty $ \cite{H-V 1}. The inequality: $\forall u\in H_{1}^{2}(M)$
\begin{equation}
\left( \int_{M}\left| u\right| ^{\frac{2n}{n-2}}dv_{\mathbf{g}}\right) ^{%
\frac{n-2}{n}}\leq K(n,2){{}^{2}}\int_{M}\left| \nabla u\right| _{\mathbf{g}%
}^{2}dv_{g}+B_{0}(\mathbf{g})\int_{M}u{{}^{2}}dv_{\mathbf{g}} 
\end{equation}
is then sharp with respect to both the first and second constants, in the sense that none of them can be 
lowered. If the value of the best constant $A_{min}=K(n,2)^{2}$ is known and independent of the 
manifold $(M,\mathbf{g})$, on the other hand, $B_{0}(\mathbf{g})$, as the notation indicates, depends 
on the geometry and its study is difficult; it is for this purpose that "critical functions" were introduced 
by E.Hebey and M.Vaugon \cite{H-V 2}. When there shall be no risk of confusion, these constants will be 
denoted by $K$ et $B_{0}.$

As a remark, note that because of the compacity of the inclusion $H_{1}^{2}(M) \subset L^{p}(M)$ for 
$p<2^{*}$, standard variational methods and elliptic theory give rapidly existence of minimizing solutions 
of the equation $\triangle _{\mathbf{g}}u+hu=f.u^{p-1}$ when $\triangle _{\mathbf{g}}+h$ is a 
coercive operator. The case $p=2^{*}$ is therefore already a limit case. (Very little is known for $p>2^{*}$ 
without additional hypothesis, like e.g. invariance by symetry, see \cite{F1}.)

The best constants in the Sobolev embedding appeared in the study of equations $\E$ when Th. Aubin 
proved the following theorem:
\begin{theo}
[Aubin] For any riemannian manifold $(M,\g)$ of dimension $n\geqslant 3$, any function $h$ such that 
$\triangle _{\mathbf{g}}+h$ is a coercive operator, and any function $f$ such that 
$\underset{M}{Sup} f>0$, 
one always has 
$$\lambda _{h,f,\mathbf{g}}\leq \frac{1}{K(n,2){{}^{2}}(\stackunder{M}{Sup}f)^{\frac{n-2}{n}}}.$$
Furthermore, if this inequality is strict, then there exists a minimizing solution for $\E$.
\end{theo}
 
 \textbf{This theorem is the starting point of all this work. It proves the existence of minimizing solutions 
 to equation $\E$ under the hypothesis:}
 $$\lambda _{h,f,\mathbf{g}}< \frac{1}{K(n,2){{}^{2}}(\stackunder{M}{Sup}f)^{\frac{n-2}{n}}}.$$
\textbf{Our work is essentially concerned with the problem of the existence of minimizing solutions 
to these equations $\E$ in the "critical case" where}
$$\lambda _{h,f,\mathbf{g}}=\frac{1}{K(n,2){{}^{2}}(\stackunder{M}{Sup}f)^{\frac{n-2}{n}}},$$
\textbf{problem which is normally not solved by variational methods. It is for the study of this problem 
that we are now going to define the "critical functions".}

Let us first review the datas:

\textbf{Datas:} Throughout this article, $(M,\g)$ will be a compact riemannian manifold of dimension 
$n\geqslant 3$. We let $f:M\rightarrow \R$ be a fixed smooth function such that $\underset{M}{Sup} f>0$.
 Let also $h:M\rightarrow \R$ be a smooth function with the additional hypothesis that the operator 
 $\triangle _{\mathbf{g}}+h$ is coercive if $f$ is not positive on all of $M$. (Remember that continuity 
 of $h$ and $f$ is sufficient in the definitions and in most of the theorems. Also, if $f\leqslant0$ on $M$, 
 classical variational methods already give a lot of results for the existence of solutions; therefore $Sup f>0$ 
 is the most interesting case.)
 
 \begin{definition}
 With these datas, and with the above notations, we say that:
 \begin{itemize}
\item  $h$ is weakly critical for $f$ and $\mathbf{g}$ if $\lambda _{h,f,\mathbf{g}}=\frac{1}{K(n,2){{}^{2}}(%
\stackunder{M}{Sup}f)^{\frac{n-2}{n}}}$

\item  $h$ is subcritical for $f$ and $\mathbf{g}$ if $\lambda _{h,f,\mathbf{g}}<\frac{1}{K(n,2){{}^{2}}(\stackunder{M}{Sup}f)^{\frac{n-2}{n}}}$

\item  $h$ is \textbf{critical} for $f$ and $\mathbf{g}$ if $h$ 
is weakly critical and if for any function 
$k\leq h,\,k\neq h$ such that $\bigtriangleup _{\mathbf{g}}+k$ is coercive, $k$ is subcritical.
\end{itemize}
\end{definition}

Using the theorem of Th.Aubin, we can give an equivalent definition of critical functions. Indeed, using this 
theorem, it is easy to see that if $h$ is weakly critical and $\E$ has a minimizing solution $u$, then $h$ is a 
critical function; just note that for $k\leq h,\,k\neq h$, $I_{k,\g}(u)<I_{h,\g}(u)$. Therefore, we can give the following equivalent definition:
\begin{definition}
A function $h$ is critical for $f$ and $\g$ if:
\begin{itemize}
\item for any continuous function $k\leq h,\,k\neq h$ such that $\bigtriangleup _{\mathbf{g}}+k$ is coercive, 
(which is the case as soon as $k$ is close enough to $h$ in $C^{0}$), $(E_{k,f,\g})$ has a minimizing solution,
\item for any continuous function $k'\geq h,\,k'\neq h$, $(E_{k',f,\g})$ has \textbf{no} minimizing solution.
\end{itemize}
\end{definition}
Remark: if $h$ is weakly critical for a positive function $f$, necessarily, $\bigtriangleup _{\mathbf{g}}+h$ is coercive; 
just use the Sobolev inequality.

\textbf{Critical functions are thus introduced as "separating" functions giving rise to an equation having 
minimizing solutions, and functions giving rise to an equation that cannot have any such solution. We therefore 
have transformed the problem of the existence of minimizing solutions when $\lambda _{h,f,\mathbf{g}}=\frac{1}{K(n,2){{}^{2}}(%
\stackunder{M}{Sup}f)^{\frac{n-2}{n}}}$ to the problem of existence of minimizing solutions to $\E$ 
when $h$ is a critical function.}

Before passing to the theorems proved in this work, we have to give two very important properties of critical functions.

First, they transform in conformal changes of metric exactly like scalar curvature: indeed, let $u\in C^{\infty }(M),\,u>0$ and
 $\mathbf{g}^{\prime }=u^{\frac{4}{n-2}}\mathbf{g}$ a metric conformal to $\mathbf{g}$. Let also $h$ be 
 a smooth function. We set 
 $$h^{\prime }=\frac{\triangle _{\mathbf{g}}u+h.u}{u^{\frac{n+2}{n-2}}}.$$
 Then, some computations show that $h$ is critical for $f$ and $\g$ iff $h'$ is critical for $f$ and $\g'$.
 
 Second, we come back to the evaluation of $\lambda _{h,f,\mathbf{g}}$. Th. Aubin introduced, in the functional 
 $J_{h,f,\g}$ the following test functions:
 $$\psi _{k}(Q)=\left\lbrace 
\begin{array}{c}
(\frac{1}{k}+r{{}^{2}})^{-\frac{n-2}{2}}-(\frac{1}{k}+\delta {{}^{2}})^{-%
\frac{n-2}{2}}\,\,\,\,if\,r<\delta \\ 
0\,\,\,\,\,\,\,\,if\,\,\,r\geq \delta
\end{array} \right. $$
where: $\delta <injM$ (the injectivity radius of $M$), $P\,\in M$ is a fixed point, 
$k\in \Bbb{N}^{*}$, and where $r=d_{\mathbf{g}}(P,Q).$ When $dimM=n\geqslant 4$, we get, if $P$ is a point 
where $f$ is maximum on $M$:
 \begin{center}
$J_{h,f,\mathbf{g}}(\psi _{k})=\frac{1}{K(n,2){{}^{2}}(\stackunder{M}{Sup}%
f)^{\frac{n-2}{n}}}\left\{ 1+\frac{1}{n(n-4)}\left( \frac{4(n-1)}{n-2}%
h(P)-S_{\mathbf{g}}(P)+\frac{n-4}{2}\frac{\bigtriangleup _{\mathbf{g}}f(P)}{%
f(P)}\right) \frac{1}{k}\right\} +o(\frac{1}{k})$
\end{center}
 We therefore get the following important proposition:
 \begin{proposition}
 If $dimM\geq4$ and if $h$ is weakly critical for $f$ and $\g$ (thus in particular if it is critical), as $\lambda _{h,f,\mathbf{g}}=\frac{1}{K(n,2){{}^{2}}(\stackunder{M}{Sup}f)^{%
\frac{n-2}{n}}}$, necessarily, if $P$ is a point of maximum of $f$:
\[
\frac{4(n-1)}{n-2}h(P)\geq S_{\mathbf{g}}(P)-\frac{n-4}{2}\frac{\bigtriangleup _{\mathbf{g}}f(P)}{f(P)} 
\]
 \end{proposition}
 Remark: if $f$ is constant on $M$, this means that $\frac{4(n-1)}{n-2}h\geq S_{\mathbf{g}}$ on all of $M$. Note 
 also that in dimension 4, the term $\frac{\bigtriangleup _{\mathbf{g}}f(P)}{f(P)}$ disappears.
 
\section{Statement of the results}
In all what follows, we will make the following hypothesis:

\textbf{Hypothesis (H):} We now suppose that $dimM=n\geqslant 4$. We suppose that all our functions $h$ are such that 
$\bigtriangleup _{\mathbf{g}}+h$ is coercive. Also, $f$ will always be a smooth function such that $\underset{M}{Sup} f>0$. 
We will denote $Max\,f=\{x\in M/ f(x)=\underset{M}{Sup} f \}$.

Our first theorem concerns the existence of minimizing solutions to $\E$ when $h$ is critical.
\begin{theo}
If $h$ is a critical function for $f$ and $\g$, ($h,f,\g$ verifying \textbf{H}), and if for all point $P$ where $f$ is 
maximum on $M$, we have $$\frac{4(n-1)}{n-2}h(P)>S_{\mathbf{g}}(P)-\frac{n-4%
}{2}\frac{\bigtriangleup _{\mathbf{g}}f(P)}{f(P)},$$ then there exist a minimizing solution for $\E$.
\end{theo}
This theorem is an immediate consequence of the following result, more general but more technical in its statement. (Just 
take $h_{t}=h-t$ to get the theorem above.)
\begin{theorem}
Let $h$ be a weakly critical function for $f$ and $\g$, (assuming hypothesis \textbf{H}). If, for all point $P$ where $f$ 
is maximum, we have
 $$\frac{4(n-1)}{n-2}h(P)>S_{\mathbf{g}}(P)-\frac{n-4}{2}\frac{\bigtriangleup _{\mathbf{g}}f(P)}{f(P)},$$
and if there exists a family of functions $(h_{t})$, $h_{t}\lvertneqq  h$, $h_{t}$ being sub-critical for all 
$t$ in a neighbourhood of a real $t_{0}\in \Bbb{R}$, and such that $h_{t}\stackunder{t\rightarrow t_{0}}{\rightarrow }h$
 in $C^{0,\alpha }$, then there exists a minimizing solution for $\E$, and therefore, $h$ is critical for $f$ and $\g$.
\end{theorem}
E. Hebey and M. Vaugon, in the context of their study of $B_{0}(\g)$, proved this theorem in the case where $f$ is constant, 
and as them, we base our computations on the article of Djadli and Druet \cite{D-D}. 
The presence of a non-constant function $f$ on the right of equation $\E$ introduces new difficulties in the proof, 
and requires the use of very powerfull estimates concerning concentration phenomena's, called $C^0-theory$, 
due to Druet and Robert \cite{D-R}, available in \cite{DHR}; 
the use of $C^0-theory$ was kindly suggested to us by E. Hebey. 
Also, an alternate proof, not using $C^0-theory$, thus in some sense more elementary, 
but requiring the additional hypothesis that the hessian of $f$ is non-degenerate at its points of maximum on $M$, 
will, as a "byproduct", prove another very important estimate concerning these concentration phenomena's, not 
available without heavy hypothesis in the case when $f$ is a constant function; this estimate concerns the speed of 
convergence to a concentration point, (see subsection 4.2), is of independent interest, and was obtained in the
author's PHD thesis \cite{C3} to prove theorem 1.

The next natural question is of course to know if there exist critical functions. The answer, positive, will appear to be a 
consequence of theorem 1. We will say that a set $E\subset M$ is \textit{thin} if $M-E$ contains a dense open subset.
\begin{theorem}
Being given the manifold $(M,\g)$ and a non constant function $f$, there exist infinitely many 
functions $h$ critical for $f$ and $\g$, which satisfy, in each point $P$ of maximum of $f$,
 $$\frac{4(n-1)}{n-2}h(P)>S_{\mathbf{g}}(P)-\frac{n-4}{2}\frac{\bigtriangleup _{\mathbf{g}}f(P)}{f(P)}\,\,\,\,(*)$$
By theorem 1, these critical functions are such that $\E$ have minimizing solutions. 
Also, if the set of maximum points of $f$ is thin and if $\int_{M}f>0$, there exist strictly positive such critical 
functions $h$, i.e. satisfying (*).
\end{theorem}

These first theorems lead us to modify slightly our vision of critical functions. Note that in equation $\E$, there are 
three datas that one can modify: the functions $h$ and $f$, of course, but also the metric $\g$ in a conformal 
class, as, by the conformal laplacian transformation formula, the equation is changed in a similar one if we change 
$\g$ in $\g'=u^{\frac{4}{n-2}}.\g$. This lead us to the following definition:
\begin{center}
\textit{$(h,f,\g)$ is a critical triple if $h$ is a critical function for $f$ and $\g$.}
\end{center}
We shall say that the triple $(h,f,\g)$ has minimizing solutions if $\E$ has; we can also speak of weakly critical 
or sub-critical triples.We then asked ourselves the following question:

\textit{Being given two of the three datas of a triple, can one find the third to obtain a critical triple?}

For example, the problem of the existence of  critical functions can be formulated in the following manner: we 
are given the function $f$ and the metric $\g$, can we complete the triple $(.,f,\g)$ by a function $h$ to obtain a 
critical triple $(h,f,\g)$?

We adress the two other questions, first fixing $h$ and $f$ and seeking a conformal metric $\g'$, and then fixing 
the function $h$ and the metric $\g$ and seeking a function $f$. We obtain answers expressed by the following two 
theorems:
\begin{theorem}
On the manifold $(M,\g)$, let be given a function $h$ and a function $f$, satisfying (\textbf{H}). 
We suppose that the set of maximum points of $f$ is thin.
Then, there exist a metric $\g'$ conformal to $\g$ such that $(h,f,\g')$ is a critical triple. 
Moreover, we can find $\g'$ such that $(h,f,\g')$ has minimizing solutions.
\end{theorem}

This theorem was proved by E. Humbert and M. Vaugon in the case $f=cst=1$ and $M$ not conformally diffeomorphic 
to the sphere, \cite{Hu-V}. Their method works in the case 
of a non constant function $f$ and an arbitrary manifold once it is proved that we can suppose the existence of positive critical functions 
satisfying the strict inequality (*) in theorem 2, result we included in this theorem (note that, as $Sup f>0$, we can 
always find a metric $\textbf{g'}$ conformal to $\g$ such that $\int_{M}fdv_{\g'}>0$). 
 In fact, when $M$ is not conformally diffeomorphic to the sphere and $S_{\g}$ is constant, it can be 
proved that $B_{0}(\g)K(n,2)^{-2}$ is a critical (constant) function for $1$ and $\g$, and it is obviously 
positive. We will discuss weaker hypothesis for this theorem, as well as the problem of existence of positive 
critical functions in section 6.

The last question brings us to the following answer when the dimension of $M$ is greater than 5, requirement which is 
linked to the fact that $\frac{\bigtriangleup _{\mathbf{g}}f(P)}{f(P)}$ dissapears in dimension 4 in the 
inequality of Proposition 1.
\begin{theorem}
Let be given the manifold $(M,\g)$ of dimension $n\geq5$, and a function $h$ such that $\bigtriangleup _{\mathbf{g}}+h$ 
is coercive. Then, there exists a non constant function $f$ such that $(h,f,\g)$ is critical with minimizing solutions if, and 
only if, $(h,1,\g)$ is a sub-critical triple (where 1 is the constant function 1).
\end{theorem}

Note that if $(h,1,\g)$ is weakly critical, then either this triple has minimizing solutions in which case it is a critical 
triple, or there is no non-constant function $f$ such that $(h,f,\g)$ is critical with minimizing solutions (see the proof 
and what follows). 
The proof of this theorem is quite difficult, and make use of the method developped for the proof of theorem 1. 
Also, this proof brought us to make some more remarks about critical functions. First, it is easily seen, by using the 
functional $J$, that if $(h,f,\g)$ is a critical triple, then, for any constant $c>0$, $(h,c.f,\g)$ is also a critical triple. 
It would therefore be more appropriate to speak of triple $(h,[f],\g)$ where $[f]=\{c.f\,/c>0\}$ could be called 
the "class" of $f$. Note for example that we can always suppose that $Sup f=1$; also, to compare two triples $(h,f,\g)$ 
and $(h,f',\g)$, one has to suppose that $Sup f=Sup f'$. Note also that on $[f]$, the quotient  $\frac{\bigtriangleup _{\mathbf{g}}f}{f}$ 
is constant. Second, in the proof of theorem 4, we had to approximate the function $f$ by a family $(f_{t})$, unlike 
theorem 1 where we used a family $(h_{t})$ approaching $h$. This suggested another possible definition of critical 
functions, dual to the first one in the sense that we exchange the role of $h$ and $f$.
\begin{definition}
Let $(M,\g)$ be of dimension $n\geq3$ and $h$ be such that $\bigtriangleup _{\mathbf{g}}+h$ is coercive. 
We shall say that a smooth function $f$ such that $\underset{M}{Sup} f>0$ is critical for $h$ and $\g$ if:
\begin{itemize}
\item  a/: $\lambda _{h,f,\mathbf{g}}=\frac{1}{K(n,2){{}^{2}}(\stackunder{M}{Sup}f)^{\frac{n-2}{n}}}$

\item  b/: for any smooth function $f^{\prime }$ such that $Supf=Supf^{\prime }$ and $f^{\prime }\gneqq f$, \\

$\lambda_{h,f^{\prime },\mathbf{g}}<\frac{1}{K(n,2){{}^{2}}(\stackunder{M}{Sup}%
f^{\prime })^{\frac{n-2}{n}}}$

\item  Remark: if $Supf=Supf^{\prime }$ and $f^{\prime
}\lvertneqq f$, then $\lambda _{h,f^{\prime },\mathbf{g}}=\frac{1}{%
K(n,2){{}^{2}}(\stackunder{M}{Sup}f^{\prime })^{\frac{n-2}{n}}}$ as $J_{h,f^{\prime },\mathbf{g}}(w)\geqslant J_{h,f,\mathbf{g}}(w)$%
for any function $w$.
\end{itemize}
\end{definition}
It is then natural to ask if the two definitions are equivalent ($\g$ being fixed):
\begin{center}
\textit{Is $f$ critical for $h$ if, and only if, $h$ is critical for $f$ ?}
\end{center} 

This question seems quite difficult. A positive answer would justify the concept of critical triple. Remember that, 
because of proposition 1, we have in both cases, when $P$ is a point where $f$ is maximum on $M$: 
$$\frac{4(n-1)}{n-2}h(P)\geqslant S_{\mathbf{g}}(P)-\frac{n-4}{2}\frac{\bigtriangleup _{\mathbf{g}}f(P)}{f(P)}.$$ 

We obtain the following theorem:
\begin{theorem}
Let $(M,\g)$ be a compact manifold of dimension $n\geq5$, and let $h$ be a function such that  $\bigtriangleup _{\mathbf{g}}+h$ is coercive. 
Let $f$ be a smooth function such that $\underset{M}{Sup} f>0$. We suppose 
that for any point $P$ where $f$ is maximum on $M$:
$$\frac{4(n-1)}{n-2}h(P)> S_{\mathbf{g}}(P)-\frac{n-4}{2}\frac{\bigtriangleup _{\mathbf{g}}f(P)}{f(P)}.$$
Then, $f$ is critical for $h$ if, and only if, $h$ is critical for $f$.
\end{theorem}
Remark: if $1$ is critical for $h$, then every non constant function $f$, such that $Sup f=1$, is weakly critical for $h$ 
with \textit{no} minimizing solutions. Indeed, here again if a function $f$ is weakly critical for $h$ with a minimizing 
solution, then $f$ is critical.

There is an interesting consequence of theorems 4 and 5. We said in the introduction that an important application 
of equations $\E$ was the study of prescribed scalar curvature: being given a smooth function $f$ on the manifold 
$(M,\g)$, is $f$ the scalar curvature of a metric conformal to $\g$? The theorem of Th. Aubin shows that if $f$ is 
sub-critical for $S_{\g}$, then $f$ is a scalar curvature. Theorem 4 applied to $h=\frac{n-2}{4(n-1)}S_{\g}$ 
shows that:

\textit{On a compact manifold $(M,\g)$ not conformaly diffeomorphic to the sphere, there exist scalar curvatures 
of metric conformal to $\g$ that are only weakly critical, (more precisely critical).}
\\

Another application, remarked by E. Hebey, is the study of \textit{Sobolev inequality in the presence of a twist.} 
See for more details on the construction of twisted metrics the article \cite{C-H-V}.

The previous theorems all deal with manifolds of dimension at least 4, or even 5. We will give results concerning 
the dimension 3 in the last section. They are very interesting, but they are rapid generalisations of results obtained 
by O. Druet in the case $f=constant$ \cite{D}, the introduction of a non constant $f$ introducing this time no real difficulties. 
We prefer therefore to state them at the end, with no proof, sending the reader to the article of O. Druet or to our 
PHD thesis \cite{C3}, available online, for more details.

\section{The three main tools}
We want to present here the three main tools used in the proof of our various theorems. These tools were developed by 
several persons since M. Vaugon and P.L. Lions, essentially E. Hebey, O. Druet F. Robert, 
M. Struwe, E. Humbert and Z. Faget, among others. 
\subsection{The concentration point.}
To prove the existence of a solution $u>0$ to our equation 
\[
(E_{h,f,\mathbf{g}}):\,\triangle _{\mathbf{g}}u+h.u=\lambda .f.u^{\frac{n+2}{%
n-2}}, 
\]
the idea will often be to associate a family of equations having minimizing solutions $u_{t}>0$ : 
\[
E_{t}:\,\,\triangle _{\mathbf{g}}u_{t}+h_{t}.u_{t}=\lambda _{t}.f.u_{t}^{%
\frac{n+2}{n-2}} 
\]
with
\[
h_{t}\rightarrow h\,\,\,\,in \,\,\,\,C^{0,\alpha }(M) 
\]
and $\lambda _{t}\rightarrow \lambda $ a converging sequence of real numbers, in such a way that for some $u\in H_{1}^{2}$ : 
$u_{t}\rightarrow u\,$ strongly in $L^{p}$ , 
$p<2^{*}$, and $u_{t}\rightharpoondown u$ weakly in  $H_{1}^{2}$ with a constraint 
\[
\int_{M}f.u_{t}^{2^{*}}dv_{\mathbf{g}}=1. 
\]
To simplify, we will suppose that all convergences are for $t\rightarrow t_{0}=1$. The difficulty will be to prove 
that $u$ is not the trivial zero solution, as then, by the maximum principle, we have $u>0$. 
We will proceed by contradiction, and suppose $u\equiv 0$. The idea is then that, because of the condition $\int_{M}f.u_{t}^{2^{*}}=1,$ 
all the "mass" of the functions $u_{t}$, which converge to  0 in $L^{p}$ , $p<2^{*}$, concentrates around a point 
of the manifold. We thus define:
\begin{definition}
\textit{\ }$x_{0}\in M$ is a point of concentration of the sequence $(u_{t})$ if for any $\delta >0$ :
\[
\stackunder{t\rightarrow t_{0}}{\lim \sup }\int_{B(x_{0},\delta
)}u_{t}^{2^{*}}dv_{\mathbf{g}}>0 
\]
\end{definition}
It is easy to see that because $M$ is compact and we require  $\int_{M}f.u_{t}^{2^{*}}dv_{\mathbf{g}}=1$, there exist 
at least one point of concentration. We will show that there exists only one point of concentration, that it is a point where 
$f$ is maximum, and that there exist a sequence of points $x_{t}$ converging to a point $x_{0}\in M$ such that
\[
u_{t}(x_{t})=\stackunder{M}{\max }u_{t}\rightarrow +\infty , 
\]
and 
\[
u_{t}\rightarrow 0\,\,\,in\,\,\,C_{loc}^{0}(M-\{x_{0}\}). 
\]
In fact the idea is that one can do "as if" the functions $u_{t}$ have compact support in a small neigbourhood of $x_{0}$ 
when $t$ is close to $t_{0}$.

\subsection{Blow-up analysis}
Thanks to the concentration point, one brings back the study of the family $u_{t}$ converging to 0, to what happens around 
$x_{0}$. The idea of \textit{blow-up analysis} is to do a "change of scale" around $x_{0}$: we will call \textit{blow-up} 
of center $x_{t}$ and coefficient $k_{t}$ the following sequence of charts and changes of metrics. We consider, for $\delta$ 
small enough:
\[
\begin{array}{ccccc}
B(x_{t},\delta ) & \stackrel{\exp _{x_{t}}^{-1}}{\rightarrow } & \stackunder{%
}{B(0,\delta )\subset \Bbb{R}^{n}} & \stackrel{\psi _{k_{t}}}{\stackunder{}{%
\rightarrow }} & B(0,k_{t}\delta )\subset \Bbb{R}^{n} \\ 
&  & x\,\,\, & \mapsto \,\, & k_{t}x \\ 
\,\mathbf{g}\, & \rightarrow & \,\,\mathbf{g}\,_{t}=\exp _{x_{t}}^{*}\mathbf{%
g}\, & \rightarrow & \,\widetilde{\,\mathbf{g}\,}_{t}=k_{t}^{2}(\psi
_{k_{t}}^{-1})^{*}\mathbf{g}_{t}
\end{array}
\]
where $\exp _{x_{t}}^{-1}$ is the chart deduced from the exponential map in $x_{t}$. We set 
\[
\overline{u}_{t}=u_{t}\circ \exp _{x_{t}}\, 
\]
\[
\overline{f_{t}}=f\circ \exp _{x_{t}} 
\]
\[
\,\overline{h_{t}}=h_{t}\circ \exp _{x_{t}} 
\]
We have
\begin{eqnarray*}
\triangle _{\mathbf{g}_{t}}\overline{u}_{t}+\overline{h}_{t}.\overline{u}%
_{t} &=&\lambda _{t}\overline{f_{t}}.\overline{u}_{t}^{\frac{n+2}{n-2}} \\
\int_{B(0,r)}\overline{u}_{t}^{\alpha }dv_{\,\mathbf{g}\,_{t}}
&=&\int_{B(x_{t},r)}u_{t}^{\alpha }dv_{\mathbf{g}}\text{ for all }\alpha
\geq 1
\end{eqnarray*}
We then set
\begin{eqnarray*}
m_{t} &=&\stackunder{M}{Max}\,u_{t}\, \\
\,\widetilde{u}_{t} &=&m_{t}^{-1}\overline{u}_{t}\circ \psi _{k_{t}}^{-1}\,
\\
\,\widetilde{h}_{t} &=&h_{t}\circ \psi _{k_{t}}^{-1}\, \\
\,\widetilde{f}_{t} &=&\overline{f}_{t}\circ \psi _{k_{t}}^{-1}\, \\
\,\widetilde{\,\mathbf{g}\,}_{t} &=&k_{t}^{2}(\exp _{x_{t}}\circ \psi
_{k_{t}}^{-1})^{*}\mathbf{g},
\end{eqnarray*}
so in particular
\begin{eqnarray*}
\,\,\widetilde{u}_{t}(x) &=&m_{t}^{-1}\overline{u}_{t}(\frac{x}{k_{t}})\, \\
\,\widetilde{\,\mathbf{g}\,}_{t}(x) &=&\exp _{x_{t}}^{*}\mathbf{g}(\frac{x}{%
k_{t}}).
\end{eqnarray*}
Then: 
\begin{eqnarray}
(\widetilde{E}_{t})\, &:&\,\triangle _{\widetilde{\mathbf{g}}_{t}}\widetilde{%
u}_{t}+\frac{1}{k_{t}^{2}}\widetilde{h}_{t}.\widetilde{u}_{t}=\frac{m_{t}^{%
\frac{4}{n-2}}}{k_{t}^{2}}\lambda _{t}\widetilde{f}_{t}.\widetilde{u}_{t}^{%
\frac{n+2}{n-2}}   \\
and\, &:&\,\int_{B(0,k_{t}r)}\widetilde{u}_{t}^{\alpha }dv_{\widetilde{\,%
\mathbf{g}\,}_{t}}=\frac{k_{t}^{n}}{m_{t}^{\alpha }}\int_{B(x_{t},r)}u_{t}^{%
\alpha }dv_{\mathbf{g}}  \nonumber
\end{eqnarray}
We will mostly use the following parameters : we consider a sequence of points ($x_{t})$ such that: 
\[
m_{t}=\stackunder{M}{Max}\,u_{t}=u_{t}(x_{t}):=\mu _{t}^{-\frac{n-2}{2}}
\]
and
\[
\,k_{t}=\mu _{t}^{-1}.
\]
$\mu _{t}$ will appear to be a fundamental parameter in the study of concentration phenomena's. Noting $(x^{i})$ the coordinates in $\Bbb{%
R}^{n}$, one has : 
\begin{eqnarray}
(\widetilde{E}_{t})\, &:&\,\triangle _{\widetilde{\mathbf{g}}_{t}}\widetilde{%
u}_{t}+\mu _{t}^{2}.\widetilde{h}_{t}.\widetilde{u}_{t}=\lambda _{t}%
\widetilde{f}_{t}.\widetilde{u}_{t}^{\frac{n+2}{n-2}}   \\
and\, &:&\,\int_{B(0,\mu _{t}^{-1}r)}x^{i_{1}}...x^{i_{p}}.\widetilde{u}%
_{t}^{\alpha }dv_{\widetilde{\,\mathbf{g}\,}_{t}}=\mu _{t}^{-p-n+\alpha 
\frac{n-2}{2}}\int_{B(0,r)}x^{i_{1}}...x^{i_{p}}\overline{u}_{t}^{\alpha
}dv_{\,\mathbf{g}\,_{t}}  \nonumber
\end{eqnarray}
A very important result is that when $\mu _{t}\rightarrow 0$ and therefore $k_{t}\rightarrow
+\infty $, the components of $\,\widetilde{\mathbf{g}}_{t}$ converge in 
$C_{loc}^{2}$ to those of the euclidean metric, and $(\widetilde{E}_{t})$``converges'' to the equation: 
\[
\triangle _{e}\widetilde{u}=\lambda f(x_{0}).\widetilde{u}^{\frac{n+2}{n-2}}
\]
in the sense that
\[
\widetilde{u}_{t}\rightarrow \widetilde{u}\,\,\,\,in \,\,\,\,C_{loc}^{2}(\Bbb{R}^{n}).
\]
It is known, then, that
\[
\widetilde{u}=(1+\frac{\lambda f(x_{0})}{n(n-2)}\left| x\right| ^{2})^{-%
\frac{n-2}{2}}.
\]

\subsection{The iteration process}
The idea of the Möser iteration process is to multiply the equations $(E_{t})$ by succesive powers $u_{t}^k$ of the functions 
$u_{t}$ and to integrate over $M$ to obtain bounds on increasing $L^p$-norms of the $u_{t}$. To localize the study around 
the concentration point $x_{0}$, which is a maximum point for $f$, we shall in fact multiply the equations by $\eta {{}^{2}}u_{t}^{k}$ where $\eta $ is a 
cut-off function equal to 1 (resp.0) on a ball $B(x_{0},r)$ where $f\geq 0$, and equal to 0
(resp. 1) on $M\backslash B(x_{0},2r)$, and where $k\geq 1$, then integrate by part. We will therefore be able to study 
blow-up around $x_{0}$ using this method. We get after some integrations by parts, and using equation $(E_{t})$ :
\begin{equation}
\frac{4k}{(k+1)^{2}}\int_{M}\left| \nabla (\eta u_{t}^{\frac{k+1}{2}%
})\right| ^{2}=\lambda _{t}\int_{M}f\eta ^{2}u_{t}^{\frac{n+2}{n-2}%
}u_{t}^{k}+\int_{M}(\frac{2}{k+1}\left| \nabla \eta \right| ^{2}+\frac{2(k-1)%
}{(k+1)^{2}}\eta \triangle \eta -\eta ^{2}h_{t})u_{t}^{k+1}
\end{equation}
where the integrals are taken with the measure $dv_{\mathbf{g}}$. Then using Hölder inequality, 
if $f\geq 0$ on $Supp\,\eta $ we obtain: 
\[
\lambda _{t}\int_{M}f\eta ^{2}u_{t}^{\frac{n+2}{n-2}}u_{t}^{k}\leq \lambda
_{t}(\stackunder{Supp\,\eta }{Sup}f)^{\frac{n-2}{n}}.(\int_{Supp\,\eta
}fu_{t}^{\frac{2n}{n-2}})^{\frac{2}{n}}.(\int_{M}(\eta u_{t}^{\frac{k+1}{2}%
})^{\frac{2n}{n-2}})^{\frac{n-2}{n}} 
\]
Then using Sobolev inequality : 
\[
(\int_{M}(\eta u_{t}^{\frac{k+1}{2}})^{\frac{2n}{n-2}})^{\frac{n-2}{n}}\leq
K(n,2){{}^{2}}\int_{M}\left| \nabla (\eta u_{t}^{\frac{k+1}{2}})\right|
^{2}+B\int_{M}\eta u_{t}^{k+1} 
\]
with $B>0$. Therefore:

\smallskip

\begin{eqnarray*}
\frac{4k}{(k+1)^{2}}(\int_{M}(\eta u_{t}^{\frac{k+1}{2}})^{\frac{2n}{n-2}})^{%
\frac{n-2}{n}} \leq &\lambda _{t}K(n,2){{}^{2}}(\stackunder{Supp\,\eta }{Sup%
}f)^{\frac{n-2}{n}}.(\int_{Supp\,\eta }fu_{t}^{\frac{2n}{n-2}})^{\frac{2}{n}%
}.(\int_{M}(\eta u_{t}^{\frac{k+1}{2}})^{\frac{2n}{n-2}})^{\frac{n-2}{n}} \\
&+\int_{M}(\frac{4k}{(k+1){{}^{2}}}B\eta +\frac{2}{k+1}\left| \nabla \eta
\right| ^{2}+\frac{2(k-1)}{(k+1)^{2}}\eta \triangle \eta -\eta
^{2}h_{t})u_{t}^{k+1}
\end{eqnarray*}
Then:
\begin{equation}
Q(t,k,\eta ).(\int_{M}(\eta u_{t}^{\frac{k+1}{2}})^{\frac{2n}{n-2}})^{\frac{%
n-2}{n}}\leq (\frac{4k}{(k+1){{}^{2}}}B+C_{0}+C_{\eta })\int_{Supp\,\eta
}u_{t}^{k+1}\,\,\,\,\,\,\,\,\,\,\,\,\,\,\,\,\,\,\,\,\,\,\,\,\,\,\,\,\,\,\,\,%
\,\,\,\,\,\,\,\,\,\,\,\,\,\,\,\,\,\,\,\,\,\,\,  
\end{equation}
where
\[
Q(t,k,\eta )=\frac{4k}{(k+1){{}^{2}}}-\lambda _{t}K(n,2){{}^{2}}(\stackunder{%
Supp\,\eta }{Sup}f)^{\frac{n-2}{n}}.(\int_{Supp\,\eta }f.u_{t}^{2^{*}})^{%
\frac{2}{n}} 
\]
where we remind that $2^{*}=\frac{2n}{n-2}$ and where $C_{0}\,et\,C_{\eta }$ are constants independant of
 $k$ and $t$ and such that $\forall k\geq 1,\forall t:$
\[
\,\,\left\| \frac{2}{k+1}\left| \nabla \eta \right| ^{2}+\frac{2(k-1)}{%
(k+1)^{2}}\eta \triangle \eta \right\| _{L^{\infty }(M)}\leq \,C_{\eta
}\,\,and\,\,\left\| h_{t}\right\| _{L^{\infty }(M)}\leq C_{0}\,. 
\]
If the sign of $f$ changes on $Supp\,\eta $, we go back to Hölder's inequality:
\[
\lambda _{t}\int_{M}f\eta ^{2}u_{t}^{\frac{n+2}{n-2}}u_{t}^{k}\leq \lambda
_{t}(\stackunder{Supp\,\eta }{Sup}\left| f\right| ).(\int_{Supp\,\eta
}u_{t}^{\frac{2n}{n-2}})^{\frac{2}{n}}.(\int_{M}(\eta u_{t}^{\frac{k+1}{2}%
})^{\frac{2n}{n-2}})^{\frac{n-2}{n}} 
\]
to obtain (6) with: 
\begin{equation}
Q(t,k,\eta )=\frac{4k}{(k+1){{}^{2}}}-\lambda _{t}K(n,2){{}^{2}}(\stackunder{%
Supp\,\eta }{Sup}\left| f\right| ).(\int_{Supp\,\eta }u_{t}^{2^{*}})^{\frac{2%
}{n}}  
\end{equation}
One can also replace $\stackunder{Supp\,\eta }{Sup}%
\left| f\right| $ by $\stackunder{M}{Sup}f$.

The goal is to show that $(\eta u_{t})$ is bounded in $L^{\frac{k+1}{2}2^{*}}$ and therefore that we can 
extract a sub-sequence converging strongly in $L^{2^{*}}$.

\textbf{Remark}

Those three tools also work for more general equations that we can associate to $(E_{h,f,\mathbf{g}}):\,\triangle _{\mathbf{g}}u+h.u=\mu
_{h}.f.u^{\frac{n+2}{n-2}}$. like e.g.
\[
E_{t}:\,\,\triangle _{\mathbf{g}}u_{t}+h_{t}.u_{t}=\lambda
_{t}.f_{t}.u_{t}^{q_{t}-1} 
\]
where $q_{t}\rightarrow 2^{*}$ and $f_{t}\rightarrow f$ in some $L^{p}$, still with $h_{t}\rightarrow h$\thinspace in $%
\,C^{0,\alpha }(M)$ and $\lambda _{t}\rightarrow \lambda $, and where we require 
$\int_{M}f_{t}u_{t}^{q_{t}}dv_{\mathbf{g}}=1$.

$x_{0}\in M$ will then be a concentration point for $(u_{t})$\ if for any $\delta >0$: 
\[
\stackunder{t\rightarrow t_{0}}{\lim \sup }\int_{B(x_{0},\delta
)}u_{t}^{q_{t}}>0 
\]
Blow-up and iteration process are then similar. Formulas (6) et (7) become
\[
Q(t,k,\eta ).(\int_{M}(\eta u_{t}^{\frac{k+1}{2}})^{q_{t}})^{\frac{2}{q_{t}}%
}\leqslant C\int_{Supp\,\eta
}u_{t}^{k+1}
\]
where
\[
Q(t,k,\eta )=\frac{4k}{(k+1){{}^{2}}}-\lambda
_{t}(Vol_{g}(M))^{\frac{q_{t}}{2^{*}}-1}K{{}(n,2)^{2}}(\stackunder{%
Supp\,\eta }{Sup}\left| f\right| ).(\int_{Supp\,\eta }u_{t}^{q_{t}})^{\frac{%
q_{t}-2}{q_{t}}} 
\]
 \subsection{Principle of the proof of theorem 1}
 We want to prove the existence of a positive solution $u$ to the equation:
\[
(E_{h,f,\mathbf{g}}):\,\triangle _{\mathbf{g}}u+h.u=\lambda .f.u^{\frac{n+2}{%
n-2}} 
\]
As we said at the beginning of this section, we will study associated equations and minimizing solutions:
\[
E_{t}:\,\,\triangle _{\mathbf{g}}u_{t}+h_{t}.u_{t}=\lambda _{t}.f.u_{t}^{%
\frac{n+2}{n-2}} 
\]
with
\[
h_{t}\stackrel{\lvertneqq }{\rightarrow }h\,\,\,\,in\,\,\,\,C^{0,\alpha
}(M) 
\]
the $h_{t}$ being sub-critical by hypothesis. The idea is the following: we are going to introduce in the Sobolev inequality 
the equation $E_{t}$ to obtain a contradiction, when $u_{t}\rightharpoondown 0$, with the condition 
\[
\frac{4(n-1)}{n-2}h(P)>S_{\mathbf{g}}(P)-\frac{n-4}{2}\frac{\bigtriangleup _{%
\mathbf{g}}f(P)}{f(P)}
\]
when $P$ is a point of maximum of $f$; remember we will show that the concentration point is a maximum for $f$.

To simplify, let us assume that $\mathbf{g}$ is flat around $x_{0}$ and that $f\equiv 1$. 
Then $S_{\mathbf{g}}=0$ near $x_{0}$
and our hypothesis is:
\[
h(x_{0})>0=\frac{n-2}{4(n-1)}S_{\mathbf{g}}(x_{0}) 
\]
Therefore, for $t$ close to 1: $h_{t}(x_{0})>0$. But on the one hand, as $u_{t}$ is minimizing, 
\[
\lambda _{h_{t},f,\mathbf{g}=}J_{h_{t}}(u_{t}):=\lambda _{t}<K(n,2)^{-2} 
\]
and thus: 
\begin{equation}
\int_{M}\left| \nabla u_{t}\right| ^{2}dv_{\mathbf{g}}+\int_{M}h_{t}.u{%
{}_{t}^{2}}dv_{\mathbf{g}}=\lambda _{t}(\int_{M}u_{t}^{2^{*}}dv_{\mathbf{g}%
})^{\frac{2}{2^{*}}}=\lambda _{t}<K(n,2)^{-2}  
\end{equation}
because $\int_{M}u_{t}^{2^{*}}dv_{\mathbf{g}}=1$; and on the other hand, the Sobolev euclidean inequality gives
\begin{equation}
K(n,2)^{-2}(\int_{\Bbb{R}^{n}}v^{2^{*}})^{\frac{2}{2^{*}}}\leq \int_{\Bbb{R}%
^{n}}\left| \nabla v\right| ^{2}\,.  
\end{equation}
However, if $u_{t}\rightarrow 0$, we will show that there is a concentration phenomena, and, as we said, this 
enables us to do "as if" the functions $u_{t}$ had compact support in a small neighbourhood $B$ of $x_{0}$ where $h_{t}>0.$ We 
would then have because of (8) 
\[
\int_{B}\left| \nabla u_{t}\right| ^{2}<K(n,2)^{-2}(\int_{B}u_{t}^{2^{*}})^{%
\frac{2}{2^{*}}}=K(n,2)^{-2} 
\]
as $h_{t}>0=S_{\mathbf{g}}$ in $B$; and on the other side, because of (9) 
\[
\int_{B}\left| \nabla u_{t}\right| ^{2}\geq
K(n,2)^{-2}(\int_{B}u_{t}^{2^{*}})^{\frac{2}{2^{*}}}=K(n,2)^{-2} 
\]
thus a contradiction. To apply this idea, we will have to multiply the functions $u_{t}$ by cut-off functions, make developments 
of the metric and of $f$, and apply all the results concerning concentration phenomena's that we shall expose in the next section.
\\
Let us go in some more details.
We want to prove the existence of a minimizing solutions to :
\[
(E_{h,f,\mathbf{g}}):\,\triangle _{\mathbf{g}}u+h.u=\lambda .f.u^{\frac{n+2}{%
n-2}}\,\,with\,\,\int_{M}fu^{2^{*}}=1 
\]
when there exists a family $(h_{t})$ of subcritical functions:
\[
h_{t}\stackrel{\lvertneqq }{\rightarrow }h\,\,\,\,in\,\,\,\,C^{0,\alpha
}(M) 
\]
As the $(h_{t})$ are subcritical, there exist a family $u_{t}$ of minimizing solutions of the equations
\[
E_{t}:\,\,\triangle _{\mathbf{g}}u_{t}+h_{t}.u_{t}=\lambda _{t}.f.u_{t}^{%
\frac{n+2}{n-2}}\text{ with }\int_{M}fu_{t}^{2^{*}}dv_{\mathbf{g}}=1 
\]
where $\lambda_{t}\rightarrow \lambda$ are the infinimum of the associated functionals.

As $\triangle _{\mathbf{g}}+h$ is coercive, $(u_{t})$ is bounded in $H_{1}^2$, thus there is a $u\in H_{1}^2$ such that 
$u_{t} \stackrel{H_{1}^2}{\rightharpoondown }u$. $u$ is a weak solution of $\E$, thus by standard elliptic theory $u$ is in fact a 
strong solution. Then, the maximum principle tells us that:
\begin{center}
{either $u>0$, or $u\equiv 0$}
\end{center}
If $u>0$, one shows with known techniques that $u$ is in fact a true minimizing solution of $\E$, and so the theorem is proved.

So, all the dificulty is to avoid the null solution. We therefore proceed by contradiction and assume that $u_{t} \stackrel{H_{1}^2}{\rightharpoondown }0$, 
and therefore that $u_{t} \stackrel{L^p}{\rightarrow }0$ for any $p<2^{*}$. The idea is then that, because of the constraint 
$\int_{M}fu_{t}^{2^{*}}dv_{\mathbf{g}}=1$, all the "mass" of the funcions $u_{t}$ concentrates around one point. 
Remember that we define a point of concentration as a point $x \in M$ such that
$$\forall\delta >0 : \overline{lim} \int_{B(x, \delta)} u_{t}^{2^{*}} >cst>0 $$
As $M$ is compact, it is easy to see that there exists at least one point of concentration. The first very important point to prove is that, 
after extraction of a subsequence:

1/ There exists precisely one point of concentration, $x_{0}$, and it is a point where $f$ is maximum on $M$.

Then, the goal of the study of this concentration phenomenom is to get a good descrition of the behaviour of the $(u_{t})$ around $x_{0}$.
One obtains the following information (up to extraction of a subsequence ):

2/ $u_{t}\rightarrow 0$ in $C^0_{loc} (M-\{x_{0}\})$

3/ there exists a sequence of points $(x_{t})$ converging to $x_{0}$ such that $u_{t}(x_{t})=\stackunder{M}{Sup}(u_{t})$.

4/ one obtains estimates of the form $\int_{B(x_{t},\delta)}d(x_{t},.)^p u_{t}^{\alpha}\sim (Sup (u_{t}))^{-\beta}$ where $p,\alpha, \beta$ 
are positive constants.

2/ and 3/ are obtained by Möser iterative process (one multiply equation $E_{t}$ by increasing powers of $u_{t}$ and integrate, getting in 
this way bounds on increasing $L^p$-norms of the $u_{t}$); 4/ is obtained by blow-up (one transfers the equation $E_{t}$ and integrals 
of the form in 4/ in exponential charts $exp^{-1}_{x_{t}}$ and multiply the transfered metric $(exp^{-1}_{x_{t}})^{*} \g$ by some power 
$(Sup\, u_{t})^k$; when $t\rightarrow \infty $, $(Sup\, u_{t})\rightarrow \infty $, and the equation and the integrals 
"converge" giving the required informations.)

Intuitively, the image is the following:
\\
(sketch)
\\
The $u_{t}$ "concentrate" around $x_{0}$. The idea is that we can do "as if" the $u_{t}$ were with compact support around $x_{0}$. We 
obtain in fact very precise information on the shape of the $u_{t}$ around $x_{0}$.

To get a contradiction, the idea is the following: we want to contradict the hypothesis
\[
\frac{4(n-1)}{n-2}h(x_{0})>S_{\mathbf{g}}(x_{0})-\frac{n-4}{2}\frac{\bigtriangleup _{\mathbf{g}}f(x_{0})}{f(x_{0})}
\]
as $x_{0}$ is a point of maximum of $f$.
We want to use the euclidean Sobolev inequality:
$$K(n,2)^{-2}(\int_{\Bbb{R}^{n}}u_{t}^{2^{*}})^{\frac{2}{2^{*}}}\leq \int_{\Bbb{R}%
^{n}}\left| \nabla u_{t}\right| ^{2} $$
in which we would inject equation $E_{t}$, that is we would like to integrate by part the gradient term and then replace the euclidean laplacian that would appear 
by its value taken from the equation.

IF the $u_{t}$ were with compact support in a small neighbourghood around $x_{0}$, IF the metric were euclidean, IF the points of maximum 
of the $u_{t}$ were all in $x_{0}$, then we would obtain quickly a contradiction as explained at the beginning of this section.
But, the $u_{t}$ are not supported in a small neibhourghood of $x_{0}$, this requires the use of cut-off functions; the metric is not euclidean, 
this requires to expand the euclidean laplacian and the euclidean measure $dx$ with respect to the laplacian $\triangle _{\g}$ and the 
measure $dv_{\g}$; and the points of maximum of the $u_{t}$ are not in $x_{0}$, which is the main difficulty introduced by the non-constant 
function $f$. The technique is then the following:

We read everything in charts $exp^{-1}_{x_{t}}$; it is in these charts that we write the Sobolev inequality for $\eta u_{t}$ where $\eta$ 
is a cut-off function; all the integrals and functions are to be understood as read in these exponential charts.
$$K(n,2)^{-2}(\int_{\Bbb{R}^{n}}(\eta u_{t})^{2^{*}} dx)^{\frac{2}{2^{*}}}\leq \int_{\Bbb{R}%
^{n}}\left| \nabla_{euc} (\eta u_{t}) \right| ^{2}_{euc} dx $$
In this inequality we make some integrations by parts, replace the euclidean laplacian appearing by $\triangle _{\g}$ and expand $dx$ 
whith respect to $dv_{\g}$, and finaly write that $\triangle _{\mathbf{g}}u_{t}=\lambda _{t}.f.u_{t}^{\frac{n+2}{n-2}}-h_{t}.u_{t}$ 
to obtain:
$$ \int_{B(x_{t},\delta)}h_{t}u_{t}^2\leqslant \,terms\, in \, R_{ij}\int_{B(x_{t},\delta)}x^{i}x^{j}u_{t}^{2^{*}}
+\,terms\,in\,\int_{B(x_{t},\delta)}(f-Sup f)u_{t}^{2^{*}}$$
(Here $R_{ij}$ is the Ricci tensor). Now, the estimates 4/ allow us to rewrite the first two terms to obtain:
$$h(x_{0})(Sup\, u_{t})^{\frac{-4}{n-2}}+o((Sup\, u_{t})^{\frac{-4}{n-2}})\leqslant S_{\g}(x_{0})(Sup\, u_{t})^{\frac{-4}{n-2}}
+o((Sup\, u_{t})^{\frac{-4}{n-2}})+\,terms\,in\,\int_{B(x_{t},\delta)}(f-Sup f)u_{t}^{2^{*}}$$
We would therefore like to prove that the last term is equivalent to 
$-\frac{n-4}{2}\frac{\bigtriangleup _{\mathbf{g}}f(x_{0})}{f(x_{0})}(Sup\, u_{t})^{\frac{-4}{n-2}}+o((Sup\, u_{t})^{\frac{-4}{n-2}})$.
For this, we expand $f$. BUT, if we develop$f$ in $x_{t}$ to use the estimates 4/ centered in $x_{t}$, the first derivatives of $f$ appear, 
but they give integrals whose order in less than $(Sup\, u_{t})^{\frac{-4}{n-2}}$:
$$\int\partial _{i}f(x_{t})x^{i}u_{t}^{2^{*}}\sim order<(Sup\, u_{t})^{\frac{-4}{n-2}}$$
To overcome this difficulty, as suggested by H. Hebey, the idea is to use a very strong theorem of O. Druet and F. Robert, which says 
that we can do as if the $u_{t}$ are radial, therefore the above integral is 0.

The other solution is to expand $f$ in $x_{0}$, because then the first derivative are 0 as $x_{0}$ is a maximum point. But then, we have to 
translate the estimates 4/ in $x_{0}$ to have $\int_{B(x_{t},\delta)}d(x_{0},.)^p u_{t}^{\alpha}\sim (Sup (u_{t}))^{-\beta}$. 
This requires an estimate on the speed of convergence of $(x_{t})$ to $x_{0}$ whith respect to $(Sup u_{t})$. Precisely we need:
$$ d(x_{t},x_{0})\leqslant c. (Sup\, u_{t})^{\frac{-2}{n-2}}$$
Then we can relace $x_{t}$ by $x_{0}$ in the estimates and conclude. I obtain this estimate with the additional hypothesis that the Hessian 
of $f$ is non-degenerate at the points of maximum (theorem 6). This estimate is very important in the study of concentration phenomena, and has been 
studied by various authors in similar settings, often in the case $f=cste$, but then requiring hypothesis on the geometry of the manifold. 
Here it seems that the hypothesis on $f$ fixes the position of the concentration point and impose  the speed of convergence. 

\section{Proof of theorem 1}
\subsection{Setup}
Let $h$ be a weakly critical function for $f$ and $\mathbf{g}$
such that for any $P\in M$ where $f$ is maximum on $M$ we have : 
\[
\text{ }h(P)>\frac{n-2}{4(n-1)}S_{\mathbf{g}}(P)-\frac{(n-2)(n-4)}{8(n-1)}%
\frac{\bigtriangleup _{\mathbf{g}}f(P)}{f(P)} 
\]
and such that there exist a family $(h_{t}),\,h_{t}\lvertneqq
h,\,h_{t}$ sub-critical for every $t$, and satisfying $h_{t}\stackunder{%
t\rightarrow t_{0}}{\rightarrow }h\,\,$in$\,$ $C^{0,\alpha}$. To simplify, we suppose that $t_{0}=1$ and that 
$t\rightarrow 1$. Then for every $t$ :
\[
\lambda _{t}:=\lambda _{h_{t},f,\mathbf{g}}<\frac{1}{K(n,2){{}^{2}}(%
\stackunder{M}{Sup}f)^{\frac{n-2}{n}}} 
\]
and there exist a family $u_{t}$ of minimizing solutions of the equations
\[
E_{t}:\,\,\triangle _{\mathbf{g}}u_{t}+h_{t}.u_{t}=\lambda _{t}.f.u_{t}^{%
\frac{n+2}{n-2}}\text{ with }\int_{M}fu_{t}^{2^{*}}dv_{\mathbf{g}}=1 
\]
We then see, as $\triangle _{\mathbf{g}}+h$ is coercive, that the sequence $(u_{t})$ is bounded in 
$H_{1}^{2}$ (just multiply $E_{t}$ by $u_{t}$ and integrate on $M$). Thus, there exist a function $u\in H_{1}^{2}\,,\,u\geq 0$ 
such that, after extracting a subsequence, 
\begin{eqnarray*}
&&u_{t}\stackrel{H_{1}^{2}}{\rightharpoondown }u\,,\, \\
&&u_{t}\stackrel{L^{2}}{\rightarrow }u\,, \\
&&\,u_{t}\stackrel{p.p.}{\rightarrow }u\,,
\end{eqnarray*}
and we can suppose
\[
\lambda _{t}\stackrel{<}{\rightarrow }\lambda \leqslant \frac{1}{K(n,2){%
{}^{2}}(\stackunder{M}{Sup}f)^{\frac{n-2}{n}}}\,\,. 
\]
In particular
\[
u_{t}\stackrel{L^{p}}{\rightarrow }u,\,\forall p<2^{*}=\frac{2n}{n-2} 
\]
as the inclusion of $H_{1}^{2}$ in $L^{p}$ is compact $\forall
p<2^{*} $. Therefore $u$ is a weak solution of
\[
\,\triangle _{\mathbf{g}}u+h.u=\lambda .f.u^{\frac{n+2}{n-2}} 
\]
and by standard elliptic theory, $u$ is $C^{\infty }$. The maximum principle then gives us that either $u>0$ or 
$u\equiv 0$.

If $u>0$ then, using elliptic theory and iteration process, and the fact that $h$ is weakly critical, one can prove that:
\[
\lambda =\frac{1}{K(n,2){{}^{2}}(\stackunder{M}{Sup}f)^{\frac{n-2}{n}}} 
\]
and then that $u$ is a minimizing positive solution of
\[
\triangle _{\mathbf{g}}u+h.u=\frac{1}{K(n,2){{}^{2}}(\stackunder{M}{Sup}f)^{%
\frac{n-2}{n}}}.f.u^{\frac{n+2}{n-2}}\text{ with }\int_{M}fu^{2^{*}}dv_{%
\mathbf{g}}=1 
\]
and the theorem is proved.

If $u\equiv 0$, we will show that there is a concentration phenomena. All the study that follows will aim at finding a 
contradiction. From now, we suppose that we are in this case:
\[
u\equiv 0. 
\]

\subsection{Concentration phenomena}
In this section we study the behavior of a family of $C^{2,\alpha }$ solutions $(u_{t})$ of
\[
\bigtriangleup _{\mathbf{g}}u_{t}+h_{t}u_{t}=\lambda _{t}fu_{t}^{\frac{n+2}{%
n-2}}\text{ with }\int_{M}fu_{t}^{\frac{2n}{n-2}}dv_{\mathbf{g}}=1 
\]
where $f$ is a smooth function such that $\underset{M}{Sup}f>0$. We also suppose that 
$h_{t}\rightarrow h$ in $C^{0,\alpha }$ where $h$ is such that $\bigtriangleup _{\mathbf{g}}+h$ is coercive.
 The sequence $(u_{t}) $ is bounded in $H_{1}^{2}$, therefore,up to a subsequence, $u_{t}\rightharpoondown u$
weakly in $H_{1}^{2}$, and we supose that $u\equiv 0$; that is $u_{t}\rightarrow 0$ in any $L^{p}$ for $p<2^{*}$. 
We also make the following "minimal energy" hypothesis:
\[
\lambda _{t}\leq \frac{1}{K(n,2)^{2}(\stackunder{M}{Sup}f)^{\frac{n-2}{n}}} 
\]
and we can suppose that $\lambda_{t}\rightarrow \lambda$. All this hypothesis are satisfied by the $u_{t}$ of the preceding section. 
The results of this section are valid for $dimM=3$, exept 
$L^{2}$-concentration, valid for $\dim M\geq 4$. In all this text, $c,C$ are constants independant ot $t$ and $\delta $.
\begin{proposition}There exist, after extraction of a subsequence, exactly one concentration point $x_{0}$, and it is 
a point where $f$ is maximum on $M$. Moreover
\[
\forall \delta >0,\,\,\,\overline{\stackunder{t\rightarrow 1}{\lim }}%
\int_{B(x_{0},\delta )}fu_{t}^{2^{*}}dv_{\mathbf{g}}=1 
\]
\end{proposition}
\textit{Proof :}
We apply the iteration process. First, as $M$
is compact, there exist at least one point of concentration. Otherwise, we could cover $M$ by a finite number of balls 
$B(x_{i},\delta )$ such that $\stackunder{t\rightarrow 1}{\lim }\int_{B(x_{i},\delta
)}u_{t}^{2^{*}}=0$, and we would have $\stackunder{t\rightarrow 1}{\lim }\int_{M}u_{t}^{2^{*}}=0,$ which 
would contradict
\[
1=\int_{M}fu_{t}^{2^{*}}dv_{\mathbf{g}}\leq Sup\left| f\right|
\int_{M}u_{t}^{2^{*}}dv_{\mathbf{g}} 
\]
The principle of iteration process is the following: if we find, for a point $x$, a cut-off function $\eta$ equal 
to 1 around $x$ such that $Q(t,k,\eta)\geq Q>0$, we get, using formula (6) or (7), that 
$(\eta u_{t}^{\frac{k+1}{2}})$ is bounded in $L^{2^{*}}$, and therefore we can extract a 
subsequence such that $(\eta u_{t})$ converges strongly to 0 in $L^{2^{*}}$; thus $x$ cannot be a concentration 
point.

Let us prove now that we can do this for a point $x$ such that $f(x)\leq 0$. If $f(x)<0$, we choose $\delta $ small enough 
such that $f<0$
on $B(x,\delta )$ and we choose $\eta $ with support in $B(x,\delta )$.
As $(u_{t})$ is bounded in $H_{1}^{2}$
and thus in $L^{2^{*}}$, we get using formula (5), that for any $k$ such that $1\leq
k\leq 2^{*}-1$: 
\[
\frac{4k}{(k+1)^{2}}\int_{M}\left| \nabla (\eta u_{t}^{\frac{k+1}{2}%
})\right| ^{2}\leq \int_{M}(\frac{2}{k+1}\left| \nabla \eta \right| ^{2}+%
\frac{2(k-1)}{(k+1)^{2}}\eta \triangle \eta -\eta ^{2}h_{t})u_{t}^{k+1}\leq
C_{1} 
\]
where $C_{1}$ is independent of $t$. Therefore for any $k$ such that $1\leq
k\leq 2^{*}-1$ there exist $C_{2}$ independent of $t$ such that: 
\[
\int_{M}\left| \nabla (\eta u_{t}^{\frac{k+1}{2}})\right| ^{2}\leq C_{2} 
\]
Therefore $(\eta u_{t}^{\frac{k+1}{2}})$ is bounded in $H_{1}^{2}$ and, using Sobolev inequality, $(\eta
u_{t}^{\frac{k+1}{2}})$ is bounded in $L^{2^{*}}$ for any $k$ such that $1\leq k\leq 2^{*}-1$.

If $f(x)=0$, by continuity of $f$ and choosing $\delta $ small 
enough, we get in (7) that for any $k$ such that $1\leq k\leq 2^{*}-1$, $%
Q(t,k,\eta)\geq Q>0$. Therefore, as we said, here again $(\eta u_{t}^{\frac{k+1}{2}})$ is bounded in $L^{2^{*}}$, 
and therefore we can extract a 
subsequence such that $(\eta u_{t})$ converges strongly to 0 in $L^{2^{*}}$. 
Thus, when $f(x)\leq0$, $x$ cannot be a concentration point.

Now, let $x$ be a concentration point: $f(x)>0$ as we just saw. 
For $\delta >0$ such that $f\geq 0$ on $B(x,\delta )$, set
\[
\stackunder{t\rightarrow 1}{\lim \sup }\int_{B(x,\delta
)}fu_{t}^{2^{*}}=a_{\delta } 
\]
Then $a_{\delta }\leq 1$ as $\int_{M}fu_{t}^{2^{*}}=1.$ Suppose that there exist $\delta >0$ such that 
$a_{\delta }<1$. Because
\[
\lambda _{t}\stackunder{\leqslant }{\rightarrow }\lambda \leqslant \frac{1}{%
K(n,2){{}^{2}}(\stackunder{M}{Sup}f)^{\frac{n-2}{n}}} 
\]
we get 
\[
\overline{\stackunder{t\rightarrow 1}{\lim }}\,\lambda _{t}K(n,2){{}^{2}}(%
\stackunder{M}{Sup}f)^{\frac{n-2}{n}}a_{\delta }<1. 
\]
Beside, $\frac{4k}{(k+1)^{2}}\stackunder{k\stackunder{>}{\rightarrow }1%
}{\rightarrow }1$. Therefore, for $k$ close to 1 such that
\[
\overline{\stackunder{t\rightarrow 1}{\lim }}\,\lambda _{t}K(n,2){{}^{2}}(%
\stackunder{M}{Sup}f)^{\frac{n-2}{n}}a_{\delta }<\frac{4k}{(k+1)^{2}} 
\]
we get, taking $\eta $ with support in $B(x,\delta )$, that in formula (6): 
$Q(t,k,\eta)\geq Q>0$ for all $t$, where $Q$ is independent of $t$. So, as before, $x$ cannot be a concentration point, 
and we have a contradiction.
Thus $a_{\delta }=1,\,\forall \delta >0$. Therefore $x$
is the only concentration point, that we will now denote $x_{0}$.
The same reasonning shows that, necessarilly,
\[
\lambda =\frac{1}{K(n,2){{}^{2}}(\stackunder{M}{Sup}f)^{\frac{n-2}{n}}}\,\,. 
\]
In the same way, if $f(x_{0})\neq \stackunder{M}{Sup}f$, there exist $\delta >0$ such that $\stackunder{B(x_{0},\delta )}{Sup}f<\stackunder{M}{Sup}%
f $. But $\lambda _{t}\leq \frac{1}{K(n,2){{}^{2}}(\stackunder{M}{Sup}f)^{%
\frac{n-2}{n}}}$, so
\[
\overline{\stackunder{t\rightarrow 1}{\lim }}\,\lambda _{t}K(n,2){{}^{2}}(%
\stackunder{B(x_{0},\delta )}{Sup}f)^{\frac{n-2}{n}}(\int_{B(x_{0},\delta
)}fu_{t}^{2^{*}})^{\frac{2^{*}-2}{2^{*}}}<1\,\,. 
\]
Then for $k$ close enough to 1, taking $\eta $ with support in $%
B(x_{{0}},\delta )$, we get in (6): $Q(t,k,n)\geq
Q>0 $ for all $t$; and once again we have a contradiction. Therefore $f(x_{0})=\stackunder{M}{Sup}f>0$.

Note that this is the main particularity introduced by the function $f$ on the right of equation $\E$. It gives a precise 
location for the concentration point.

The next propositions concerning the concentration phenomenom are now quite standard, even though they are mostly 
published in the case $f=constant$ and often with few details. We shall therefore give possible proofs, refering to the  
books \cite{D-H} and \cite{DHR} for more information, the presence of a function $f$ introducing only slight modifications that we will indicate 
when necessary.

\begin{proposition}
$u_{t}\rightarrow 0$ in $C_{loc}^{0}(M-\{x_{0}\})$.
\end{proposition}
\textit{Proof :}
It is a typical aplication of the iteration process in standard elliptic theory.
First step: Let $q>0$ be fixed. We prove that for any $\delta >0$, there exists $C=C(\delta ,q)$ independent of 
$t$ such that for $t$ close enough to 1: 
\[
\left\| u_{t}\right\| _{L^{q}(M\backslash B(x_{0},\delta ))}\leq C\left\|
u_{t}\right\| _{L^{2}(M)}\,\,. 
\]
To apply the iteration process, we build a sequence $\eta _{1},...,\eta
_{m}$ of $m$ cut-off functions such that $\eta _{j}=0$ on $B(x_{0},\delta /2)$
and $\eta _{j}=1$ on $M\backslash B(x_{0},\delta )$ and such that
\[
M\backslash B(x_{0},\delta )\subset ...\subset \{\eta _{j+1}=1\}\subset
Supp\,\eta _{j+1}\subset \{\eta _{j}=1\}\subset ...\subset M\backslash
B(x_{0},\delta /2) 
\]
and where $m$ is chosen such that $2(\frac{2^{*}}{2})^{m}>q$. We set $q_{1}=2$ and 
$q_{j}=(\frac{2^{*}}{2})q_{j-1}$. The iteration process (6), (7), gives that
\[
Q(t,q_{j}-1,\eta _{j}).(\int_{M}(\eta _{j}u_{t}^{\frac{q_{j}}{2}})^{2^{*}})^{%
\frac{n-2}{n}}\leq (\frac{4(q_{j}-1)}{q_{j}{{}^{2}}}B+C_{0}+C_{\eta
_{j}})\int_{Supp\,\eta _{j}}u_{t}^{q_{j}}\,\,\,. 
\]
But for $j\leq m$ we have $\frac{4(q_{j}-1)}{q_{j}{{}^{2}}}\geq c>0$ and 
from proposition 2, $\int_{Supp\,\eta _{j}}u_{t}^{2^{*}}\rightarrow 0$, therefore in (7), 
\[
Q(t,q_{j}-1,\eta _{j})\geq c>0,\,\forall j. 
\]
Thus there exists a neighborhood $V_{j}$ of 1 and a constant $C_{j}>0$ such that for $t\in V_{j}$: 
\[
(\int_{M}(\eta _{j}u_{t}^{\frac{q_{j}}{2}})^{2^{*}})^{\frac{n-2}{n}}\leq
C_{j}\int_{Supp\,\eta _{j}}u_{t}^{q_{j}}\,\,\,\,. 
\]
Then by construction of the $\eta _{j}$ we have
\[
(\int_{\{\eta _{j}=1\}}u_{t}^{q_{j}\frac{2^{*}}{2}})^{\frac{n-2}{n}}\leq
C_{j}\int_{\{\eta _{j-1}=1\}}u_{t}^{q_{j}}\, 
\]
and thus
\[
\left\| u_{t}\right\| _{L^{q}(M\backslash B(x_{0},\delta ))}\leq
C(\prod\limits_{j=1}^{m}C_{j})\left\| u_{t}\right\| _{L^{2}(M)}\,\,\forall
t\in V_{1}\cap ...\cap V_{m}\,. 
\]

Second step: By Gilbarg-Trudinger theorem (8.25) \cite{G-T}, we have : if $u$ is solution of an equation 
$E:\,\,\triangle _{\mathbf{g}}u+h.u=F$, where $\triangle _{\mathbf{g}}+h$ is coercive, and if 
$\omega \subset \subset \omega ^{\prime }$ are two open set, for $r>1,\,q>n/2\,:$%
\[
\stackunder{\omega }{Sup}\,u\leq c\left\| u\right\| _{L^{r}(\omega ^{\prime
})}+c^{\prime }\left\| F\right\| _{L^{q}(\omega ^{\prime })}\,\,.
\]
This theorem is also an application of the iteration process. We apply it to 
$E_{t}:\,\,\triangle _{\mathbf{g}}u_{t}+h_{t}.u_{t}=\lambda _{t}.f.u_{t}^{\frac{n+2}{n-2}}$ 
and to $\omega \subset \subset
\omega ^{\prime }\subset M\backslash \{x_{0}\}$.

Then with the first step applied to $q\frac{n+2}{n-2}$,
and chosing
\[
\omega =M\backslash B(x_{0},\delta ),\,\omega ^{\prime }=M\backslash
B(x_{0},\delta /2),\,r=2,\,q>n/2 
\]
we obtain
\begin{eqnarray*}
\stackunder{M\backslash B(x_{0},\delta )}{Sup}u_{t} & \leq &c\left\|
u_{t}\right\| _{L^{2}(\omega ^{\prime })}+c^{\prime }\lambda _{t}^{q}\left\|
u_{t}\right\| _{L^{q\frac{n+2}{n-2}}(\omega ^{\prime })}^{\frac{n+2}{n-2}}
\\ 
& \leq &c\left\| u_{t}\right\| _{L^{2}(M)}+c^{\prime \prime }\left\|
u_{t}\right\| _{L^{2}(M)}^{\frac{n+2}{n-2}}
\end{eqnarray*}
But $\left\| u_{t}\right\| _{L^{2}(M)}\rightarrow 0$, thus the result.

We recall now the notations of subsection (3.2): we consider a sequence of points $(x_{t})$ such that
\[
m_{t}=\stackunder{M}{Max}\,u_{t}=u_{t}(x_{t}):=\mu _{t}^{-\frac{n-2}{2}}. 
\]
From proposition 3, $x_{t}\rightarrow x_{0}$ and $\mu
_{t}\rightarrow 0$. Remember that $\overline{u}_{t},\overline{f}_{t},%
\overline{h}_{t},\,\mathbf{g}\,_{t}$ are the functions and the metric "viewed" in the chart $\exp _{x_{t}}^{-1}$, and $\,\,\widetilde{%
u}_{t}\,,\,\widetilde{h}_{t}\,,\,\widetilde{f}_{t},\widetilde{\,\mathbf{g}\,}%
_{t}$ are the functions and the metric after
blow-up. From now, all the blow-up's will be made on balls $B(x_{t},\delta )$ where $%
f\geq 0$, which is possible as $f(x_{0})>0$.

\begin{proposition}

$\forall R>0$ : $\stackunder{t\rightarrow 1}{\lim }\int_{B(x_{t},R\mu
_{t})}fu_{t}^{2^{*}}dv_{\mathbf{g}}=1-\varepsilon _{R}$ where $\varepsilon
_{R}\stackunder{R\rightarrow +\infty }{\rightarrow }0.$
\end{proposition}

\textit{Proof:} 
This is a direct application of blow-up analysis
in $x_{t}$ with $k_{t}=\mu _{t}^{-1}$:
\begin{center}
$\widetilde{u}_{t}\rightarrow \widetilde{u}=(1+\frac{\lambda f(x_{0})}{n(n-2)%
}\left| x\right| ^{2})^{-\frac{n-2}{2}}\,=(1+\frac{f(x_{0})^{\frac{2}{n}}}{%
K(n,2)^{2}n(n-2)}\left| x\right| ^{2})^{-\frac{n-2}{2}}$ in $\,C_{loc}^{2}(%
\Bbb{R}^{n})$.
\end{center}
Then:
\begin{center}
$\int_{B(x_{t},R\mu _{t})}fu_{t}^{2^{*}}dv_{\mathbf{g}}=\int_{B(0,R)}%
\widetilde{f}_{t}.\widetilde{u}_{t}^{2^{*}}dv_{\widetilde{\,\mathbf{g}\,}%
_{t}}\stackunder{t\rightarrow 1}{\rightarrow }f(x_{0})(\int_{B(0,R)}%
\widetilde{u}^{2^{*}}dx)\stackunder{R\rightarrow \infty }{\rightarrow }1$
\end{center}

\begin{proposition}\textbf{Weak estimates, first part.}

$\exists C>0$ such that $\forall x\in M:\,d_{\mathbf{g}}(x,x_{t})^{\frac{n-2}{2}}u_{t}(x)\leq C$.
\end{proposition}

\textit{Proof :}
Define $w_{t}(x)=\,d_{\mathbf{g}}(x,x_{t})^{%
\frac{n-2}{2}}u_{t}(x)$. We want to prove that there exists $C>0$ such that $\stackunder{M}{Sup}\,w_{t}\leq C$. 
By contradiction, we suppose that (for a subsequence) $\stackunder{M}{Sup}\,w_{t}\rightarrow +\infty $. 
Let $y_{t}$ be a point where $w_{t}$ is maximum. $M$ being compact, $d_{\mathbf{g}}(x,x_{t})$ 
is bounded, therefore $u_{t}(y_{t})\rightarrow \infty $, and thus from proposition 3, $y_{t}\rightarrow x_{0}$. 
Besides, the definition of $\mu _{t}$ gives: 
\[
\frac{d_{\mathbf{g}}(y_{t},x_{t})}{\mu _{t}}\rightarrow +\infty \,\,. 
\]
We now do a \textit{blow-up} of center $y_{t}$ and coefficient $k_{t}=u_{t}(y_{t})^{\frac{2}{n-2}}$ 
and with $m_{t}=u_{t}(y_{t})$. We obtain (taking the notation of 3.2 for this case) :
\[
\,\triangle _{\widetilde{\mathbf{g}}_{t}}\widetilde{u}_{t}+u_{t}(y_{t})^{-%
\frac{4}{n-2}}.\widetilde{h}_{t}.\widetilde{u}_{t}=\lambda _{t}\widetilde{f}%
_{t}.v_{t}^{\frac{n+2}{n-2}}\,\,. 
\]
If $x\in B(0,2):$%
\begin{eqnarray*}
d_{\mathbf{g}}(x_{t},\exp _{y_{t}}(u_{t}(y_{t})^{-\frac{2}{n-2}}x) & \geq& d_{%
\mathbf{g}}(y_{t},x_{t})-2u_{t}(y_{t})^{-\frac{2}{n-2}} \\ 
& \geq& u_{t}(y_{t})^{-\frac{2}{n-2}}(w_{t}(y_{t})^{\frac{2}{n-2}}-2)\sim d_{%
\mathbf{g}}(y_{t},x_{t})
\end{eqnarray*}
as $w_{t}(y_{t})\rightarrow \infty $ and $u_{t}(y_{t})\rightarrow \infty $.
Therefore, for $t$ close to 1: 
\[
d_{\mathbf{g}}(x_{t},\exp _{y_{t}}(u_{t}(y_{t})^{-\frac{2}{n-2}}x)\geq \frac{%
1}{2}d_{\mathbf{g}}(y_{t},x_{t})\,\,. 
\]
By consequence, for any $R>0$ and $t$ close to 1: 
\[
B(y_{t},2u_{t}(y_{t})^{-\frac{2}{n-2}})\cap B(x_{t},R\mu _{t})=\emptyset 
\]
Therefore, by proposition 4,
\begin{eqnarray*}
\int_{B(0,2)}\widetilde{f}_{t}.\widetilde{u}_{t}^{2^{*}}dv_{\widetilde{\,%
\mathbf{g}\,}_{t}} & =\int_{B(y_{t},2u_{t}(y_{t})^{-\frac{2}{n-2}%
})}fu_{t}^{2^{*}}dv_{\mathbf{g}} & \leq \int_{M\backslash B(x_{t},R\mu
_{t})}fu_{t}^{2^{*}}dv_{\mathbf{g}} \\ 
&  & \leq \int_{M}fu_{t}^{2^{*}}dv_{\mathbf{g}}-\int_{B(x_{t},R\mu
_{t})}fu_{t}^{2^{*}}dv_{\mathbf{g}} \\ 
&  & \stackunder{t\rightarrow 1,R\rightarrow \infty }{\longrightarrow }0
\end{eqnarray*}
But the iteration process then gives that for $1\leq k\leq 2^{*}-1:$%
\[
\int_{B(0,1)}\widetilde{u}_{t}^{\frac{k+1}{2}2^{*}}dv_{\widetilde{\,\mathbf{g%
}\,}_{t}}\rightarrow 0 
\]
and by iteration we obtain that $\forall p\geq 1:$%
\[
\int_{B(0,1)}\widetilde{u}_{t}^{p}dv_{\widetilde{\,\mathbf{g}\,}%
_{t}}\rightarrow 0 
\]
We deduce that $\left\| \widetilde{u}_{t}\right\| _{L^{\infty
}(B(0,1))}\rightarrow 0$ whereas $\widetilde{u}_{t}(0)=1$. Thus a contradiction.

\begin{proposition}\textbf{Weak estimates, second part.}

$\forall \varepsilon >0$ , $\exists R>0$ such that $\forall t,\,\forall x\in
M: $%
\[
\,d_{\mathbf{g}}(x,x_{t})\geq R\mu _{t}\,\Rightarrow \,\,d_{\mathbf{g}%
}(x,x_{t})^{\frac{n-2}{2}}u_{t}(x)\leq \varepsilon . 
\]
\end{proposition}
 
 \textit{Proof :}
 We use the same method, supposing the existence of a $\varepsilon _{0}>0$ and $y_{t}\in M$ such that 
\[
\stackunder{t\rightarrow 1}{\lim }\frac{d_{\mathbf{g}}(y_{t},x_{t})}{\mu _{t}%
}=+\infty \text{ et }w_{t}(y_{t})=\,d_{\mathbf{g}}(y_{t},x_{t})^{\frac{n-2}{2%
}}u_{t}(y_{t})\geq \varepsilon _{0} 
\]
We do a blow-up of center $y_{t}$ and coefficient $k_{t}=u_{t}(y_{t})^{\frac{2}{n-2}}$ and with 
$m_{t}=u_{t}(y_{t}).$

Then, with these hypothesis, if $x\in B(0,\frac{1}{2}\varepsilon _{0}^{%
\frac{2}{n-2}}):$%
\begin{eqnarray*}
d_{\mathbf{g}}(x_{t},\exp _{y_{t}}(u_{t}(y_{t})^{-\frac{2}{n-2}}x) && \geq d_{%
\mathbf{g}}(y_{t},x_{t})-\frac{1}{2}\varepsilon _{0}^{\frac{2}{n-2}%
}u_{t}(y_{t})^{-\frac{2}{n-2}} \\ 
& &\geq \frac{1}{2}d_{\mathbf{g}}(y_{t},x_{t})
\end{eqnarray*}
Therefore for any $R>0$ and $t$ close to 1: 
\[
B(y_{t},\frac{1}{2}\varepsilon _{0}^{\frac{2}{n-2}}u_{t}(y_{t})^{-\frac{2}{%
n-2}})\cap B(x_{t},R\mu _{t})=\emptyset 
\]
Therefore, as previously: 
\[
\int_{B(0,\frac{1}{2}\varepsilon _{0}^{\frac{2}{n-2}})}\widetilde{f}_{t}.%
\widetilde{u}_{t}^{2^{*}}dv_{\widetilde{\,\mathbf{g}\,}_{t}}\rightarrow 0 
\]
and we obtain in the same way a contradiction.

\begin{proposition} \textbf{$L^2$-concentration.} 

If $\dim M\geq 4,$ $\forall \delta >0\,:\,$%
\[
\stackunder{t\rightarrow 1}{\lim }\frac{\int_{B(x_{0},\delta )}u_{t}^{2}dv_{%
\mathbf{g}}}{\int_{M}u_{t}^{2}dv_{\mathbf{g}}}=1 
\]
\end{proposition}
 
 \textit{Proof :}
 We first use the two first step of the proof of proposition 3 to show that there exists $c>0$ such that: 
\[
\stackunder{M\backslash B(x_{0},\delta )}{Sup}u_{t}\leq c\left\|
u_{t}\right\| _{L^{2}(M)}\,\,. 
\]
Indeed, going over what we did there,: 
\begin{eqnarray*}
\stackunder{M\backslash B(x_{0},\delta )}{Sup}u_{t} && \leq c\left\|
u_{t}\right\| _{L^{2}(\omega ^{\prime })}+c^{\prime }\lambda _{t}^{q}\left\|
u_{t}^{\frac{n+2}{n-2}}\right\| _{L^{q}(\omega ^{\prime })} \\ 
&& \leq c\left\| u_{t}\right\| _{L^{2}(M)}+c^{\prime }\lambda _{t}^{q}%
\stackunder{\omega ^{\prime }}{Sup}(u_{t}^{\frac{n+2}{n-2}-1})\left\|
u_{t}\right\| _{L^{q}(\omega ^{\prime })}\\
&&\leq c^{\prime \prime }\left\|
u_{t}\right\| _{L^{2}(M)}
\end{eqnarray*}
as we know now that $\stackunder{\omega ^{\prime }}{Sup}(u_{t}^{\frac{%
n+2}{n-2}-1})\rightarrow 0$ and that, on the other hand, the first step of the proof of proposition 3 gives 
$\left\| u_{t}\right\| _{L^{q}(\omega ^{\prime })}\leq C\left\|u_{t}\right\| _{L^{2}(M)}$

Third step: Using this: 
\begin{eqnarray}
\left\| u_{t}\right\| _{L^{2}(M\backslash B(x_{0},\delta ))}^{2} && \leq 
\stackunder{M\backslash B(x_{0},\delta )}{Sup}u_{t}.\int_{M\backslash
B(x_{0},\delta )}u_{t}\nonumber \\ 
&& \leq c\left\| u_{t}\right\| _{L^{2}(M)}\left\| u_{t}\right\| _{L^{1}(M)}
\end{eqnarray}
We now want to prove that
\begin{equation}
\left\| u_{t}\right\| _{L^{1}(M)}\leq c\left\| u_{t}\right\|
_{L^{2^{*}-1}(M)}^{2^{*}-1}\,\,.  
\end{equation}
If $h>0$, we get the result by integrating equation $E_{t}$. Otherwise, as $\lambda _{h,f,\,\mathbf{g}\,}>0$, 
for any $q\in ]2,2^{*}[$, there exists $\varphi >0$ solution of $\triangle _{\mathbf{g}}\varphi +h\varphi
=\lambda _{h,f,\,\mathbf{g}\,}.f.\varphi ^{q-1}$. We set
\[
\mathbf{g}^{\prime }=\varphi ^{\frac{4}{n-2}}\mathbf{g}\text{ and }\overline{%
h_{t}}=\frac{\triangle _{\mathbf{g}}\varphi +h_{t}\varphi }{\varphi ^{\frac{%
n+2}{n-2}}} 
\]
Then for $t$ close to 1
$$\overline{h_{t}}=\varphi ^{q-2^{*}}-(h-h_{t})\varphi ^{2-2^{*}}\geq
\varepsilon _{0}>0 $$

Besides, by conformal invariance, and using $E_{t}$, we have: 
\[
\triangle _{\mathbf{g}^{\prime }}\overline{u_{t}}+\overline{h_{t}}.\overline{%
u_{t}}=\lambda _{t}f.\overline{u_{t}}^{\frac{n+2}{n-2}} 
\]
where $\overline{u_{t}}=\varphi ^{-1}.u_{t}$. Integrating, we obtain: 
\[
\varepsilon _{0}\int_{M}\overline{u_{t}}dv_{\mathbf{g}^{\prime }}\leq
\lambda _{t}Supf\int_{M}\overline{u_{t}}^{\frac{n+2}{n-2}}dv_{\mathbf{g}%
^{\prime }} 
\]
and thus there exists $C>0$ such that for $t$ close to 1 
\[
\left\| u_{t}\right\| _{L^{1}(M)}\leq C\left\| u_{t}\right\|
_{L^{2^{*}-1}(M)}^{2^{*}-1} 
\]
where the norms are now relative to $dv_{\mathbf{g}}$.

Fourth step: We conclude using Hölder's inequality. If $n=\dim M\geq 6$ : 
\[
\left\| u_{t}\right\| _{L^{2^{*}-1}(M)}^{2^{*}-1}\leq \left\| u_{t}\right\|
_{L^{2}(M)}^{\frac{n+2}{n-2}}Vol_{\mathbf{g}}(M)^{\frac{n-6}{2(n-2)}}\,. 
\]
With (10) and (11), we obtain :
\[
\stackunder{t\rightarrow 1}{\lim }\frac{\left\| u_{t}\right\|
_{L^{2}(M\backslash B(x_{0},\delta ))}^{2}}{\left\| u_{t}\right\|
_{L^{2}(M)}^{2}}=0 
\]
which proves the result. If $n=5$, Hölder's inequality gives: 
\[
\left\| u_{t}\right\| _{L^{2^{*}-1}(M)}^{2^{*}-1}\leq \left\| u_{t}\right\|
_{L^{2}(M)}^{\frac{3}{2}}\left\| u_{t}\right\| _{L^{2}(M)}^{\frac{5}{6}} 
\]
and we also conclude using (10) and (11). If now $n=4$,
we have to use proposition 6 and the associated blow-up. We have
\[
\frac{\left\| u_{t}\right\| _{L^{3}(M)}^{3}}{\left\| u_{t}\right\|
_{L^{2}(M)}}\leq \left\| u_{t}\right\| _{L^{\infty }(M\backslash
B(x_{0},\delta ))}\left\| u_{t}\right\| _{L^{2}(M)}+\frac{\int_{B(0,\delta
\mu _{t}^{-1})}\widetilde{u}_{t}^{3}dv_{\widetilde{\,\mathbf{g}\,}_{t}}}{%
(\int_{B(0,\delta \mu _{t}^{-1})}\widetilde{u}_{t}^{2}dv_{\widetilde{\,%
\mathbf{g}\,}_{t}})^{\frac{1}{2}}}\,\,. 
\]
Then for any $R>0$, using Hölder's inequality and proposition 6, we obtain: 
\[
\int_{B(0,\delta \mu _{t}^{-1})}\widetilde{u}_{t}^{3}dv_{\widetilde{\,%
\mathbf{g}\,}_{t}}\leq \int_{B(0,R)}\widetilde{u}_{t}^{3}dv_{\widetilde{\,%
\mathbf{g}\,}_{t}}+\varepsilon _{R}(\int_{B(0,\delta \mu _{t}^{-1})}%
\widetilde{u}_{t}^{2}dv_{\widetilde{\,\mathbf{g}\,}_{t}})^{\frac{1}{2}}\,\,. 
\]
It follows that for any $R,R^{\prime }>0$,

\[
\stackunder{t\rightarrow 1}{\lim \sup }\frac{\left\| u_{t}\right\|
_{L^{3}(M)}^{3}}{\left\| u_{t}\right\| _{L^{2}(M)}}\leq \varepsilon _{R}+%
\frac{\int_{B(0,R)}\widetilde{u}^{3}dx}{(\int_{B(0,R^{\prime })}\widetilde{u}%
^{2}dx)^{\frac{1}{2}}}\,\,. 
\]
As $\widetilde{u}\in L^{3}(\Bbb{R}^{4})$ and $\stackunder{R^{\prime
}\rightarrow \infty }{\lim }\int_{B(0,R^{\prime })}\widetilde{u}%
^{2}dx=+\infty $, we finally get
\[
\stackunder{t\rightarrow 1}{\lim \sup }\frac{\left\| u_{t}\right\|
_{L^{3}(M)}^{3}}{\left\| u_{t}\right\| _{L^{2}(M)}}=0 
\]
and we conclude once again using (10) and (11).

\begin{proposition} \textbf{Strong estimates.}

For any $\nu$, $0< \nu < n-2$, there exists a constant $C(\nu )>0$ such that
$$
\forall x\in M:\,d_{\mathbf{g}}(x,x_{t})^{n-2-\nu }\mu _{t}^{-\frac{n-2}{2}%
+\nu }u_{t}(x)\leq C(\nu )  
$$
\end{proposition}

\textit{Proof :} The proof requires the use of the Green function and of the weak estimates. The idea is due to 
O. Druet and F. Robert \cite{D-R}. We recall first the property of the Green function.
If $\bigtriangleup _{\mathbf{g}}+h$ is a coercive operator, there exists a unique function (at least $C^{2}$ 
with our hypothesis) 
\[
G_{h}:M\times M\backslash \{(x,x),x\in M\}\rightarrow \Bbb{R} 
\]
symetric and positive, such that in the sense of distributions, we have: $\forall x\in M$%
\begin{equation}
\bigtriangleup _{\mathbf{g},y}G_{h}(x,y)+h(y)G_{h}(x,y)=\delta _{x}  
\end{equation}
Furthermore, there exists $c>0,\,\rho >0$ such that $\forall (x,y)$ with $0<d_{%
\mathbf{g}}(x,y)<\rho :$%
\begin{equation}
\frac{c}{d_{\mathbf{g}}(x,y)^{n-2}}\leq G_{h}(x,y)\leq \frac{c^{-1}}{d_{%
\mathbf{g}}(x,y)^{n-2}}  
\end{equation}

\begin{equation}
\frac{\left| \nabla _{y}G_{h}(x,y)\right| }{G_{h}(x,y)}\geq \frac{c}{d_{%
\mathbf{g}}(x,y)}  
\end{equation}

\begin{center}
$c$ and $\rho $ vary continuously with $h$ 
\begin{equation}
G_{h}(x,y)d_{\mathbf{g}}(x,y)^{n-2}\rightarrow \frac{1}{(n-2)\omega _{n-1}}%
\text{ when }d_{\mathbf{g}}(x,y)\rightarrow 0  
\end{equation}
\end{center}

To prove these strong estimates, it is sufficient, considering (13), to prove that 
$\mu _{t}^{\frac{n-2}{2}-(n-2)(1-\nu)}u_{t}(x)\leq c^{\prime }G_{h}^{1-\nu}(x,x_{t})$, (just change $\nu$ 
by $(n-2)\nu$). 
First, notice that, using for example the weak estimates, the strong estimates are true in any ball 
$B(x_{t},R\mu _{t})$ where $R$ is fixed. We therefore have to prove the estimates in the manifold with 
boundary $M\backslash B(x_{t},R\mu _{t})$ whose boundary is $b(M\backslash B(x_{t},R\mu
_{t}))=bB(x_{t},R\mu _{t})$. 
For $\nu $ small, there exists $\varepsilon _{0}>0$ such that he operator
\[
\bigtriangleup _{\mathbf{g}}+\frac{h-2\varepsilon _{0}}{1-\nu } 
\]
is still coercive; let $\widetilde{G}$ be its Green function.
To prove our esimate, we apply the maximum principle to : 
\[
L_{t}\varphi =\bigtriangleup _{\mathbf{g}}\varphi +h_{t}\varphi -\lambda
_{t}fu_{t}^{2^{*}-2}\varphi 
\]
and to $x\longmapsto \widetilde{G}^{1-\nu}(x,x_{t})-c\mu _{t}^{\frac{n-2}{2}-(n-2)(1-\nu)}u_{t}(x)$.
As $L_{t}u_{t}=0$ with $u_{t}>0$, $L_{t}$ satisfies the maximum principle (see \cite{B-N-V}). 

Using (12) and the fact that $\delta _{x_{t}}(x)=0$ on $M\backslash
B(x_{t},R\mu _{t})$ some computations give that $\forall x\in M\backslash B(x_{t},R\mu _{t}):$%
\[
\frac{L_{t}\widetilde{G}^{1-\nu }}{\widetilde{G}^{1-\nu }}%
(x,x_{t})=2\varepsilon _{0}+h_{t}(x)-h(x)-\lambda
_{t}f(x)u_{t}(x)^{2^{*}-2}+\nu (1-\nu )\left| \frac{\nabla \widetilde{G}}{%
\widetilde{G}}\right| ^{2}(x,x_{t}) 
\]
But for $t$ close to 1, $h_{t}-h\geq -\varepsilon _{0}$ as $%
h_{t}\rightarrow h$ in $C^{0}$. Therefore
\begin{equation}
\frac{L_{t}\widetilde{G}^{1-\nu }}{\widetilde{G}^{1-\nu }}(x,x_{t})\geq
\varepsilon _{0}-\lambda _{t}f(x)u_{t}(x)^{2^{*}-2}+\nu (1-\nu )\left| \frac{%
\nabla \widetilde{G}}{\widetilde{G}}\right| ^{2}(x,x_{t})  
\end{equation}
We now separate $M\backslash B(x_{t},R\mu _{t})$ in two parts using a ball 
$B(x_{t},\rho )$ where $\rho >0$ is as in (13) and (14).
For $t$ close to 1, $\rho >R\mu _{t}$. $R>0$ will be fixed later.

1/:As $u_{t}\rightarrow 0$ in $C_{loc}^{0}(M\backslash \{x_{0}\})$,
(16) gives for $t$ close to 1: 
\[
\forall x\in M\backslash B(x_{t},\rho ):\,L_{t}\widetilde{G}^{1-\nu
}(x,x_{t})\geq 0. 
\]

2/: Using the weak estimates (second part), in $B(x_{t},\rho
)\backslash B(x_{t},R\mu _{t})$ :
\[
d_{\mathbf{g}}(x,x_{t})^{2}u_{t}(x)^{2^{*}-2}\leq \varepsilon _{R} 
\]
where $\varepsilon _{R}\stackunder{R\rightarrow \infty }{\rightarrow }0$.
Then, with (14) et (16), for $R$ big enough: 
\begin{eqnarray*}
\frac{L_{t}\widetilde{G}^{1-\nu }}{\widetilde{G}^{1-\nu }}(x,x_{t}) && \geq
\varepsilon _{0}-\lambda _{t}f(x)u_{t}(x)^{2^{*}-2}+\nu (1-\nu )\frac{c}{d_{%
\mathbf{g}}(x,x_{t})^{2}} \\ 
&& \geq \varepsilon _{0}-\lambda _{t}(\stackunder{B(x_{t},\rho )}{Sup}f).%
\frac{\varepsilon _{R}}{d_{\mathbf{g}}(x,x_{t})^{2}}+\nu (1-\nu )\frac{c}{d_{%
\mathbf{g}}(x,x_{t})^{2}} \\ 
&& \geq \varepsilon _{0}+\frac{c^{\prime }}{d_{\mathbf{g}}(x,x_{t})^{2}}\geq 0
\end{eqnarray*}

We have proved that in $M\backslash B(x_{t},R\mu _{t})$ and for any constant $C_{t}>0$ 
which can depend of $t$ :
\[
L_{t}(C_{t}.\widetilde{G}^{1-\nu }(x,x_{t}))=C_{t}.L_{t}\widetilde{G}^{1-\nu
}(x,x_{t})\geq 0=L_{t}u_{t} 
\]
At last, on the boundary $b(M\backslash B(x_{t},R\mu _{t}))$, using (13), we obtain :
\[
\widetilde{G}^{1-\nu }(x,x_{t})\geq \frac{c}{d_{\mathbf{g}%
}(x,x_{t})^{(n-2)(1-\nu )}}=\frac{c}{(R\mu _{t})^{(n-2)(1-\nu )}}\,\,. 
\]
So, if we let $C_{t}=c^{-1}R^{(n-2)(1-\nu )}\mu _{t}^{(n-2)(1-\nu )-\frac{%
n-2}{2}},$ we have for $x\in bB(x_{t},R\mu _{t})=b(M\backslash B(x_{t},R\mu
_{t})):$%
\[
C_{t}.\widetilde{G}^{1-\nu }(x,x_{t})\geq \mu _{t}^{-\frac{n-2}{2}%
}=Sup\,u_{t}\geq u_{t}(x) 
\]
Therefore, by the maximum principle :
\[
C_{t}.\widetilde{G}^{1-\nu }(x,x_{t})\geq u_{t}(x)\text{ in }M\backslash
B(x_{t},R\mu _{t}) 
\]
which can be rewriten
\[
\widetilde{G}^{1-\nu }(x,x_{t})\geq C_{t}^{-1}u_{t}(x)=c\text{ }\mu _{t}^{%
\frac{n-2}{2}-(n-2)(1-\nu )}u_{t}(x) 
\]
and therefore, using (13) : 
\[
d_{\mathbf{g}}(x,x_{t})^{(n-2)(1-\nu )}\mu _{t}^{\frac{n-2}{2}-(n-2)(1-\nu
)}u_{t}(x)\leq c 
\]
which gives the strong estimates by changing $\nu $ in $(n-2)\nu $.

\begin{proposition}\textbf{Corollary: Strong $L^p$-concentration.}

$\forall R>0$, $%
\forall \delta >0$\ and $\forall p>\frac{n}{n-2}$
$$
\stackunder{t\rightarrow 1}{\lim }\frac{\int_{B(x_{t},R\mu
_{t})}u_{t}^{p}dv_{\mathbf{g}}}{\int_{B(x_{t},\delta )}u_{t}^{p}dv_{\mathbf{g%
}}}=1-\varepsilon _{R}  \text{ where } \varepsilon _{R}\stackunder{%
R\rightarrow +\infty }{\rightarrow }0.
$$
\end{proposition}

\textit{Proof :}
Just apply the strong estimates to a blow-up in $x_{t}.$ By blow-up formulae
\begin{eqnarray*}
\int_{M}u_{t}^{p}dv_{\mathbf{g}} && \geq \int_{B(x_{t},\mu _{t})}u_{t}^{p}dv_{%
\mathbf{g}}=\mu _{t}^{n-\frac{n-2}{2}p}\int_{B(0,1)}\widetilde{u}_{t}^{p}dv_{%
\widetilde{g}_{t}} \\ 
& &\geq C\mu _{t}^{n-\frac{n-2}{2}p}
\end{eqnarray*}
On the other hand, by the strong estimates: 
\begin{eqnarray*}
\int_{M\backslash B(x_{t},R\mu _{t})}u_{t}^{p}dv_{\mathbf{g}} && \leq C\mu
_{t}^{p\frac{n-2}{2}}\int_{M\backslash B(x_{t},R\mu _{t})}d_{\mathbf{g}%
}(y_{t},x)^{(2-n)p}dv_{\mathbf{g}} \\ 
&& \leq C\mu _{t}^{n-p\frac{n-2}{2}}R^{n+(2-n)p}
\end{eqnarray*}
as soon as $p>\frac{n}{n-2}$. Dividing, we obtain the corollary.
\\

At this point, to carry on the proof of theorem 1, we need a powerfull extension of the strong estimates, called 
$C^0-theory$, which is in fact a complete control of the sequence 
$d_{\mathbf{g}}(x,x_{t})^{n-2}\mu _{t}^{-\frac{n-2}{2}}u_{t}(x)$; it is expressed by the next theorem of Druet 
and Robert, and proved in arbitrary energy in \cite{DHR}.

Another approach, also accessible at this point and originally used in the author's PHD thesis, 
is to prove another very important estimate 
concerning the "speed" of convergence of $(x_{t})$ to $x_{0}$, but it requires the additional hypothesis that the Hessian of $f$ is non-degenerate 
at the points of maximum of $f$; it will be our theorem 6, whose proof is independent of the theorem of Druet-Robert, 
only requiring the results up to proposition 9, and appears as a byproduct of an alternative proof of theorem 1. It is 
however of independent interest, as it is a very important estimate concerning concentration phenomena's which has been 
studied by various authors.

We now state the theorem of Druet and Robert and refer for its proof to \cite{DHR}, the function $f$ introducing no difficulties. 
It says first that one can take $\nu =0$ in the strong estimates, but also that one has somehow the reverse estimate.
\\

\begin{theo} [Druet, Robert]

For any $\ve>0$, there exist $\delta_{\ve}>0$ such that, up to a subsequence, for any $t$ and any 
$x\in B(x_{0},\delta_{\ve})$ :
$$(1-\ve) B_{t}(x) \leq u_{t}(x) \leq (1+\ve) B_{t}(x)$$
where 
$$ B_{t}(x)=\mu_{t}^{-\frac{n-2}{2}}\Bigl( 1+\frac{\lambda f(x_{0})}{n(n-2)} \frac{d_{\g}(x_{t},x)^2}{\mu_{t}^2}\Bigr) ^{-\frac{n-2}{2}} $$
is the "standard bubble".
\end{theo}
Note that in the proof of theorem 1, we will need the minoration:\\
 $(1-\ve) B_{t}(x) \leq u_{t}(x)$, which is 
a stronger result than $u_{t}(x) \leq (1+\ve) B_{t}(x)$ which must first be proved to get the minoration.

Finally, we come to our main result concerning the concentration phenomenom, which is the "missing link" between 
the sequence $(x_{t})$ and $x_{0}$.
\begin{theorem}\textbf{"Second fundamental estimate".}
Suppose that $dimM\geq5$ and that the hessian of the function $f$ is non-degenerate at each of its points of maximum. Then, there 
exist a constant $C$ such that for all $t$ :
$$\frac{d_{\mathbf{g}}(x_{t},x_{0})}{\mu _{t}}\leq C .$$
Moreover, if for each point $P$ of maximum of $f$ we have
$$h(P)=\frac{n-2}{4(n-1)}S_{\mathbf{g}}(P)-\frac{(n-2)(n-4)}{8(n-1)}\frac{%
\bigtriangleup _{\mathbf{g}}f(P)}{f(P)} ,$$
then more precisely
$$\frac{d_{\mathbf{g}}(x_{t},x_{0})}{\mu _{t}}\rightarrow 0 .$$
\end{theorem}
To understand the significance of this theorem, note that the weak and strong estimates, the strong $L^p$-concentration 
and the estimates in the theorem of Druet-Robert, are "centered" in $x_{t}$. Theorem 6 allows one to "translate" 
these estimates in $x_{0}$ in the sense that one can now replace $x_{t}$ by $x_{0}$. This estimate, called by Zoé Faget "second fundamental 
estimate", (the "first one" being the strong estimate), joined with the estimates of $C^0-theory$ presented in 
the theorem of Druet and Robert above, 
gives a complete description of the behavior of a sequence of solutions of equations 
$\bigtriangleup _{\mathbf{g}}u_{t}+h_{t}u_{t}=\lambda _{t}fu_{t}^{\frac{n+2}{n-2}}$ 
in the spirit of the study of Palais-Smale sequences associated to these equations. It has been studied, for example, by 
Druet and Robert in the case $f=constant=1$ in \cite{D-R} where they require strong hypothesis on the shape of the functions 
$h_{t}$ and on the geometry of the manifold near the concentration point, or by Hebey in the euclidean setting.
Intuitively, it seems that our hypothesis on $f$ "fixes" the position of the concentration point, and so we get a control 
on the distance between $x_{t}$ and $x_{0}$. Also, our method seems to be applicable to other settings, 
see e.g.\cite{F2} and \cite{C-H-V}.

\subsection{Proof of theorem 1}
We now apply the principle exposed in 3.4.

Remember that $\overline{u}_{t},\overline{f}_{t},%
\overline{h}_{t},\,\mathbf{g}\,_{t}$ are the functions and the metric "viewed" in the chart $\exp _{x_{t}}^{-1}$, and $\,\,\widetilde{%
u}_{t}\,,\,\widetilde{h}_{t}\,,\,\widetilde{f}_{t},\widetilde{\,\mathbf{g}\,}%
_{t}$ are the functions and the metric after
blow-up with center $x_{t}$ and coefficient $k_{t}=\mu_{t}^{-1}$. From now, all the blow-up's will be made on balls $B(x_{t},\delta )$ where $%
f\geq 0$, which is possible as $f(x_{0})>0$.

Let also $\eta $ be a cut-off function on $\Bbb{R}^{n}$ equal to 1 on the euclidean ball $B(0,\delta /2)$, 
and equal to 0 on $\Bbb{R}^{n}\backslash B(0,\delta )$, $0\leq \eta \leq 1$ with $\left|
\nabla \eta \right| \leq C.\delta ^{-1}$ where $\delta $ is chosen small enough to have $f\geq 0$ on the balls 
$B(x_{t},\delta )$. 
The Sobolev inequality gives on the one hand
\begin{equation}
(\int_{B(0,\delta )}(\eta \overline{u}_{t})^{2^{*}}dx)^{\frac{2}{2^{*}}}\leq
K(n,2)^{2}\int_{B(0,\delta )}\left| \nabla (\eta \overline{u}_{t})\right|
_{e}^{2}dx\,\, 
\end{equation}
where $\left| .\right| _{e}$ is the euclidean metric of associated measure $dx$.

On the other hand, integration by part gives, noting that $\left|
\nabla \eta \right| =$ $\Delta \eta =0$ on $B(0,\delta /2)$ : 
\[
\int_{B(0,\delta )}\left| \nabla (\eta \overline{u}_{t})\right|
_{e}^{2}dx\leq \int_{B(0,\delta )}\eta ^{2}\overline{u}_{t}\bigtriangleup
_{e}\overline{u}_{t}dx+C.\delta ^{-2}\int_{B(0,\delta )\backslash B(0,\delta
/2)}\overline{u}_{t}^{2}dx 
\]
Noting $\,\mathbf{g}\,_{t}^{ij}$ the components of  $\,\mathbf{g}\,_{t}$
and $\Gamma (\,\mathbf{g}\,_{t})_{ij}^{k}$ the associated Christoffel symbols, we write : 
\[
\bigtriangleup _{e}\overline{u}_{t}=\bigtriangleup _{\mathbf{g}_{t}}%
\overline{u}_{t}+(\,\mathbf{g}\,_{t}^{ij}-\delta ^{ij})\partial _{ij}%
\overline{u}_{t}-\,\mathbf{g}\,_{t}^{ij}\Gamma (\,\mathbf{g}%
\,_{t})_{ij}^{k}\partial _{k}\overline{u}_{t} 
\]
We get from this inequallity, using 
using this expression of the laplacian, equation $E_{t}:\,\,\triangle _{%
\mathbf{g}}u_{t}+h_{t}.u_{t}=\lambda _{t}.f.u_{t}^{\frac{n+2}{n-2}}$ ``viewed''
in the chart exp$_{x_{t}}^{-1}$, and using the fact that $\left| \nabla \eta \right|
=\Delta \eta =0$ on $B(0,\delta /2)$ and with some integration by parts: 
\begin{eqnarray*}
\int_{B(0,\delta )}\left| \nabla (\eta \overline{u}_{t})\right|
_{e}^{2}dx\leq & &\lambda _{t}\int_{B(0,\delta )}\eta ^{2}\overline{f}_{t}%
\overline{u}_{t}^{2^{*}}dx-\int_{B(0,\delta )}\eta ^{2}\overline{h}_{t}%
\overline{u}_{t}^{2}dx\\
&&+C.\delta ^{-2}\int_{B(0,\delta )\backslash B(0,\delta
/2)}\overline{u}_{t}^{2}dx \\ 
&& -\int_{B(0,\delta )}\eta ^{2}(\,\mathbf{g}\,_{t}^{ij}-\delta
^{ij})\partial _{i}\overline{u}_{t}\partial _{j}\overline{u}_{t}dx\\
&&+\frac{1}{2}\int_{B(0,\delta )}(\partial _{k}(\,\mathbf{g}\,_{t}^{ij}\Gamma (\,\mathbf{g%
}\,_{t})_{ij}^{k}+\partial _{ij}\,\mathbf{g}\,_{t}^{ij})(\eta \overline{u}%
_{t}^{2})dx\,.
\end{eqnarray*}
Using the Sobolev inequality (17) and the fact that 
 $\lambda _{t}\leq \frac{1}{K(n,2){{}^{2}}(%
\stackunder{M}{Sup}f)^{\frac{n-2}{n}}}$, we obtain at last: 

\begin{equation}
\int_{B(0,\delta )}\overline{h}_{t}(\eta \overline{u}_{t})^{2}dx\leq A_{t}
+B_{t}+C_{t}+C.\delta ^{-2}\int_{B(0,\delta )\backslash B(0,\delta /2)}\overline{u}%
_{t}^{2}dx
\end{equation}
where:

$B_{t}=\frac{1}{2}\int_{B(0,\delta )}(\partial _{k}(\,\mathbf{g}%
\,_{t}^{ij}\Gamma (\,\mathbf{g}\,_{t})_{ij}^{k}+\partial _{ij}\,\mathbf{g}%
\,_{t}^{ij})(\eta \overline{u}_{t}^{2})dx$

$C_{t}=\left| \int_{B(0,\delta )}\eta ^{2}(\,\mathbf{g}\,_{t}^{ij}-\delta
^{ij})\partial _{i}\overline{u}_{t}\partial _{j}\overline{u}_{t}dx\right| $

$A_{t}=\frac{1}{K(n,2){{}^{2}}(\stackunder{M}{Sup}f)^{\frac{n-2}{n}}}%
\int_{B(0,\delta )}\overline{f}_{t}\eta ^{2}\overline{u}_{t}^{2^{*}}dx-\frac{%
1}{K(n,2){{}^{2}}}(\int_{B(0,\delta )}(\eta \overline{u}_{t})^{2^{*}}dx)^{%
\frac{2}{2^{*}}}$

These computations were developed in the article of Djadli and Druet \cite{D-D}. Our goal is to use $L{{}^{2}}$-concentration (proposition 7) 
to obtain a contradiction; we shall divide (18) by $%
\int_{B(0,\delta )}\overline{u}_{t}^{2}dx$ and take the limit when $t\rightarrow
t_{0}=1$.

$L{{}^{2}}$-concentration first gives : 
\[
\frac{C.\delta ^{-2}\int_{B(0,\delta )\backslash B(0,\delta /2)}\overline{u}%
_{t}^{2}dx}{\int_{B(0,\delta )}\overline{u}_{t}^{2}dx}\stackunder{%
t\rightarrow 1}{\rightarrow }0\,\,. 
\]

Z.Djadli and O.Druet \cite{D-D} showed (see also [10] for full details):
\[
\stackunder{t\rightarrow 1}{\overline{\lim }}\frac{C_{t}}{\int_{B(0,\delta )}%
\overline{u}_{t}^{2}dx}\leq \varepsilon _{\delta }\text{ where }\varepsilon
_{\delta }\rightarrow 0\text{ when }\delta \rightarrow 0\,\,.
\]

Furthermore, as $x_{t}\rightarrow x_{0}$ we have $\stackunder{t\rightarrow 1}{%
\lim }(\partial _{k}(\,\mathbf{g}\,_{t}^{ij}\Gamma (\,\mathbf{g}%
\,_{t})_{ij}^{k}+\partial _{ij}\,\mathbf{g}\,_{t}^{ij})(0)=\frac{1}{3}S_{\,%
\mathbf{g}\,}(x_{0})$, therefore, using $L{{}^{2}}$-concentration : 
\[
\stackunder{t\rightarrow 1}{\overline{\lim }}\frac{B_{t}}{\int_{B(0,\delta )}%
\overline{u}_{t}^{2}dx}=\frac{1}{6}S_{\,\mathbf{g}\,}(x_{0})+\varepsilon
_{\delta }\,\,. 
\]

It is the expression $A_{t}$ which will give $\frac{n-2}{4(n-1)}%
S_{\,\mathbf{g}\,}(x_{0})$ $-\frac{1}{6}S_{\,\mathbf{g}\,}(x_{0})$ and $\frac{%
(n-2)(n-4)}{8(n-1)}\frac{\bigtriangleup _{\mathbf{g}}f(x_{0})}{f(x_{0})}$.

By Hölder's inequality: 
\[
\int_{B(0,\delta )}\overline{f}_{t}\eta ^{2}\overline{u}_{t}^{2^{*}}dx\leq
(\int_{B(0,\delta )}\overline{f}_{t}\overline{u}_{t}^{2^{*}}dx)^{\frac{2}{n}%
}(\int_{B(0,\delta )}\overline{f}_{t}(\eta \overline{u}_{t})^{2^{*}}dx)^{%
\frac{n-2}{n}}\, 
\]
Beside : 
\[
\,dx\leq (1+\frac{1}{6}Ric(x_{t})_{ij}x^{i}x^{j}+C\left| x\right| ^{3})dv_{\,%
\mathbf{g}\,_{t}} 
\]
Using this development and $(1+x)^\alpha \leq 1+\alpha x$ for $0<\alpha \leq1$: 
\[
(\int_{B(0,\delta )}\overline{f}_{t}\overline{u}_{t}^{2^{*}}dx)^{\frac{2}{n}%
}\leq (\int_{B(0,\delta )}\overline{f}_{t}\overline{u}_{t}^{2^{*}}dv_{\,%
\mathbf{g}\,_{t}})^{\frac{2}{n}}+\frac{1}{(\int_{B(0,\delta )}\overline{f}%
_{t}\overline{u}_{t}^{2^{*}}dv_{\,\mathbf{g}\,_{t}})^{\frac{n-2}{n}}}\frac{2%
}{n}\{S_{t}\}+C\{S_{t}\}^{2} 
\]
where 
\[
\{S_{t}\}=\frac{1}{6}Ric(x_{t})_{ij}\int_{B(0,\delta )}x^{i}x^{j}\overline{f}%
_{t}\overline{u}_{t}^{2^{*}}dv_{\,\mathbf{g}\,_{t}}+C\int_{B(0,\delta
)}\left| x\right| ^{3}\overline{u}_{t}^{2^{*}}dv_{\,\mathbf{g}\,_{t}}\,\,. 
\]
We deduce
\[
A_{t}\leq \frac{1}{K(n,2){{}^{2}}(\stackunder{M}{Sup}f)^{\frac{n-2}{n}}}%
(A_{t}^{1}+A_{t}^{2}) 
\]
where
\[
A_{t}^{1}=(\int_{B(0,\delta )}\overline{f}_{t}\overline{u}_{t}^{2^{*}}dv_{\,%
\mathbf{g}\,_{t}})^{\frac{2}{n}}(\int_{B(0,\delta )}\overline{f}_{t}(\eta 
\overline{u}_{t})^{2^{*}}dx)^{\frac{n-2}{n}}\,-(Supf.\int_{B(0,\delta
)}(\eta \overline{u}_{t})^{2^{*}}dx)^{\frac{n-2}{n}} 
\]
and
$$
A_{t}^{2}=\frac{2(\int_{B(0,\delta )}\overline{f}_{t}(\eta \overline{u}%
_{t})^{2^{*}}dx)^{\frac{n-2}{n}}}{n(\int_{B(0,\delta )}\overline{f}_{t}%
\overline{u}_{t}^{2^{*}}dv_{\,\mathbf{g}\,_{t}})^{\frac{n-2}{n}}}
\{\frac{1}{6}Ric(x_{t})_{ij}\int_{B(0,\delta )}x^{i}x^{j}\overline{f}_{t}%
\overline{u}_{t}^{2^{*}}dv_{\,\mathbf{g}\,_{t}}
+C\int_{B(0,\delta )}\left|
x\right| ^{3}\overline{u}_{t}^{2^{*}}dv_{\,\mathbf{g}\,_{t}}\}(1+\varepsilon
_{\delta }) 
$$
as $\{S_{t}\}\rightarrow 0$ when $\delta \rightarrow 0\,$uniformly in $t$.  $A_{t}^{2}$ will give, by 
developing the metric, $S_{\,\mathbf{g}%
\,}(x_{0})\,$while $A_{t}^{1}$ will give, by developing $f$, $-\frac{(n-2)(n-4)}{8(n-1)}%
\frac{\bigtriangleup _{\mathbf{g}}f(x_{0})}{f(x_{0})}$.

Note that for any $\alpha \in H_{1}^{2}(B(x_{0},2\delta )):$%
\[
\stackunder{t\rightarrow 1}{\lim }\frac{\int_{B(x_{t},\delta )}\alpha dx}{%
\int_{B(x_{t},\delta )}\alpha dv_{\,\mathbf{g}\,_{t}}}=1+O(\delta
^{2})=1+\varepsilon _{\delta } 
\]

We start by studying $A_{t}^{2}$:

1/: We have $\stackunder{t\rightarrow 1}{\lim }\frac{(\int_{B(0,\delta )}%
\overline{f}_{t}(\eta \overline{u}_{t})^{2^{*}}dx)^{\frac{n-2}{n}}}{%
(\int_{B(0,\delta )}\overline{f}_{t}\overline{u}_{t}^{2^{*}}dv_{\,\mathbf{g}%
\,_{t}})^{\frac{n-2}{n}}}=1+\varepsilon _{\delta }$

2/: Using the weak estimates (proposition 5), $\left| x\right| ^{2}%
\overline{u}_{t}^{2^{*}}\leq c\overline{u}_{t}^{2}$, from where we get: 
\[
\frac{\int_{B(0,\delta )}\left| x\right| ^{3}\overline{u}_{t}^{2^{*}}dv_{\,%
\mathbf{g}\,_{t}}}{\int_{B(0,\delta )}\overline{u}_{t}^{2}dv_{\,\mathbf{g}%
\,_{t}}}\leq C.\varepsilon _{\delta }\,\,. 
\]

3/: Using the blow-up formula's we write: for all $R>0:$%
\begin{eqnarray*}
\int_{B(0,\delta )}x^{i}x^{j}\overline{f}_{t}\overline{u}_{t}^{2^{*}}dv_{\,%
\mathbf{g}\,_{t}} && =\int_{B(0,R\mu _{t})}x^{i}x^{j}\overline{f}_{t}%
\overline{u}_{t}^{2^{*}}dv_{\,\mathbf{g}\,_{t}}+\int_{B(0,\delta )\backslash
B(0,R\mu _{t})}x^{i}x^{j}\overline{f}_{t}\overline{u}_{t}^{2^{*}}dv_{\,%
\mathbf{g}\,_{t}} \\ 
&& =\mu _{t}^{2}\int_{B(0,R)}x^{i}x^{j}\widetilde{f_{t}}\widetilde{u}%
_{t}^{2^{*}}dv_{\widetilde{\,\mathbf{g}\,}_{t}}+\mu _{t}^{2}\int_{B(0,\delta
\mu _{t}^{-1})\backslash B(0,R)}x^{i}x^{j}\widetilde{f_{t}}\widetilde{u}%
_{t}^{2^{*}}dv_{\widetilde{\,\mathbf{g}\,}_{t}}
\end{eqnarray*}
and
\[
\int_{B(0,\delta )}\overline{u}_{t}^{2}dv_{\mathbf{g}_{t}}=\mu
_{t}^{2}\int_{B(0,\delta \mu _{t}^{-1})}\widetilde{u}_{t}^{2}dv_{\widetilde{%
\mathbf{g}}_{t}}\,\,. 
\]
Using the weak estimates again, we get : 
\[
\int_{B(0,\delta \mu _{t}^{-1})\backslash B(0,R)}x^{i}x^{j}\widetilde{f_{t}}%
\widetilde{u}_{t}^{2^{*}}dv_{\widetilde{\mathbf{g}}_{t}}\leq \varepsilon
_{R}.\int_{B(0,\delta \mu _{t}^{-1})\backslash B(0,R)}\widetilde{u}%
_{t}^{2}dv_{\widetilde{\mathbf{g}}_{t}} 
\]
thus:
\[
\frac{\int_{B(0,\delta \mu _{t}^{-1})\backslash B(0,R)}x^{i}x^{j}\widetilde{%
f_{t}}\widetilde{u}_{t}^{2^{*}}dv_{\widetilde{\mathbf{g}}_{t}}}{%
\int_{B(0,\delta \mu _{t}^{-1})}\widetilde{u}_{t}^{2}dv_{\widetilde{\mathbf{g%
}}_{t}}}\leq \varepsilon _{R} 
\]
where $\varepsilon _{R}\rightarrow 0$ when $R\rightarrow +\infty $.
Now, if $i\neq j:$%
\[
\stackunder{t\rightarrow 1}{\overline{\lim }}\frac{|\int_{B(0,R)}x^{i}x^{j}%
\widetilde{f_{t}}\widetilde{u}_{t}^{2^{*}}dv_{\widetilde{\mathbf{g}}_{t}}|}{%
\int_{B(0,\delta \mu _{t}^{-1})}\widetilde{u}_{t}^{2}dv_{\widetilde{\mathbf{g%
}}_{t}}}\leq \stackunder{t\rightarrow 1}{\overline{\lim }}\frac{%
|\int_{B(0,R)}x^{i}x^{j}\widetilde{f_{t}}\widetilde{u}_{t}^{2^{*}}dv_{%
\widetilde{\mathbf{g}}_{t}}|}{\int_{B(0,R)}\widetilde{u}_{t}^{2}dv_{%
\widetilde{\mathbf{g}}_{t}}}=0 
\]
because
\[
\widetilde{u}_{t}\rightarrow \widetilde{u}=(1+\frac{f(x_{0})^{\frac{2}{n}}}{%
K(n,2)^{2}n(n-2)}\left| x\right| ^{2})^{-\frac{n-2}{2}}\,\,\,in\,\,%
\,C^{0}(B(0,R)) 
\]
and $\widetilde{u}$ is radial (see subsection 3.2).

If $i=j:$%
\[
\frac{\int_{B(0,R)}x^{i}x^{i}\widetilde{f_{t}}\widetilde{u}_{t}^{2^{*}}dv_{%
\widetilde{\mathbf{g}}_{t}}}{\int_{B(0,\delta \mu _{t}^{-1})}\widetilde{u}%
_{t}^{2}dv_{\widetilde{\mathbf{g}}_{t}}}=\frac{\int_{B(0,R)}(x^{i}){{}^{2}}%
\widetilde{f_{t}}\widetilde{u}_{t}^{2^{*}}dv_{\widetilde{\mathbf{g}}_{t}}}{%
\int_{B(0,R)}\widetilde{u}_{t}^{2}dv_{\widetilde{\mathbf{g}}_{t}}}.\frac{%
\int_{B(0,R)}\widetilde{u}_{t}^{2}dv_{\widetilde{\mathbf{g}}_{t}}}{%
\int_{B(0,\delta \mu _{t}^{-1})}\widetilde{u}_{t}^{2}dv_{\widetilde{\mathbf{g%
}}_{t}}} 
\]
But as soon as $n>4$, using strong $L{{}^{2}}$-concentration 
(proposition 9), we obtain: 
\[
\stackunder{R\rightarrow \infty }{\lim }\,\stackunder{t\rightarrow 1}{%
\overline{\lim }}\frac{\int_{B(0,R)}\widetilde{u}_{t}^{2}dv_{\widetilde{%
\mathbf{g}}_{t}}}{\int_{B(0,\delta \mu _{t}^{-1})}\widetilde{u}_{t}^{2}dv_{%
\widetilde{\mathbf{g}}_{t}}}=1 
\]
therefore

\[
\stackunder{R\rightarrow \infty }{\lim }\,\stackunder{t\rightarrow 1}{%
\overline{\lim }}\frac{\int_{B(0,R)}x^{i}x^{i}\widetilde{f_{t}}\widetilde{u}%
_{t}^{2^{*}}dv_{\widetilde{\mathbf{g}}_{t}}}{\int_{B(0,\delta \mu _{t}^{-1})}%
\widetilde{u}_{t}^{2}dv_{\widetilde{\mathbf{g}}_{t}}}=f(x_{0})\frac{\int_{%
\Bbb{R}^{n}}(x^{i})^{2}.\widetilde{u}^{2}dx}{\int_{\Bbb{R}^{n}}\widetilde{u}%
^{2}dx}=f(x_{0})^{\frac{n-2}{n}}K(n,2)^{2}\frac{n-4}{4(n-1)} 
\]
and thus
\[
\stackunder{t\rightarrow 1}{\overline{\lim }}\frac{1}{f(x_{0})^{\frac{n-2}{n}%
}K(n,2)^{2}}\frac{A_{t}^{2}}{\int_{B(0,\delta )}\overline{u}_{t}^{2}dv_{%
\mathbf{g}_{t}}}=\frac{n-4}{12(n-1)}S_{\mathbf{g}}(x_{0})+\varepsilon
_{\delta } 
\]
which, with $\stackunder{t\rightarrow 1}{\overline{\lim }}\frac{B_{t}}{%
\int_{B(0,\delta )}\overline{u}_{t}^{2}dx}=\frac{1}{6}S_{\mathbf{g}%
}(x_{0})+\varepsilon _{\delta }$ gives 
\[
\stackunder{t\rightarrow 1}{\overline{\lim }}(\frac{1}{f(x_{0})^{\frac{n-2}{n%
}}K(n,2)^{2}}\frac{A_{t}^{2}}{\int_{B(0,\delta )}\overline{u}_{t}^{2}dv_{%
\mathbf{g}_{t}}}+\frac{B_{t}}{\int_{B(0,\delta )}\overline{u}_{t}^{2}dx})=%
\frac{n-2}{4(n-1)}S_{\mathbf{g}}(x_{0})+\varepsilon _{\delta } 
\]

If $n=4$ we write: 
\[
\stackunder{R\rightarrow \infty }{\lim }\stackunder{t\rightarrow 1}{%
\overline{\lim }}\frac{\int_{B(0,R)}x^{i}x^{i}\widetilde{f_{t}}\widetilde{u}%
_{t}^{2^{*}}dv_{\widetilde{\mathbf{g}}_{t}}}{\int_{B(0,\delta \mu _{t}^{-1})}%
\widetilde{u}_{t}^{2}dv_{\widetilde{\mathbf{g}}_{t}}}\leq f(x_{0})^{\frac{n-2%
}{n}}K(n,2)^{2}\frac{n-4}{4(n-1)} 
\]
and we get the conclusion by distinguishing two cases, $S_{\mathbf{g}}(x_{0})<0$ or $S_{%
\mathbf{g}}(x_{0})\geq 0$, the proof being finished as $%
\frac{\bigtriangleup _{\mathbf{g}}f(x_{0})}{f(x_{0})}$ does not appear in dimension 4 (see the end of the proof).

Let us now consider $A^1_{t}$.
\[
A_{t}^{1}=(\int_{B(0,\delta )}\overline{f}_{t}\overline{u}_{t}^{2^{*}}dv_{%
\mathbf{g}_{t}})^{\frac{2}{n}}(\int_{B(0,\delta )}\overline{f}_{t}(\eta 
\overline{u}_{t})^{2^{*}}dx)^{\frac{n-2}{n}}\,-(Supf.\int_{B(0,\delta
)}(\eta \overline{u}_{t})^{2^{*}}dx)^{\frac{n-2}{n}}\,\,. 
\]
We write $f_{t}=f(x_{0})+g_{t}$. Remembering that $f(x_{0})=Sup f$, we have $g_{t}(x_{0})=0$ and 
$g_{t}\leq0$. Using $(1+x)^\alpha \leq 1+\alpha x$ for $0<\alpha \leq1$:
\[
(\int_{B(0,\delta )}\overline{f}_{t}(\eta \overline{u}_{t})^{2^{*}}dx)^{%
\frac{n-2}{n}}\leq (\int_{B(0,\delta )}f(x_{0})(\eta \overline{u}%
_{t})^{2^{*}}dx)^{\frac{n-2}{n}}+\frac{n-2}{n}\frac{\int_{B(0,\delta )}\overline{g}_{t}(\eta \overline{u}_{t})^{2^{*}}dx}
{(\int_{B(0,\delta)}f(x_{0})(\eta \overline{u}_{t})^{2^{*}}dx)^{\frac{2}{n}}}%
\]
where $\overline{g}_{t}$ is $g_{t}$ in the exponential chart in $x_{t}$. 
We now use the theorem of Druet and Robert to write in $B(0,\delta)$: 
$$\overline{u}_{t}\geq (1-\ve_{\delta})\overline{B}_{t},$$
where $\overline{B}_{t}$ is $B_{t}$ in the exponential chart in $x_{t}$.
Because $g_{t}\leq0$, we have:
$$\int_{B(0,\delta )}\overline{g}_{t}(\eta \overline{u}_{t})^{2^{*}}dx\leq
(1-\ve_{\delta})\int_{B(0,\delta )}\overline{g}_{t}(\eta \overline{B}_{t})^{2^{*}}dx.$$
Combining this with the expansion above and the fact that 
$\int_{B(0,\delta )}\overline{f}_{t}\overline{u}_{t}^{2^{*}}dv_{\mathbf{g}_{t}}\leq1$, we obtain:
$$A^1_{t}\leq(1-\ve_{\delta})\frac{n-2}{n}\frac{\int_{B(0,\delta )}\overline{g}_{t}(\eta \overline{B}_{t})^{2^{*}}dx}
{(\int_{B(0,\delta)}f(x_{0})(\eta \overline{u}_{t})^{2^{*}}dx)^{\frac{2}{n}}}$$
We now expand $g_{t}$ noting that $\partial _{i}\overline{g}_{t}=\partial _{i}\overline{f}_{t}$ and 
$\partial _{i}\overline{g}_{t}=\partial _{i}\overline{f}_{t}$.
\[
\overline{g}_{t}(x)\leq g_{t}(x_{t})+x^{i}\partial _{i}\overline{f}_{t}(x_{t})+\frac{1}{2}\partial _{kl}\overline{f}%
_{t}(x_{t}).x^{k}x^{l}+c\left|x\right| ^{3}
\]
Thus
\begin{eqnarray*}
\int_{B(0,\delta )}\overline{g}_{t}(\eta \overline{B}_{t})^{2^{*}}dx\leq
&& g_{t}(x_{t})\int_{B(0,\delta)}(\eta \overline{B}_{t})^{2^{*}}dx\\
&&+\partial _{i}\overline{f}_{t}(x_{t})\int_{B(0,\delta)}x^{i}(\eta \overline{B}_{t})^{2^{*}}dx\\
&&+\frac{1}{2}\partial _{kl}\overline{f}_{t}(x_{t})\int_{B(0,\delta)}x^{k}x^{l}(\eta \overline{B}_{t})^{2^{*}}dx\\
&&+C\int_{B(0,\delta )}\left| x\right| ^{3}(\eta \overline{B}_{t})^{2^{*}}dx
\end{eqnarray*}
Now, first $g_{t}(x_{t})\leq0$, and second, and this is the main point for which we need the theorem of Druet and Robert (see the reason 
at the beginning of the next section), as $\overline{B}_{t}$ is radial, we have
$$\partial _{i}\overline{f}_{t}(x_{t})\int_{B(0,\delta)}x^{i}(\eta \overline{B}_{t})^{2^{*}}dx=0.$$
Therefore, introducing all this in the last inequality for $A^1_{t}$, we have
\[
\stackunder{t\rightarrow 1}{\overline{\lim }}\frac{A_{t}^{1}}{%
\int_{B(0,\delta )}\overline{u}_{t}^{2}dv_{\mathbf{g}_{t}}}\leq 
\]
$$\frac{n-2}{n}(1-\varepsilon _{\delta })\stackunder{t\rightarrow 1}{\overline{\lim }}
\frac{\frac{1}{2}\partial _{kl}\overline{f}_{t}(x_{t})\int_{B(0,\delta
)}x^{k}x^{l}(\eta \overline{B}_{t})^{2^{*}}dv_{%
\mathbf{g}_{t}}+C\int_{B(0,\delta )}\left| x\right| ^{3}(\eta 
\overline{B}_{t})^{2^{*}}dv_{\mathbf{g}_{t}}}{\int_{B(0,\delta )}\overline{u}%
_{t}^{2}dv_{\mathbf{g}_{t}}}$$
where we have replaced $dx$ by $dv_{\mathbf{g}_{t}}$ using the remark made at the beginning of the study of $A^2_{t}$.

Now, as for $A_{t}^{2}$, we write:
\begin{eqnarray*}
\stackunder{t\rightarrow 1}{\overline{\lim }}\frac{\int_{B(0,\delta
)}x^{k}x^{l}(\eta \overline{B}_{t})^{2^{*}}dv_{\mathbf{g}_{t}}}{%
\int_{B(0,\delta )}\overline{u}_{t}^{2}dv_{\mathbf{g}_{t}}} &=&f(x_{0})^{%
\frac{-2}{n}}K(n,2)^{2}\frac{n-4}{4(n-1)}\text{ if }k=l \\
&=&0\text{ if }k\neq l
\end{eqnarray*}
and therefore
\[
\frac{1}{K(n,2)^{2}f(x_{0})^{\frac{n-2}{n}}}\frac{n-2}{n}\stackunder{%
t\rightarrow 1}{\overline{\lim }}\frac{\frac{1}{2}\partial _{kl}\overline{f}_{t}(x_{t})\int_{B(0,\delta
)}x^{k}x^{l}(\eta \overline{B}_{t})^{2^{*}}dv_{\mathbf{g}_{t}}}{%
\int_{B(0,\delta )}\overline{u}_{t}^{2}dv_{\mathbf{g}_{t}}}= 
\]
\[
=\frac{1}{f(x_{0})}\frac{(n-2)(n-4)}{4(n-1)}\sum_{l}\frac{1}{2}\partial_{ll}\overline{f}_{1}(0)
\]
\[
=-\frac{(n-2)(n-4)}{8(n-1)}\frac{\bigtriangleup _{\mathbf{g}}f(x_{0})}{f(x_{0})}
\]
as $\bigtriangleup _{\mathbf{g}}f(x_{0})=-\sum_{l}\partial _{ll}\overline{f}%
_{1}(0)$ in the exponential chart in $x_{0}$.
Also
\[
\frac{\int_{B(0,\delta )}\left| x\right| ^{3}\overline{u}_{t}^{2^{*}}dv_{\,%
\mathbf{g}\,_{t}}}{\int_{B(0,\delta )}\overline{u}_{t}^{2}dv_{\,\mathbf{g}%
\,_{t}}}\leq C.\varepsilon _{\delta }\,\,. 
\]
Thus, we have proved that dividing inequality (18) by $\int_{B(0,\delta )}\overline{u}_{t}^{2}dv_{\mathbf{g}_{t}}$ 
and letting $t$ go to 1, we get
\[
h(x_{0})+\varepsilon _{\delta }\leq \frac{n-2}{4(n-1)}S_{\mathbf{g}}(x_{0})-%
\frac{(n-2)(n-4)}{8(n-1)}\frac{\bigtriangleup _{\mathbf{g}}f(x_{0})}{f(x_{0})%
}+\varepsilon _{\delta }\,\,. 
\]
Letting $\delta $ tend to 0: 
\[
h(x_{0})\leq \frac{n-2}{4(n-1)}S_{\mathbf{g}}(x_{0})-\frac{(n-2)(n-4)}{8(n-1)%
}\frac{\bigtriangleup _{\mathbf{g}}f(x_{0})}{f(x_{0})} 
\]
which contradict our hypothesis: 
\[
h(x_{0})>\frac{n-2}{4(n-1)}S_{\mathbf{g}}(x_{0})-\frac{(n-2)(n-4)}{8(n-1)}%
\frac{\bigtriangleup _{\mathbf{g}}f(x_{0})}{f(x_{0})} 
\]
when $x_{0}$ is a point of maximum of $f$. This prove that $u\not\equiv 0$, and therefore $u_{t}\rightarrow u>0$, a minimizing solution for $\E$,
and thus the weakly critical function $h$ is in fact critical.

\subsection{Alternate proof, proof of the fundamental estimate}
As we saw in the last part of the proof, the difficulty introduced by the presence of the function $f$ is to control 
the first derivatives of $f$, $\partial _{i}f(x_{t})$, as blow-up gives 
\[
\int_{B(0,\delta )}\partial _{i}f(x_{t})x^{i}\overline{u}_{t}^{2^{*}}dv_{%
\mathbf{g}_{t}}=\mu _{t}\int_{B(0,\delta \mu _{t}^{-1})}\partial
_{i}f(x_{t})x^{i}\widetilde{u}_{t}^{2^{*}}dv_{\widetilde{\mathbf{g}}_{t}} 
\]
to be divided by
\[
\mu _{t}^{2}\int_{B(0,\delta \mu _{t}^{-1})}\widetilde{u}_{t}^{2}dv_{%
\widetilde{\mathbf{g}}_{t}}\,\,, 
\]
and it would be necessary to control $\frac{\partial
_{i}f(x_{t})}{\mu _{t}}$, which seems to be difficult. But thanks to the theorem of Druet and Robert, we can replace 
$u(t)$ by $B(t)$ near $x_{t}$, and after blow-up 
$$\mu _{t}\int_{B(0,\delta \mu _{t}^{-1})}\partial
_{i}f(x_{t})x^{i}\widetilde{B}_{t}^{2^{*}}dv_{\widetilde{\mathbf{g}}_{t}} =0$$
as $\widetilde{B}_{t}$ is radial. Of course, the proof is then short, but the proof of the theorem of Druet and Robert 
is quite involved, even though the strong estimates (proposition 8) is the first step. 

The other way to get over the problem of the first derivatives of $f$ is to expand $f$ in $x_{0}$ as then 
 $\partial _{i}f(x_{0})=0$ because $x_{0}$ is a point of maximum of $f$. 
But then, one has to transpose the weak and strong estimates from $x_{t}$ to $x_{0}$,
which, as we said in the section about concentration phenomenom, requires to prove the following estimate: 
\[
\frac{d_{g}(x_{t},x_{0})}{\mu _{t}}\leq C\,\,. 
\]
As we said, this estimate is important and of independent interest, as it gives a complete description of the sequence 
$(u_{t})$. This is why we give this alternate proof of theorem 1, even though it requires an additional hypothesis.
This proof, which gives at the same time the proof of theorem 1 and of the estimate, is, we think, interesting, and is available 
directly after proposition 9, i.e it does not require the theorem of Druet and Robert.

We now make the hypothesis that the hessian of $f$ is nondegenerate at its points of maximum. 
We also suppose now that $dim M \geq 5$, even though our proof gives theorem 1 in dimension 4.

Let us note $x_{0}(t)=\exp _{x_{t}}^{-1}(x_{0})=(x_{0}^{1}(t),...,x_{0}^{n}(t))$, 
which is possible as soon as $t$ is close enough to 1 for a fixed radius $\delta$. Then $x_{0}(t)\rightarrow 0$
 when $t\rightarrow 1$. The point $x_{0}(t)$ is a locally strict maximum of $\overline{f}_{t}$. 
 We will let $\delta$ go to 0 at the end of the reasoning, after having taken the limit when $t\rightarrow 1$. 

The expansion of $\overline{f}_{t}$ in $x_{0}(t)$ gives: 
\[
\overline{f}_{t}(x)\leq f(x_{0})+\frac{1}{2}\partial _{kl}\overline{f}%
_{t}(x_{0}(t)).(x^{k}-x_{0}^{k}(t))(x^{l}-x_{0}^{l}(t))+c\left|
x-x_{0}(t)\right| ^{3}:=f(x_{0})+T_{t} 
\]
($T_{t}$ like Taylor) where ($\partial _{kl}\overline{f}_{t}(x_{0})$) is a negative definite matrix 
(we shall write $<0$). $c,C$ will always be constants independent of $t$ and $\delta $.
Remember that
\[
A_{t}^{1}=(\int_{B(0,\delta )}\overline{f}_{t}\overline{u}_{t}^{2^{*}}dv_{%
\mathbf{g}_{t}})^{\frac{2}{n}}(\int_{B(0,\delta )}\overline{f}_{t}(\eta 
\overline{u}_{t})^{2^{*}}dx)^{\frac{n-2}{n}}\,-(Supf.\int_{B(0,\delta
)}(\eta \overline{u}_{t})^{2^{*}}dx)^{\frac{n-2}{n}}\,\,. 
\]
Introducing the expansion of $\overline{f}_{t}$ in $x_{0}(t)$, 
and using again the fact that $(1+x)^\alpha \leq 1+\alpha x$ for $0<\alpha \leq1$, we get: 
\[
(\int_{B(0,\delta )}\overline{f}_{t}(\eta \overline{u}_{t})^{2^{*}}dx)^{%
\frac{n-2}{n}}\leq (\int_{B(0,\delta )}f(x_{0})(\eta \overline{u}%
_{t})^{2^{*}}dx)^{\frac{n-2}{n}}+\frac{\frac{n-2}{n}}{(\int_{B(0,\delta
)}f(x_{0})(\eta \overline{u}_{t})^{2^{*}}dx)^{\frac{2}{n}}}%
\{F_{t}\} 
\]
where
\[
\{F_{t}\}=\frac{1}{2}\partial _{kl}\overline{f}_{t}(x_{0}(t))\int_{B(0,%
\delta )}(x^{k}-x_{0}^{k}(t))(x^{l}-x_{0}^{l}(t))(\eta \overline{u}%
_{t})^{2^{*}}dx+C\int_{B(0,\delta )}\left| x-x_{0}(t)\right| ^{3}(\eta 
\overline{u}_{t})^{2^{*}}dx 
\]
from where, remembering that $\stackunder{M}{Sup}f=f(x_{0})$ and that $%
\int_{B(0,\delta )}\overline{f}_{t}\overline{u}_{t}^{2^{*}}dv_{\mathbf{g}%
_{t}}\leq 1$: 
\begin{equation}
A_{t}^{1}\leq \frac{n-2}{n}\frac{(\int_{B(0,\delta )}\overline{f}_{t}%
\overline{u}_{t}^{2^{*}}dv_{\mathbf{g}_{t}})^{\frac{2}{n}}}{%
(\int_{B(0,\delta )}f(x_{0})(\eta \overline{u}_{t})^{2^{*}}dx)^{\frac{2}{n}}}%
\{F_{t}\}
\end{equation}
Therefore, we obtain:
\[
\stackunder{t\rightarrow 1}{\overline{\lim }}\frac{A_{t}^{1}}{%
\int_{B(0,\delta )}\overline{u}_{t}^{2}dv_{\mathbf{g}_{t}}}\leq 
\]
$\frac{n-2}{n}(1+\varepsilon _{\delta })\stackunder{t\rightarrow 1}{\overline{\lim }}
\frac{\frac{1}{2}\partial _{kl}\overline{f}_{t}(x_{0})\int_{B(0,\delta
)}(x^{k}-x_{0}^{k}(t))(x^{l}-x_{0}^{l}(t))(\eta \overline{u}_{t})^{2^{*}}dv_{%
\mathbf{g}_{t}}+C\int_{B(0,\delta )}\left| x-x_{0}(t)\right| ^{3}(\eta 
\overline{u}_{t})^{2^{*}}dv_{\mathbf{g}_{t}}}{\int_{B(0,\delta )}\overline{u}%
_{t}^{2}dv_{\mathbf{g}_{t}}}$
\\

where we write $\partial _{kl}\overline{f}_{t}(x_{0})$ for $\partial
_{kl}\overline{f}_{t}(x_{0}(t)).$ Considering the expansion
\[
\overline{f}_{t}(x)\leq f(x_{0})+\frac{1}{2}\partial _{kl}\overline{f}%
_{t}(x_{0}(t)).(x^{k}-x_{0}^{k}(t))(x^{l}-x_{0}^{l}(t))+c\left|
x-x_{0}(t)\right| ^{3},
\]
note that by the regularity of $\exp _{x_{t}}^{-1}\circ \exp _{x_{0}}$ with respect to all the variables, 
we can suppose that $c$ is independent of $t$. Moreover: 
\begin{eqnarray*}
c\left| x-x_{0}(t)\right| ^{3} & &\leq c^{\prime }\left| x-x_{0}(t)\right|
\sum_{k}(x^{k}-x_{0}^{k}(t))^{2} \\ 
& &\leq 2\delta c^{\prime }\sum_{k}(x^{k}-x_{0}^{k}(t))^{2}
\end{eqnarray*}
where we remind that $\delta $ is the radius of the ball of integration. 
We can then write:
\[
\overline{f}_{t}(x)\leq f(x_{0})+(\frac{1}{2}\partial _{kl}\overline{f}%
_{t}(x_{0}(t))+\delta C_{kl})(x^{k}-x_{0}^{k}(t))(x^{l}-x_{0}^{l}(t)) 
\]
where $C_{kl}=c\delta _{kl}=c$ if $k=l$ and $C_{kl}=0$ if $k\neq l$ ($\delta
_{kl}$ is the Kr\"{o}necker symbol) is independent of $t$.

We introduce one more notation:
\[
D_{kl}(t,\delta )=\frac{1}{2}\partial _{kl}\overline{f}_{t}(x_{0}(t))+\delta
C_{kl}\,\,. 
\]
Then:

1/: $\stackunder{\delta \rightarrow 0}{\lim }\stackunder{t\rightarrow 1}{%
\lim }D_{kl}(t,\delta )=\frac{1}{2}\partial _{kl}\overline{f}_{1}(x_{0}(1))$
where $\overline{f}_{1}=f\circ \exp _{x_{0}}^{-1}$ and \\
$x_{0}(1)=0=\exp_{x_{0}}^{-1}(x_{0})$.

2/: for any $\delta $ small enough and for all $t$ close to 1, $D_{kl}(t,\delta)$ is still negative definite.

$D_{kl}(t,\delta )$ is the hessian of $f$ in $x_{0}(t)$ perturbated on its diagonal by the third order terms. 
It is for the second point that we need the hypothesis that the hessian of $f$ is non degenerate. Thus
\[
\frac{1}{2}\partial _{kl}\overline{f}_{t}(x_{0})\int_{B(0,\delta
)}(x^{k}-x_{0}^{k}(t))(x^{l}-x_{0}^{l}(t))(\eta \overline{u}_{t})^{2^{*}}dv_{%
\mathbf{g}_{t}}+C\int_{B(0,\delta )}\left| x-x_{0}(t)\right| ^{3}(\eta 
\overline{u}_{t})^{2^{*}}dv_{\mathbf{g}_{t}}\leq 
\]
\[
D_{kl}(t,\delta )\int_{B(0,\delta
)}(x^{k}-x_{0}^{k}(t))(x^{l}-x_{0}^{l}(t))(\eta \overline{u}_{t})^{2^{*}}dv_{%
\mathbf{g}_{t}}\,\,. 
\]
Let
\[
\{F_{t}^{\prime }\}=D_{kl}(t,\delta )\int_{B(0,\delta
)}(x^{k}-x_{0}^{k}(t))(x^{l}-x_{0}^{l}(t))(\eta \overline{u}_{t})^{2^{*}}dv_{%
\mathbf{g}_{t}} 
\]
We have
\[
\stackunder{t\rightarrow 1}{\overline{\lim }}\frac{A_{t}^{1}}{%
\int_{B(0,\delta )}v_{t}^{2}dv_{\mathbf{g}_{t}}}\leq \frac{n-2}{n}%
\stackunder{t\rightarrow 1}{\overline{\lim }}\frac{D_{kl}(t,\delta
)\int_{B(0,\delta )}(x^{k}-x_{0}^{k}(t))(x^{l}-x_{0}^{l}(t))(\eta \overline{u%
}_{t})^{2^{*}}dv_{\mathbf{g}_{t}}}{\int_{B(0,\delta )}\overline{u}%
_{t}^{2}dv_{\mathbf{g}_{t}}}(1+\varepsilon _{\delta }) 
\]

In the expansion of $D_{kl}(t,\delta
)(x^{k}-x_{0}^{k}(t))(x^{l}-x_{0}^{l}(t))$, we are interested by the first term, 
i.e $D_{kl}(t,\delta )x^{k}x^{l}$ (look back how we obtained $%
S_{g}(x_{0})$ in $A_{t}^{2}$), and we are going to show that the other terms can be neglected. 
The idea is to reorganize the expansion of $%
\{F_{t}^{\prime }\}$ and to use the fact $D_{kl}(t,\delta )$ is a negative bilinear form: 
\begin{eqnarray*}
\{F_{t}^{\prime }\}=&&D_{kl}(t,\delta )\int_{B(0,\delta
)}x^{k}x^{l}(\eta \overline{u}_{t})^{2^{*}}dv_{\mathbf{g}_{t}}+D_{kl}(t,%
\delta )x_{0}^{k}(t)x_{0}^{l}(t)\int_{B(0,\delta )}(\eta \overline{u}%
_{t})^{2^{*}}dv_{\mathbf{g}_{t}} \\ 
& &-D_{kl}(t,\delta )\int_{B(0,\delta
)}(x^{k}x_{0}^{l}(t)+x^{l}x_{0}^{k}(t))(\eta \overline{u}_{t})^{2^{*}}dv_{%
\mathbf{g}_{t}}\,\,.
\end{eqnarray*}
We rewrite the two last terms (suppressing some $\delta $ et $%
t $ and all integral being taken with respect to $dv_{\mathbf{g}_{t}}$):

\[
D_{kl}.x_{0}^{k}x_{0}^{l}\int_{B(0,\delta )}(\eta \overline{u}%
_{t})^{2^{*}}dv_{\mathbf{g}_{t}}-D_{kl}\int_{B(0,\delta
)}(x^{k}x_{0}^{l}+x^{l}x_{0}^{k})(\eta \overline{u}_{t})^{2^{*}}dv_{\mathbf{g%
}_{t}}= 
\]

\[
D_{kl}\Big[ x_{0}^{k}x_{0}^{l}\int_{B(0,\delta )}(\eta \overline{u}%
_{t})^{2^{*}}-x_{0}^{l}\int_{B(0,\delta )}x^{k}(\eta \overline{u}%
_{t})^{2^{*}}-x_{0}^{k}\int_{B(0,\delta )}x^{l}(\eta \overline{u}%
_{t})^{2^{*}}\Big]= 
\]
\[
D_{kl}\Big[x_{0}^{k}(\int_{B(0,\delta )}(\eta \overline{u}_{t})^{2^{*}})^{\frac{1%
}{2}}.x_{0}^{l}(\int_{B(0,\delta )}(\eta \overline{u}_{t})^{2^{*}})^{\frac{1%
}{2}} 
\]
\[
-x_{0}^{l}(\int_{B(0,\delta )}(\eta \overline{u}_{t})^{2^{*}})^{\frac{1}{2}}%
\frac{\int_{B(0,\delta )}x^{k}(\eta \overline{u}_{t})^{2^{*}}}{%
(\int_{B(0,\delta )}(\eta \overline{u}_{t})^{2^{*}})^{\frac{1}{2}}}%
-x_{0}^{k}(\int_{B(0,\delta )}(\eta \overline{u}_{t})^{2^{*}})^{\frac{1}{2}}%
\frac{\int_{B(0,\delta )}x^{l}(\eta \overline{u}_{t})^{2^{*}}}{%
(\int_{B(0,\delta )}(\eta \overline{u}_{t})^{2^{*}})^{\frac{1}{2}}}\Big]\,\,. 
\]
Thus, setting (sorry): 
\[
\varepsilon ^{k}(t)=\int_{B(0,\delta )}x^{k}(\eta \overline{u}%
_{t})^{2^{*}}dv_{\mathbf{g}_{t}} 
\]
\[
z_{t}=(\int_{B(0,\delta )}(\eta \overline{u}_{t})^{2^{*}}dv_{\mathbf{g}%
_{t}})^{\frac{1}{2}} 
\]
the expression above becomes:
\begin{center}
\[
D_{kl}.x_{0}^{k}x_{0}^{l}\int_{B(0,\delta )}(\eta \overline{u}%
_{t})^{2^{*}}dv_{\mathbf{g}_{t}}-D_{kl}\int_{B(0,\delta
)}(x^{k}x_{0}^{l}+x^{l}x_{0}^{k})(\eta \overline{u}_{t})^{2^{*}}dv_{\mathbf{g%
}_{t}}= 
\]
\[
=D_{kl}\Big[x_{0}^{k}(t).z_{t}.x_{0}^{l}(t).z_{t}-x_{0}^{l}(t).z_{t}.\frac{%
\varepsilon ^{k}(t)}{z_{t}}-x_{0}^{k}(t).z_{t}.\frac{\varepsilon ^{l}(t)}{z_{t}}\Big] 
\]
\end{center}
\[
=D_{kl}\Big[(x_{0}^{k}(t).z_{t}-\frac{\varepsilon ^{k}(t)}{z_{t}}%
)(x_{0}^{l}(t).z_{t}-\frac{\varepsilon ^{l}(t)}{z_{t}})-\frac{\varepsilon
^{k}(t)\varepsilon ^{l}(t)}{z_{t}^{2}}\Big] 
\]
By this method of reorganization of the hessian, we have obtained: 
\[
\frac{1}{2}\partial _{kl}\overline{f}_{t}(x_{0})\int_{B(0,\delta
)}(x^{k}-x_{0}^{k}(t))(x^{l}-x_{0}^{l}(t))(\eta \overline{u}_{t})^{2^{*}}dv_{%
\mathbf{g}_{t}}+C\int_{B(0,\delta )}\left| x-x_{0}(t)\right| ^{3}(\eta 
\overline{u}_{t})^{2^{*}}dv_{\mathbf{g}_{t}}\leq 
\]
\begin{center}
\[
D_{kl}(t,\delta )\int_{B(0,\delta )}x^{k}x^{l}(\eta \overline{u}%
_{t})^{2^{*}}dv_{\mathbf{g}_{t}}+D_{kl}(t,\delta )(x_{0}^{k}(t).z_{t}-\frac{%
\varepsilon ^{k}(t)}{z_{t}})(x_{0}^{l}(t).z_{t}-\frac{\varepsilon ^{l}(t)}{%
z_{t}})-D_{kl}(t,\delta )\frac{\varepsilon ^{k}(t)\varepsilon ^{l}(t)}{%
z_{t}^{2}} 
\]
\end{center}
\[
\leq D_{kl}(t,\delta )\int_{B(0,\delta )}x^{k}x^{l}(\eta \overline{u}%
_{t})^{2^{*}}dv_{\mathbf{g}_{t}}-D_{kl}(t,\delta )\frac{\varepsilon
^{k}(t)\varepsilon ^{l}(t)}{z_{t}^{2}} 
\]
because, and that is the fundamental point : 
\[
D_{kl}(t,\delta )\omega ^{k}\omega ^{l}\leq 0\text{ \thinspace }\forall
\omega =(\omega ^{1},...,\omega ^{n}) 
\]
which allows to suppress from the inequality $$D_{kl}(t,\delta
)(x_{0}^{k}(t).z_{t}-\frac{\varepsilon ^{k}(t)}{z_{t}})(x_{0}^{l}(t).z_{t}-%
\frac{\varepsilon ^{l}(t)}{z_{t}})$$ It is this term that will give us the estimate $\frac{d_{\mathbf{g}}(x_{t},x_{0})}{\mu _{t}}\leq C$ 
(see below).

We have therefore obtained:
\[
\stackunder{t\rightarrow 1}{\overline{\lim }}\frac{A_{t}^{1}}{%
\int_{B(0,\delta )}\overline{u}_{t}^{2}dv_{\mathbf{g}_{t}}}\leq \frac{n-2}{n}%
\stackunder{t\rightarrow 1}{\overline{\lim }}\frac{D_{kl}(t,\delta
)\int_{B(0,\delta )}x^{k}x^{l}(\eta \overline{u}_{t})^{2^{*}}dv_{\mathbf{g}%
_{t}}-D_{kl}(t,\delta )\frac{\varepsilon ^{k}(t)\varepsilon ^{l}(t)}{%
z_{t}^{2}}}{\int_{B(0,\delta )}\overline{u}_{t}^{2}dv_{\mathbf{g}_{t}}}%
(1+\varepsilon _{\delta }) 
\]
Now, as for $A_{t}^{2}$, we write:
\begin{eqnarray*}
\stackunder{t\rightarrow 1}{\overline{\lim }}\frac{\int_{B(0,\delta
)}x^{k}x^{l}(\eta \overline{u}_{t})^{2^{*}}dv_{\mathbf{g}_{t}}}{%
\int_{B(0,\delta )}\overline{u}_{t}^{2}dv_{\mathbf{g}_{t}}} &=&f(x_{0})^{%
\frac{-2}{n}}K(n,2)^{2}\frac{n-4}{4(n-1)}\text{ if }k=l \\
&=&0\text{ if }k\neq l
\end{eqnarray*}
and therefore
\[
\frac{1}{K(n,2)^{2}f(x_{0})^{\frac{n-2}{n}}}\frac{n-2}{n}\stackunder{%
t\rightarrow 1}{\overline{\lim }}\frac{D_{kl}(t,\delta )\int_{B(0,\delta
)}x^{k}x^{l}(\eta \overline{u}_{t})^{2^{*}}dv_{\mathbf{g}_{t}}}{%
\int_{B(0,\delta )}\overline{u}_{t}^{2}dv_{\mathbf{g}_{t}}}= 
\]
\[
=\frac{1}{f(x_{0})}\frac{(n-2)(n-4)}{4(n-1)}\sum_{l}(\frac{1}{2}\partial
_{ll}\overline{f}_{1}(0)+c_{ll}\delta ) 
\]
\[
=-\frac{(n-2)(n-4)}{8(n-1)}\frac{\bigtriangleup _{\mathbf{g}}f(x_{0})}{%
f(x_{0})}+\varepsilon _{\delta } 
\]
as $\bigtriangleup _{\mathbf{g}}f(x_{0})=-\sum_{l}\partial _{ll}\overline{f}%
_{1}(0)$ in the exponential chart in $x_{0}$.

At last, let us show that the residual term can be neglected. 
\[
\left| \varepsilon ^{k}(t)\varepsilon ^{l}(t)\right| \leq \frac{1}{2}%
(\varepsilon ^{k}(t)^{2}+\varepsilon ^{l}(t)^{2}) 
\]
But 
\[
\varepsilon ^{k}(t)^{2}=(\int_{B(0,\delta )}x^{k}(\eta \overline{u}%
_{t})^{2^{*}}dv_{\mathbf{g}_{t}})^{2} 
\]
\[
=(\int_{B(0,R\mu _{t})}x^{k}(\eta \overline{u}_{t})^{2^{*}}dv_{\mathbf{g}%
_{t}}+\int_{B(0,\delta )\backslash B(0,R\mu _{t})}x^{k}(\eta \overline{u}%
_{t})^{2^{*}}dv_{\mathbf{g}_{t}})^{2} 
\]
\[
\leq 2(\int_{B(0,R\mu _{t})}x^{k}(\eta \overline{u}_{t})^{2^{*}}dv_{\mathbf{g%
}_{t}})^{2}+2(\int_{B(0,\delta )\backslash B(0,R\mu _{t})}x^{k}(\eta 
\overline{u}_{t})^{2^{*}}dv_{\mathbf{g}_{t}})^{2} 
\]
The blow-up formula's give, for a fixed $R$ : 
\[
\frac{(\int_{B(0,R\mu _{t})}x^{k}(\eta \overline{u}_{t})^{2^{*}}dv_{\mathbf{g%
}_{t}})^{2}}{\int_{B(0,\delta )}\overline{u}_{t}^{2}dv_{\mathbf{g}_{t}}}\leq 
\frac{(\mu _{t}\int_{B(0,R)}x^{k}\widetilde{u}_{t}^{2^{*}}dv_{\widetilde{%
\mathbf{g}}_{t}})^{2}}{\mu _{t}^{2}\int_{B(0,R)}\widetilde{u}_{t}^{2}dv_{%
\widetilde{\mathbf{g}}_{t}}.}\stackunder{t\rightarrow 1}{\rightarrow }\frac{%
(\int_{B(0,R)}x^{k}\widetilde{u}^{2^{*}}dx)^{2}}{\int_{B(0,R)}\widetilde{u}%
^{2}dx}=0 
\]
because $\widetilde{u}$ is radial.

At last, using the weak estimates: $\,d_{\mathbf{g}}(x,x_{t})^{\frac{%
n-2}{2}}u_{t}(x)\leq \varepsilon $ if $d_{\mathbf{g}%
}(x,x_{t})\geq R\mu _{t}$, and using the Hölder's inequality: 
\begin{eqnarray*}
(\int_{B(0,\delta )\backslash B(0,R\mu _{t})}x^{k}(\eta \overline{u}%
_{t})^{2^{*}}dv_{\mathbf{g}_{t}})^{2} && \leq \varepsilon
_{R}^{2}(\int_{B(0,\delta )\backslash B(0,R\mu _{t})}\overline{u}_{t}^{2%
\frac{n-1}{n-2}}dv_{\mathbf{g}_{t}})^{2} \\ 
&& \leq \varepsilon _{R}^{2}(\int_{B(0,\delta )\backslash B(0,R\mu _{t})}%
\overline{u}_{t}^{2}dv_{\mathbf{g}_{t}})(\int_{B(0,\delta )\backslash
B(0,R\mu _{t})}\overline{u}_{t}^{\frac{2n}{n-2}}dv_{\mathbf{g}_{t}})
\end{eqnarray*}
therefore 
\[
\frac{(\int_{B(0,\delta )\backslash B(0,R\mu _{t})}x^{k}(\eta \overline{u}%
_{t})^{2^{*}}dv_{\mathbf{g}_{t}})^{2}}{\int_{B(0,\delta )}\overline{u}%
_{t}^{2}dv_{\mathbf{g}_{t}}}\leq \varepsilon _{R}^{2}(\int_{B(0,\delta
)\backslash B(0,R\mu _{t})}\overline{u}_{t}^{\frac{2n}{n-2}}dv_{\mathbf{g}%
_{t}})\leq c\varepsilon _{R}^{2} 
\]
where $\varepsilon _{R}\rightarrow 0$ when $R\rightarrow \infty $. Remarking that because $x_{0}$ 
is a concentration point: 
\[
z_{t}^{2}=\int_{B(0,\delta )}(\eta \overline{u}_{t})^{2^{*}}dv_{\mathbf{g}%
_{t}}\geq \int_{B(x_{0},\delta /4)}u_{t}^{2^{*}}dv_{\mathbf{g}}\geq c>0\text{
} 
\]
we have obtained
\[
\frac{\left| \varepsilon ^{k}(t)\varepsilon ^{l}(t)\right| }{%
z_{t}^{2}\int_{B(0,\delta )}\overline{u}_{t}^{2}dv_{\mathbf{g}_{t}}}%
\stackunder{t\rightarrow 1}{\rightarrow }0 
\]
We have therefore obtained once again that
\[
h(x_{0})+\varepsilon _{\delta }\leq \frac{n-2}{4(n-1)}S_{\mathbf{g}}(x_{0})-%
\frac{(n-2)(n-4)}{8(n-1)}\frac{\bigtriangleup _{\mathbf{g}}f(x_{0})}{f(x_{0})%
}+\varepsilon _{\delta }\,\,. 
\]
Letting $\delta $ tend to 0: 
\[
h(x_{0})\leq \frac{n-2}{4(n-1)}S_{\mathbf{g}}(x_{0})-\frac{(n-2)(n-4)}{8(n-1)%
}\frac{\bigtriangleup _{\mathbf{g}}f(x_{0})}{f(x_{0})} 
\]
which contradict our hypothesis: 
\[
h(x_{0})>\frac{n-2}{4(n-1)}S_{\mathbf{g}}(x_{0})-\frac{(n-2)(n-4)}{8(n-1)}%
\frac{\bigtriangleup _{\mathbf{g}}f(x_{0})}{f(x_{0})} 
\]
when $x_{0}$ is a point of maximum of $f$. This prove that $u_{t}\rightarrow u>0$, a minimizing solution for $\E$,
and therefore the weakly critical function $h$ is in fact critical.

We now prove the estimate
\[
\frac{d_{\mathbf{g}}(x_{t},x_{0})}{\mu _{t}}\leq C 
\]

Going back to the computations above, we have obtained:
\begin{eqnarray}
h(x_{0}) &\leq &\frac{n-2}{4(n-1)}S_{\mathbf{g}}(x_{0})-\frac{(n-2)(n-4)}{%
8(n-1)}\frac{\bigtriangleup _{\mathbf{g}}f(x_{0})}{f(x_{0})}+\varepsilon
_{\delta }  \\
&&+\stackunder{t\rightarrow 1}{\overline{\lim }}\frac{n-2}{n}\frac{%
D_{kl}(t,\delta )(x_{0}^{k}(t).z_{t}-\frac{\varepsilon ^{k}(t)}{z_{t}}%
)(x_{0}^{l}(t).z_{t}-\frac{\varepsilon ^{l}(t)}{z_{t}})}{\int_{B(0,\delta )}%
\overline{u}_{t}^{2}dv_{\mathbf{g}_{t}}}  \nonumber
\end{eqnarray}
where $D_{kl}(t,\delta )$ is negative definite for $t$ close to 1
and for all $\delta $ small enough, and where we remind that $%
x_{0}(t)=\exp _{x_{t}}^{-1}(x_{0})=(x_{0}^{1}(t),...,x_{0}^{n}(t))$. So, there exists a $\lambda >0$ such that 
$\forall \omega \in \Bbb{R}^{n}:$%
\[
D_{kl}(t,\delta )\omega ^{k}\omega ^{l}\leq -\lambda \sum_{k}\left| \omega
^{k}\right| ^{2} 
\]
and so
\[
D_{kl}(t,\delta )\frac{(x_{0}^{k}(t).z_{t}-\frac{\varepsilon ^{k}(t)}{z_{t}}%
)(x_{0}^{l}(t).z_{t}-\frac{\varepsilon ^{l}(t)}{z_{t}})}{(\int_{B(0,\delta )}%
\overline{u}_{t}^{2}dv_{\mathbf{g}_{t}})^{\frac{1}{2}}(\int_{B(0,\delta )}%
\overline{u}_{t}^{2}dv_{\mathbf{g}_{t}})^{\frac{1}{2}}}\leq 
\]
\[
-\lambda \sum_{k}\left| \frac{x_{0}^{k}(t).z_{t}}{(\int_{B(0,\delta )}%
\overline{u}_{t}^{2}dv_{\mathbf{g}_{t}})^{\frac{1}{2}}}-\frac{\varepsilon
^{k}(t)}{z_{t}(\int_{B(0,\delta )}\overline{u}_{t}^{2}dv_{\mathbf{g}_{t}})^{%
\frac{1}{2}}}\right| ^{2}\,\,. 
\]
Moreover, we already proved that: 
\[
\frac{\varepsilon ^{k}(t)^{2}}{z_{t}^{2}\int_{B(0,\delta )}\overline{u}%
_{t}^{2}dv_{\mathbf{g}_{t}}}\stackunder{t\rightarrow 1}{\rightarrow }0 
\]
as we also have $z_{t}=(\int_{B(0,\delta )}\overline{u}_{t}^{2^{*}}dv_{%
\mathbf{g}_{t}})^{\frac{1}{2}}$, and therefore as $x_{0}$ is a concentration point: 
\[
0<c\leq \lim \inf \,z_{t}\leq \lim \sup \,z_{t}\leq c^{\prime }<+\infty . 
\]
Therefore, necessarilly, because of (20), for all $k$, there exists a constant 
$C>0$ such that for $t\rightarrow 1$ : 
\[
\,\frac{x_{0}^{k}(t)}{(\int_{B(0,\delta )}\overline{u}_{t}^{2}dv_{\mathbf{g}%
_{t}})^{\frac{1}{2}}}\leq C 
\]
Now
\[
\int_{B(0,\delta )}\overline{u}_{t}^{2}dv_{\mathbf{g}_{t}}=\mu
_{t}^{2}\int_{B(0,\delta \mu _{t}^{-1})}\widetilde{u}_{t}^{2}dv_{\widetilde{%
\mathbf{g}}_{t}}. 
\]
But the strong estimates give that 
\[
\stackunder{t\rightarrow 1}{\overline{\lim }}\int_{B(0,\delta \mu _{t}^{-1})}%
\widetilde{u}_{t}^{2}dv_{\widetilde{\mathbf{g}}_{t}}<+\infty 
\]
therefore
\[
\int_{B(0,\delta )}\overline{u}_{t}^{2}dv_{\mathbf{g}_{t}}\sim C\mu _{t}^{2} 
\]
from where we have
\[
\forall k:\,\frac{x_{0}^{k}(t)}{\mu _{t}}\leq C^{\prime } 
\]
and so 
\[
\frac{d_{\mathbf{g}}(x_{t},x_{0})}{\mu _{t}}\leq C 
\]

If we have furthermore that at the points of maximum of $f$ : 
\[
h(P)=\frac{n-2}{4(n-1)}S_{\mathbf{g}}(P)-\frac{(n-2)(n-4)}{8(n-1)}\frac{%
\bigtriangleup _{\mathbf{g}}f(P)}{f(P)} 
\]
then we have more precisely that 
\[
\frac{d_{\mathbf{g}}(x_{t},x_{0})}{\mu _{t}}\rightarrow 0 
\]

\textit{Remark:} Note that when concentration occurs we have: 
\[
h(x_{0})\leq \frac{n-2}{4(n-1)}S_{\mathbf{g}}(x_{0})-\frac{(n-2)(n-4)}{8(n-1)%
}\frac{\bigtriangleup _{\mathbf{g}}f(x_{0})}{f(x_{0})} 
\]

\section{Critical triple 1: existence of critical functions}
The idea to prove the existence of critical functions (theorem 2), is to find, being given the manifold $(M,\g)$ 
and the function $f$, a subcritical function $h_{0}$ and a weakly critical function $h_{1}$ and then to join these 
two functions by a continuous path; theorem 1 then shows that this path must "cross" the set of critical functions.

Note first that, by the sharp Sobolev inequality (2), $B_{0}(\mathbf{g})K(n,2)^{-2}$ is a weakly critical function 
for any manifold $(M,\g)$ and any function $f$. Also, it is known that 
$$B_{0}(\mathbf{g})K(n,2)^{-2}\geq \frac{n-2}{4(n-1)} \underset{M}{Sup}\, S_{\g}$$
Therefore, for any $\alpha >0$, and for any point $P$ where $f$ is maximum on $M$, we have
$$B_{0}(\mathbf{g})K(n,2)^{-2}+ \alpha \geq \frac{n-2}{4(n-1)}S_{\mathbf{g}}(P)-\frac{(n-2)(n-4)}{8(n-1)%
}\frac{\bigtriangleup _{\mathbf{g}}f(P)}{f(P)}.$$

Now, we are going to modify the weakly critical function $B_{0}(\mathbf{g})K(n,2)^{-2}+ \alpha$ by the test 
functions presented in the introduction. They can be seen under the following form: 
 for any $x\in M$ and any $\delta >0$
small enough, there exists a sequence of functions $(\psi _{k})$ with compact support in 
$B(x,\delta )$ such that for any function $h$: 
\[
J_{h,1,\mathbf{g}}(\psi _{k})=\frac{\int_{M}\left| \nabla \psi _{k}\right|
^{2}dv_{\mathbf{g}}+\int_{M}h.\psi _{k}{{}^{2}}dv_{\mathbf{g}}}{\left(
\int_{M}\left| \psi _{k}\right| ^{\frac{2n}{n-2}}dv_{\mathbf{g}}\right) ^{%
\frac{n-2}{n}}}\stackunder{k\rightarrow \infty }{\rightarrow }\frac{1}{K(n,2)%
{{}^{2}}} 
\]
and
\[
\int_{M}\psi _{k}^{\frac{2n}{n-2}}dv_{\mathbf{g}}=1 
\]
this last condition being obtained by multiplying the functions in the introduction by suitable constants. We will 
use the functional $J$ here, as 
\[
\int_{M}f.\psi _{k}^{\frac{2n}{n-2}}dv_{\mathbf{g}}\neq 1\,\,. 
\]
Let then $\psi _{k}$ be one of these functions, where $k$ and $B(x,\delta )$
will be fixed later. We consider, for $t>0$ the sequence
\[
h_{t}=B_{0}(\mathbf{g})K(n,2)^{-2}+\alpha -t.\psi _{k}^{\frac{4}{n-2}}\,\,. 
\]
First, we seek a condition for $\triangle _{\mathbf{g}%
}+h_{t} $ to be coercive. Noting $B_{0}K^{-2}=B_{0}(\mathbf{g}%
)K(n,2)^{-2}$, and taking all integrals for the measure $dv_{\mathbf{g}}$, we have for $u\in H_{1}^{2}$%
\begin{eqnarray*}
\int_{M}(\left| \nabla u\right| _{\mathbf{g}}^{2}+h_{t}u^{2}) && 
=\int_{M}(\left| \nabla u\right| _{\mathbf{g}}^{2}+B_{0}K^{-2}.u^{2})-(t-\alpha)%
\int_{M}\psi _{k}^{\frac{4}{n-2}}.u^{2} \\ 
&& \geq K(n,2)^{-2}(\int_{M}u^{\frac{2n}{n-2}})^{\frac{n-2}{n}}-(t-\alpha)\int_{M}\psi
_{k}^{\frac{4}{n-2}}.u^{2}
\end{eqnarray*}
by Sobolev inequality. But using Hölder's inequality:
\[
\int_{M}\psi _{k}^{\frac{4}{n-2}}.u^{2}\leq (\int_{M}\psi _{k}^{\frac{2n}{n-2%
}})^{\frac{n-2}{n}}(\int_{M}u^{\frac{2n}{n-2}})^{\frac{2}{n}}=(\int_{M}u^{%
\frac{2n}{n-2}})^{\frac{2}{n}} 
\]
as $\int_{M}\psi _{k}^{\frac{2n}{n-2}}=1$. Thus, using Hölder's inequality again to get the existence of a constant
 $C>0$ such that 
\[
C\int_{M}u^{2}\leq (\int_{M}u^{\frac{2n}{n-2}})^{\frac{n-2}{n}} 
\]
we have as soon as $K(n,2)^{-2}-(t-\alpha)>0$
\begin{eqnarray*}
\int_{M}(\left| \nabla u\right| _{\mathbf{g}}^{2}+h_{t}u^{2})& & \geq
(K(n,2)^{-2}-(t-\alpha))(\int_{M}u^{\frac{2n}{n-2}})^{\frac{n-2}{n}} \\ 
& &\geq (K(n,2)^{-2}-(t-\alpha))C\int_{M}u^{2}\,\,\,.
\end{eqnarray*}
So $\triangle _{\mathbf{g}}+h_{t}$ is coercive as soon as $t-\alpha <K(n,2)^{-2}$; 
we then fix $t_{1}$ such that $\alpha <t_{1}<K(n,2)^{-2}+ \alpha$. 

We now want to fix 
$\psi _{k}$ so that $h_{t_{1}}$ is subcritical for $f$. We pick first $x$ close enough to a point $x_{0}$ 
of maximum of $f$
and $\delta $ small enough such that $f>0$ on $B(x,\delta )$, to obtain : 
\begin{eqnarray*}
J_{h_{t_{1}},f,\mathbf{g}}(\psi _{k}) && =\frac{\int_{M}\left| \nabla \psi
_{k}\right| ^{2}+\int_{M}B_{0}K^{-2}.\psi _{k}{{}^{2}}-(t_{1}-\alpha)\int_{M}\psi
_{k}^{\frac{2n}{n-2}}}{\left( \int_{M}f\left| \psi _{k}\right| ^{\frac{2n}{%
n-2}}\right) ^{\frac{n-2}{n}}} \\ 
&& \leq \frac{J_{B_{0}K^{-2},1}(\psi _{k})}{(\stackunder{B(x,\delta )}{Inf}%
f)^{\frac{n-2}{n}}}-\frac{t_{1}-\alpha}{(\stackunder{B(x,\delta )}{Sup}f)^{\frac{n-2%
}{n}}} \\ 
&& \leq \frac{J_{B_{0}K^{-2},1}(\psi _{k})}{(\stackunder{B(x,\delta )}{Inf}%
f)^{\frac{n-2}{n}}}-\frac{t_{1}-\alpha}{(\stackunder{M}{Sup}f)^{\frac{n-2}{n}}}%
\,\,\,.
\end{eqnarray*}
For any $\varepsilon >0$, by continuity of $f$, we can choose $x$
close enough to a point of maximum $x_{0}$ and $\delta $ small enough such that $B(x,\delta
)\cap \left\{ x/f(x)=Maxf\right\} =\emptyset $ and
\[
\frac{1}{(\stackunder{B(x,\delta )}{Inf}f)^{\frac{n-2}{n}}}\leq \frac{1}{(%
\stackunder{M}{Sup}f)^{\frac{n-2}{n}}}+\varepsilon 
\]
$x$ and $\delta $ being fixed, we can now choose $k$ large enough to have
\[
J_{B_{0}K^{-2},1}(\psi _{k})\leq K(n,2)^{-2}+\varepsilon \,\,. 
\]
Therfore, choosing $\varepsilon $ small enough, we see that because$\frac{t_{1}-\alpha}{(\stackunder{M}{Sup}f)^{\frac{n-2}{n}}}>0$ : 
\[
J_{h_{t_{1}},f,\mathbf{g}}(\psi _{k})<\frac{1}{K(n,2)^{-2}(\stackunder{M}{Sup%
}f)^{\frac{n-2}{n}}} 
\]
and therefore $h_{t_{1}}$ is subcritical for $f$. We now set:
\[
t_{0}=Inf\{t\leq t_{1}/\lambda _{h_{t}}<\frac{1}{K(n,2){{}^{2}}(\stackunder{M%
}{Sup}f)^{\frac{n-2}{n}}}\}\text{ }. 
\]
Then $t_{0}\geq 0$, and
\[
\lambda _{h_{t_{0}}}=\frac{1}{K(n,2){{}^{2}}(\stackunder{M}{Sup}f)^{\frac{n-2%
}{n}}}\,\,\,and\,\,\,\lambda _{h_{t}}<\frac{1}{K(n,2){{}^{2}}(\stackunder{M}{%
Sup}f)^{\frac{n-2}{n}}}\,\,\,if\,\,\,t>t_{0}. 
\]
Furthermore $\forall t$, $t_{0}\leq t\leq t_{1}$,

\[
\frac{4(n-1)}{n-2}h_{t_{0}}(P)>S_{\mathbf{g}}(P)-\frac{n-4}{2}\frac{%
\bigtriangleup _{\mathbf{g}}f(P)}{f(P)}\, for \, P\in \left\{ x/f(x)=Maxf\right\} 
\]
because $B(x,\delta )\cap \left\{ x/f(x)=Maxf\right\} =\emptyset $. At last, $%
h_{t}\stackunder{t\rightarrow t_{0}}{\rightarrow }h_{t_{0}}$ in 
$C^{0,\alpha }$, and $\triangle _{\mathbf{g}}+h_{t_{0}}$ is coercive. Therefore by theorem 1, $h_{t_{0}}$ 
is critical and $\E$ has minimizing solutions.

Now, we prove that if $\left\{ x/f(x)=Maxf\right\}$ is thin and if $\int_{M}f>0$, there exist positive critical functions. 
We start again with $h=B_{0}(\mathbf{g})K(n,2)^{-2}+\alpha$, with $\alpha>0$.
For all $P$ where $f$ is maximum on $M$ : 
\[
\frac{4(n-1)}{n-2}B_{0}(\mathbf{g})K(n,2)^{-2}+\alpha>S_{\mathbf{g}}(P)-\frac{n-4}{2%
}\frac{\bigtriangleup _{\mathbf{g}}f(P)}{f(P)} 
\]
as
\[
B_{0}(\mathbf{g})K(n,2)^{-2}\geq \frac{(n-2)}{4(n-1)}MaxS_{\mathbf{g}}\,\,. 
\]
As $f$ is not constant, there exist $\eta $ with support in $M\backslash \left\{ x/f(x)=Maxf\right\} $ 
and such that $0\leq \eta \leq 1$. Let
\[
c=\left( \int_{M}fdv_{\mathbf{g}}\right) ^{-\frac{n-2}{n}} 
\]
That's where we need $\int_{M}fdv_{\mathbf{g}}>0$. We have $\int
fc^{2^{*}}dv_{\mathbf{g}}=1$.
For $t\in \Bbb{R}^{+}$ we set
\[
h_{t}=B_{0}K^{-2}+\alpha -t\eta \,\,. 
\]
Then $h_{t}=B_{0}K^{-2}+\alpha $ on $\left\{ x/f(x)=Maxf\right\} $, and
\[
I_{h_{t}}(c)=\int_{M}(B_{0}K^{-2}+\alpha)c{{}^{2}}dv_{\mathbf{g}}-c{{}^{2}}%
t\int_{M}\eta dv_{\mathbf{g}}=\left( \int_{M}fdv_{\mathbf{g}}\right) ^{-%
\frac{2}{2^{*}}}((B_{0}K^{-2}+\alpha){Vol}_{\mathbf{g}}(M)-t\int_{M}\eta dv_{\mathbf{g%
}}) 
\]
So, if $t$ is large enough, 
\[
I_{h_{t}}(c)<\frac{1}{K(n,2){{}^{2}}(\stackunder{M}{Sup}f)^{\frac{n-2}{n}}}%
\,\,. 
\]
We also want $h_{t}$ to be positive on $M$.
By the definition of $h_{t}$ and because $\stackunder{M}{Sup}\,\eta =1$, it is the case if 
\begin{equation}
t<B_{0}(\mathbf{g})K(n,2)^{-2}+\alpha.  
\end{equation}
But we also want that
\[
I_{h_{t}}(c)<\frac{1}{K(n,2){{}^{2}}(\stackunder{M}{Sup}f)^{\frac{n-2}{n}}} 
\]
which requires
\begin{equation}
t>\frac{1}{\int_{M}\eta dv_{\mathbf{g}}}\Bigl( (B_{0}K^{-2}+\alpha){Vol}_{\mathbf{g}}(M)-%
\frac{\left( \int_{M}fdv_{\mathbf{g}}\right) ^{\frac{n-2}{n}}}{K(n,2){{}^{2}}%
(\stackunder{M}{Sup}f)^{\frac{n-2}{n}}}\Bigr)\,\,. 
\end{equation}
We can find such a $t$ if : 
\[
\frac{1}{\int_{M}\eta dv_{\mathbf{g}}}\Bigl( (B_{0}K^{-2}+\alpha){Vol}_{\mathbf{g}}(M)-%
\frac{\left( \int_{M}fdv_{\mathbf{g}}\right) ^{\frac{n-2}{n}}}{K(n,2){{}^{2}}%
(\stackunder{M}{Sup}f)^{\frac{n-2}{n}}}\Bigr)<B_{0}K^{-2}+\alpha
\]
which can be writen 
\begin{equation}
\int_{M}\eta dv_{\mathbf{g}}>{Vol}_{\mathbf{g}}(M)-\frac{K^{-2}\left( \int_{M}fdv_{%
\mathbf{g}}\right) ^{\frac{n-2}{n}}}{(B_{0}K^{-2}+\alpha)(\stackunder{M}{Sup}%
f)^{\frac{n-2}{n}}}\,\,.  
\end{equation}
Remember that we want $\eta $ to have support in $M\backslash \left\{
x/f(x)=Maxf\right\} $ with $0\leq \eta \leq 1$. But we made the hypothesis that $\left\{ x/f(x)=Maxf\right\} $, 
the set of maximum points of $f$, is a thin set. We can therefore find such a function $\eta $ with
$\int_{M}\eta dv_{\mathbf{g}}$ as close as we want to ${Vol}_{\mathbf{g}}(M)$.
As
\[
\frac{\left( \int_{M}fdv_{\mathbf{g}}\right) ^{\frac{n-2}{n}}}{B_{0}(\mathbf{%
g})(\stackunder{M}{Sup}f)^{\frac{n-2}{n}}}>0 
\]
we can find $\eta $ satisfying (23) and a real $t$, denoted $t_{1}$, satisfying (21) and (22).

On the set $\left\{ x/f(x)=Maxf\right\} $, $h_{t}=B_{0}K^{-2}+\alpha$, so $\forall P\in
\left\{ x/f(x)=Maxf\right\} $:

\[
\frac{4(n-1)}{n-2}h_{t}(P)>S_{\mathbf{g}}(P)-\frac{n-4}{2}\frac{%
\bigtriangleup _{\mathbf{g}}f(P)}{f(P)}. 
\]
We then set: 
\[
t_{0}=Inf\{t\leq t_{1}\,/\,\lambda _{h_{t}}<\frac{1}{K(n,2){{}^{2}}(%
\stackunder{M}{Sup}f)^{\frac{n-2}{n}}}\}\text{ }. 
\]
Necessarilly, $t_{0}<t_{1}$. We remind that (see section 1):
\[
\lambda _{h,f,\mathbf{g}}=\lambda _{h}=\stackunder{w\in \mathcal{H}_{f}}{%
\inf }I_{h}(w) 
\]
Therefore
\[
\lambda _{h_{t_{0}}}=\frac{1}{K(n,2){{}^{2}}(\stackunder{M}{Sup}f)^{\frac{n-2%
}{n}}}\,\,\,and\,\,\,\lambda _{h_{t}}<\frac{1}{K(n,2){{}^{2}}(\stackunder{M}{%
Sup}f)^{\frac{n-2}{n}}}\,\,\,if\,\,\,t>t_{0}. 
\]
Furthermore $\forall t$, $t_{0}\leq t\leq t_{1}$, $h_{t}>0$ on $M$ and

\[
\frac{4(n-1)}{n-2}h_{t_{0}}(P)>S_{\mathbf{g}}(P)-\frac{n-4}{2}\frac{%
\bigtriangleup _{\mathbf{g}}f(P)}{f(P)}\,\,\,for\,\,\,P\in \left\{
x/f(x)=Maxf\right\} . 
\]
At last $h_{t}\stackunder{t\rightarrow t_{0}}{\rightarrow }h_{t_{0}}$ in $C^{0}$, and as 
$h_{t_{0}}>0$, $\triangle _{\mathbf{g}}+h_{t_{0}}$ is coercive. 
Therefore by theorem 1, $h_{t_{0}}$ is critical and $\E$ has minimizing solutions.

\textit{Remark:} The precceding proofs also show, by replacing $B_{0}K^{-2}$ by $h$, that if 
$(h,f,\mathbf{g})$ is weakly critical, and if
\[
\frac{4(n-1)}{n-2}h(P)>S_{\mathbf{g}}(P)-\frac{n-4}{2}\frac{\bigtriangleup _{%
\mathbf{g}}f(P)}{f(P)}\text{ \thinspace for \thinspace }P\in \left\{
x/f(x)=Maxf\right\} 
\]
then there exists $h^{\prime }$ $\leqslant h$ such that$ (h^{\prime },f,g)$ is critical.

If we only have
\[
\frac{4(n-1)}{n-2}h(P)\geq S_{\mathbf{g}}(P)-\frac{n-4}{2}\frac{%
\bigtriangleup _{\mathbf{g}}f(P)}{f(P)}\text{ \thinspace for \thinspace }%
P\in \left\{ x/f(x)=Maxf\right\} 
\]
then for any $\varepsilon >0$ there exists $h^{\prime }$ $\leq h+\varepsilon $ such that $(h^{\prime },f,g)$ 
is critical.

Weaker hypothesis are sufficient to prove the existence of positive critical functions: for example, it suffices that 
the \textit{boundary} of the set $Max f$ is a set of null measure; see \cite{C3} for full details.

\section{Critical triple 2}
We want to prove here theorem 3.
This theorem lies on the transformation formula for a critical function in a conformal change of metric (seen at 
the end of the introduction): 
\begin{center}
($h^{\prime },f$\textit{,}$\mathbf{g}^{\prime }=u^{\frac{4}{n-2}}\mathbf{g)}$
\textit{ is critical if and only if (}$h=h^{\prime }u^{\frac{4}{n-2}}-%
\frac{\triangle _{\mathbf{g}}u}{u}$\textit{,}$f$\textit{,}$\mathbf{g)}$ 
\textit{is critical.}
\end{center}
We set, for $u\in C_{+}^{\infty }(M)=\{u\in C^{\infty }(M)\,/\,u>0\}$ : 
\[
F_{h^{\prime }}(u)=h^{\prime }u^{\frac{4}{n-2}}-\frac{\bigtriangleup _{%
\mathbf{g}}u}{u} 
\]
Then:
\begin{center}
$(h^{\prime },f,\mathbf{g}^{\prime })$\textit{\ is critical if and only if }
$(F_{h^{\prime }}(u),f,\mathbf{g})$\textit{\ is critical.}
\end{center}
To prove the theorem, we therefore have to prove the existence of a function $h$ such that:

1/: $\triangle _{\mathbf{g}}u+h.u=h^{\prime }u^{\frac{n+2}{n-2}}$ has a solution $u>0$, and

2/: $(h,f,\mathbf{g})$ is critical.

Indeed, in this case $h=F_{h^{\prime }}(u)$ and $h^{\prime }$ is critical for 
$f$ and $\mathbf{g}^{\prime }=u^{\frac{4}{n-2}}\mathbf{g}$.

E. Humbert et M. Vaugon proved this theorem in the case $f=cste$
and for a manifold not conformaly diffeomorphic to the sphere \cite{Hu-V}. Their method lies on the fact that for such a manifold ,
after a first conformal change of metric, $B_{0}(%
\mathbf{g})K(n,2)^{-2}$ is a critical function, (we will denote these two constants $K$ et $B_{0}$). In fact, 
a careful study of their proof shows that what is needed is in fact that $B_{0}K^{-2}$ is positive. But we proved 
in the previous section the existence of positive critical functions under a geometric hypothesis concerning $f$. 
Remark that our proof will work on the sphere, but only for a non-constant function $f$.

The principle of the proof of E. Humbert and M. Vaugon is the following. We know that there exists a sequence 
($h_{t}$) of sub-critical functions for $f$ and $\mathbf{g}$ such that $h_{t}%
\stackrel{C^{2}}{\rightarrow }h$ where $(h,f,\mathbf{g})$ is critical and such that for any point $P$ where $f$ is maximum on $M$
\[
\frac{4(n-1)}{n-2}h(P)>S_{\mathbf{g}}(P)-\frac{n-4}{2}\frac{\bigtriangleup _{%
\mathbf{g}}f(P)}{f(P)}. 
\]
For a sequence $q_{t}\rightarrow 2^{*},\,q_{t}<2^{*}$ we build a sequence
$u_{t}>0$ of solutions of 
\[
\triangle _{\mathbf{g}}u+h.u=h^{\prime }u^{q_{t}-1}\,\,\,with\,\,\,\int
h^{\prime }u_{t}^{q_{t}}dv_{\mathbf{g}}\leq C\text{ independant of }t 
\]
such that $u_{t}\stackrel{H_{1}^{2}}{\rightharpoondown }u\geqslant 0$. Here again, if $u>0$, 
then $u$ is solution (up to a multiplicative constant) of $\triangle _{\mathbf{g}}u+h.u=h^{\prime }u^{\frac{n+2}{n-2}}$
and we are done.

Now, if $u=0$, one shows that the $u_{t}$ concentrate and that using this phenomenom, one can find a $t_{0}$ 
close to 1 ( if e.g. $t\rightarrow 1$) and a real $s$ large, such that $F_{h^{\prime
}}(u_{t_{0}})$ is sub-critical and $F_{h^{\prime }}(u_{t_{0}}^{s})$
is weaklly critical, with furthermore
\[
\frac{4(n-1)}{n-2}F_{h^{\prime }}(u_{t_{0}}^{s})(P)>S_{\mathbf{g}}(P)-\frac{n-4}{2}\frac{%
\bigtriangleup _{\mathbf{g}}f(P)}{f(P)} 
\]
at any point $P$ where $f$ is maximum. Then, considering the path 
$t\rightarrow F_{h^{\prime }}(u_{t_{0}}^{ts})$ and using theorem 1, we get the existence of a critical 
function on this path. It is to obtain the conditions on $F_{h^{\prime }}(u_{t_{0}}^{s})$ at the maximum points 
of $f$ that we need the existence of positive critical functions.

We will now give the scheme of the proof, refering for complete details to the article of E. Humbert and M. Vaugon or to our PHD 
thesis available online, 
and we will only indicate the modifications due to our function $f$ and the necessity of positive critical functions.

First, we said that we will need positive critical functions. Their existence was proved under the hypothesis that 
$Max f$ is thin and that $\int_{M}fdv_{\mathbf{g}}>0$. But $\stackunder{M}{Sup}f>0$,
so, after making if necessary a first conformal change of metric, we can suppose that 
$\int_{M}fdv_{\mathbf{g}}>0$, and we supposed in the hypothesis of theorem 3 that $Max f$ is thin, and therefore 
we can suppose that we have positive critical function for $f$ and $\mathbf{g}$.

Then, we fix some (more) notations:
\begin{center}
\[
J_{h,h^{\prime },\mathbf{g,}q}(w)=\frac{\int_{M}\left| \nabla w\right|
^{2}dv_{\mathbf{g}}+\int_{M}h.w{{}^{2}}dv_{\mathbf{g}}}{\left(
\int_{M}h^{\prime }\left| w\right| ^{q}dv_{\mathbf{g}}\right) ^{\frac{2}{q}}}
\]
\end{center}

\[
\stackunder{w\in \mathcal{H}_{h^{\prime },q}^{+}}{\inf }J_{h,h^{\prime },%
\mathbf{g,}q}(w):=\lambda _{h,h^{\prime },\mathbf{g,}q} 
\]
where
\[
\mathcal{H}_{h^{\prime },q}^{+}=\{w\in
H_{1}^{2}(M)\,/\,\,\,w>0\,\,and\,\,\int_{M}h^{\prime }.w^{q}dv_{\mathbf{g}%
}>0\}. 
\]
and
\[
\Omega _{h,h^{\prime },\mathbf{g,}q}=\{u\in \mathcal{H}_{h^{\prime
},q}^{+}/\,J_{h,h^{\prime },\mathbf{g,}q}(u)=\lambda _{h,h^{\prime },\mathbf{%
g,}q}\,\,and\,\,\int_{M}h^{\prime }.w^{q}dv_{\mathbf{g}}=(\lambda
_{h,h^{\prime },\mathbf{g,}q}\,\,)^{\frac{q}{q-2}}\}\,\,. 
\]

Let ($h_{t}$) be a sequence of sub-critical functions for $f$ and $%
\mathbf{g}$ such that $h_{t}\stackrel{C^{2}}{\rightarrow }h$ where $(h,f,%
\mathbf{g})$ is critical, with $\triangle _{\mathbf{g}}+h_{t}$ coercive. We know that we can find such a 
sequence with $h_{t}>0$ et $h>0$, and also
\[
\frac{4(n-1)}{n-2}h_{t}(P)>S_{\mathbf{g}}(P)-\frac{n-4}{2}\frac{\bigtriangleup _{\mathbf{g}%
}f(P)}{f(P)} 
\]
for all $P\in Max f$. But here, we can say more, and that is where the existence of positive critical functions is 
crucial. Indeed, for any constant $c>0$, if $\mathbf{g}^{\prime }=c%
\mathbf{g}$, then $S_{\mathbf{g}^{\prime }}=c^{-1}S_{\mathbf{g}}$ and $%
\triangle _{\mathbf{g}^{\prime }}=c^{-1}\triangle _{\mathbf{g}}$ and by the transformation formula 
for critical functions:

\begin{center}
\textit{\ }$h$\textit{\ is (sub-, weakly) critical for }$f$\textit{\
and }$\mathbf{g}$\textit{\ if and onlu if }$c^{-1}h$\textit{\ is (sub-,
weakly) critical for }$f$\textit{\ and }$\mathbf{g}^{\prime }$.
\end{center}

Therefore, up to multiplying $\mathbf{g}$ by a constant, we can, for any constant $C>0$, suppose : 
\begin{eqnarray*}
h_{t} &>&C\,\,\,\,on\, M \\
\frac{4(n-1)}{n-2}h_{t}(P)-S_{\mathbf{g}}(P)+\frac{n-4}{2}\frac{\bigtriangleup _{\mathbf{g}%
}f(P)}{f(P)} &>&C\,\,\,\,\forall P \in Max f
\end{eqnarray*}
and $(h,f,\mathbf{g)}$ has minimizing solutions.

We can now follow the method exposed above; we only give the scheme of the proof.

\textit{First step:} Thanks to the compacity of the inclusion $H_{1}^{2}\subset L^{q}$, it is known 
that  $\forall q<2^{*}$ and $\forall u\in \Omega _{h,h^{\prime },\mathbf{g,}q}$, $u$ is solution of 
$\triangle _{\mathbf{g}}u+h.u=h^{\prime }u^{q-1}$. Using this fact, in the first step, one proves the following:

There exist sequences $(q_{i}),(t_{i}),$ such that
\begin{eqnarray*}
\,2 &<&q_{i}<2^{*} \\
\,q_{i} &\rightarrow &2^{*} \\
\,t_{i} &\rightarrow &1 \\
\,h_{t_{i}} &\rightarrow &h\,
\end{eqnarray*}
and a sequence ($v_{i})\in \Omega _{h_{t_{i}},h^{\prime },\mathbf{g,}q_{i}}$
such that ($F_{h^{\prime }}(v_{i}),f,\mathbf{g)\,}$ is sub-critical.

We note 
\[
J_{i}=J_{h_{t_{i}},h^{\prime },\mathbf{g,}q_{i}} 
\]
$\,$and
\[
\lambda _{i}=\lambda _{h_{t_{i}},h^{\prime },\mathbf{g,}q_{i}}\,\,\,. 
\]
Then
\[
J_{i}(v_{i})=\lambda _{i}\,\,\,and\,\,\int h^{\prime }v_{i}^{q_{i}}dv_{%
\mathbf{g}}=\lambda _{i}^{\frac{q_{i}}{q_{i}-2}} 
\]
and $v_{i}$ is a positive solution of 
\[
\triangle _{\mathbf{g}}v_{i}+h_{t_{i}}.v_{i}=h^{\prime }v_{i}^{q_{i}-1}. 
\]
The sequence ($v_{i})$ is bounded in $H_{1}^{2}$ and thus there exists $v\in
H_{1}^{2}$ such that
\[
v_{i}\stackrel{H_{1}^{2}}{\rightharpoondown }v,\,\,v_{i}\stackrel{L^{2}}{%
\rightarrow }v\,\,\,et\,\,\,\,v_{i}\stackrel{L^{2^{*}-2}}{\rightarrow }v. 
\]
Once again, we have two possibilities: $v\equiv 0$ or $v>0$.

\textit{Second step:}

If $v>0$, as we said above, the proof is over: up to a subsequence, $v_{i}%
\stackrel{C^{2}}{\rightarrow }v$ and so on one hand $F_{h^{\prime
}}(v_{i})\rightarrow F_{h^{\prime }}(v)$, and on the other hand
\[
F_{h^{\prime }}(v_{i})=h_{t_{i}}+h^{\prime }(v_{i}^{\frac{4}{n-2}%
}-v_{i}^{q_{i}-2})\rightarrow h\, ,
\]
that is, $F_{h^{\prime }}(v)=h$ which is critical for $f$ and $\mathbf{g}$ with minimizing solutions. 
Thus $h^{\prime }$ is critical for $f$ and $\mathbf{g}^{\prime }=v^{\frac{4}{n-2}}\mathbf{g}$, with 
minimizing solutions.

\textit{The rest of the proof is therefore concerned with the case }$v\equiv 0$\textit{.}

\textit{Third step:} One proves that there is a concentration phenomenom:

a/: One first shows that: 
\[
0<c\leqslant \overline{\lim }\,\lambda _{i}\leqslant
K^{-2}(Sup_{M}\,h^{\prime })^{-\frac{n-2}{n}}\,\,. 
\]

b/: Second, one shows that: 
\[
0<\lambda ^{\frac{n}{2}}(Sup_{M}\,h^{\prime })^{-1}\leqslant \overline{\lim }%
\,\int_{M}v_{i}^{q_{i}}dv_{\mathbf{g}}\leqslant K^{2^{*}}\lambda ^{\frac{%
n2^{*}}{4}}\leqslant K^{-n}(Sup_{M}\,h^{\prime })^{-\frac{n}{2}} 
\]
where $\lambda >0$ is such that, after extraction, $\lambda_{i} \rightarrow \lambda$.

c/: We say that $x\in M$ is a concentration point if
\[
\forall r>0:\,\overline{\lim }\,\int_{B(x,r)}v_{i}^{q_{i}}dv_{\mathbf{g}}>0. 
\]

Using a/ and b/, and method analogous to section 4.2, one gets the following:

First, as $M$ is compact, there exists at least one concentration point $x\in M$.

Then, using the iteration process, one shows that 
\[
\overline{\lim }\,\int_{B(x,r)}v_{i}^{q_{i}}dv_{\mathbf{g}}\geqslant
K^{-n}(Sup_{M}\,h^{\prime })^{-\frac{n}{2}}. 
\]

d/: Therefore using the method of section 4.2, we get:

1/: $\overline{\lim }\,\int_{B(x,r)}v_{i}^{q_{i}}dv_{\mathbf{g}%
}=K^{-n}(Sup_{M}\,h^{\prime })^{-\frac{n}{2}}$ , $\forall r>0$

2/: $x$ is the only concentration point, denoted $x_{0}$

3/: $\lambda =K^{-2}(Sup_{M}\,h^{\prime })^{-\frac{n-2}{n}}$

4/: $x_{0}$ is a point of maximum of $h^{\prime }$

5/: $v_{i}\rightarrow 0$ in $C_{loc}^{2}(M-\{x_{0}\})$
\\

\textit{Fourth step:}

We know now that the sequence $(v_{i})$ concentrates in $x_{0}$ and that for any 
$i$ $F_{h^{\prime }}^{\prime }(v_{i})$ is sub-critical for $f$ and $\textbf{g}$. We would like to find a 
$v_{i_{0}}$ , a function $v>0$ and a continuous path from
$v_{i_{0}}$ to $v$ such that $F_{h^{\prime }}(v)$ is weakly critical for $f$ and $\textbf{g}$ and such that
\[
\frac{4(n-1)}{n-2}F_{h^{\prime }}(v)(P)>S_{\mathbf{g}}(P)-\frac{n-4}{2}\frac{\bigtriangleup _{%
\mathbf{g}}f(P)}{f(P)} 
\]
for all $P \in Max f$ . Then, the theorem 1 will tell us that on the path $u_{t}$ from $v_{i_{0}}$ to $v$ 
there exists a $u_{t}$ such that $F_{h^{\prime }}(u_{t})$ is critical for $f $ and $\mathbf{g}$ .

That is where we are going to use the existence of positive critical functions.

Let $s\geqslant 1$ and let $v$ be a positive function. Then 
\[
\bigtriangleup _{\mathbf{g}}(v^{s})=sv^{s-1}\bigtriangleup _{\mathbf{g}%
}v-s(s-1)v^{s-2}\left| \nabla v\right| _{\mathbf{g}}^{2}\,\,. 
\]
Thus
\[
F_{h^{\prime }}(v_{i}^{s})=h^{\prime }v_{i}^{s\frac{4}{n-2}%
}+sh_{t_{i}}-sh^{\prime }v_{i}^{q_{i}-2}+s(s-1)\frac{\left| \nabla v\right|
_{\mathbf{g}}^{2}}{v_{i}^{2}} 
\]
and therefore
\[
F_{h^{\prime }}(v_{i}^{s})\geqslant sh_{t_{i}}+h^{\prime }(v_{i}^{s\frac{4}{%
n-2}}-sv_{i}^{q_{i}-2})\,\,. 
\]
Now:

On $\{x\in M\,/\,h^{\prime }(x)\leqslant 0\}$ :

$v_{i}\rightarrow 0$ uniformly because $x_{0} \in Max \, h^{\prime }$ and $h^{\prime }(x_{0})>0$ as we 
have supposed that $\bigtriangleup _{\mathbf{g}}+h^{\prime }$ is coercive. Furthermore if 
$s\geqslant 1$ then $s\frac{4}{n-2}\geqslant q_{i}-2$. Thus, for $i$ large enough 
\[
F_{h^{\prime }}(v_{i}^{s})\geqslant sh_{t_{i}}\text{ on }\{x\in
M\,/\,h^{\prime }(x)\leqslant 0\} 
\]

On $\{x\in M\,/\,h^{\prime }(x)>0\}$ :

We consider the function of a real variable defined for $x\geqslant 0$ by 
\[
\beta _{i,s}(x)=x^{s\frac{4}{n-2}}-sx^{q_{i}-2}=x^{q_{i}-2}(x^{s\frac{4}{n-2}%
-q_{i}+2}-s). 
\]
An easy study of this function shows that
\[
for\,\,\,x\geqslant 0:\beta _{i,s}(x)\geqslant -s. 
\]
But 
\[
F_{h^{\prime }}(v_{i}^{s})\geqslant sh_{t_{i}}+h^{\prime }\beta
_{i,s}(v_{i}) 
\]
therefore 
\[
F_{h^{\prime }}(v_{i}^{s})\geqslant sh_{t_{i}}-sh^{\prime }\,\,\,\text{on}%
\,\,\,\{x\in M\,/\,h^{\prime }(x)>0\}. 
\]
We can therefore write : 
\[
F_{h^{\prime }}(v_{i}^{s})\geqslant s(h_{t_{i}}-\stackunder{M}{Sup}\,h^{\prime
})\text{ on }\{x\in M\,/\,h^{\prime }(x)>0\}\text{.} 
\]
We now use our work from the beginning of the proof, that is that, for any $C>0$, we can suppose that: 
\[
h_{t}>C\,\,on\, M 
\]
and 
\[
\frac{4(n-1)}{n-2}h_{t}(P)-S_{\mathbf{g}}(P)+\frac{n-4}{2}\frac{\bigtriangleup _{\mathbf{g}%
}f(P)}{f(P)}>C, \,\,\,\, \forall P\in Max f\,. 
\]
Then, first, if we suppose that $h>\stackunder{M}{Sup}\,h^{\prime }$ on $M$, we see that for $i$ and $s$ 
large enough : 
\begin{equation}
F_{h^{\prime }}(v_{i}^{s})\geqslant B_{0}(\mathbf{g})K(n,2)^{-2}  
\end{equation}
and therefore $F_{h^{\prime }}(v_{i}^{s})$ is weakly critical for $f$ and $\mathbf{g}$ . 
Beside, for all $t\in [1,s]$ we also have
\[
F_{h^{\prime }}(v_{i}^{t})\geqslant t(h_{t_{i}}-\stackunder{M}{Sup}\,h^{\prime
})\text{ }\geqslant h_{t_{i}}-\stackunder{M}{Sup}\,h^{\prime }>0 
\]
so $\bigtriangleup _{\mathbf{g}}+F_{h^{\prime }}(v_{i}^{t})$ is coercive.

Secondly, if we also suppose that
\[
\frac{4(n-1)}{n-2}h_{t}(P)-S_{\mathbf{g}}(P)+\frac{n-4}{2}\frac{\bigtriangleup _{\mathbf{g}%
}f(P)}{f(P)}>\frac{4(n-1)}{n-2}\stackunder{M}{Sup}\ h^{\prime }\,\,\,\,\, \forall P \in Max f
\]
we have for all $t\in [1,s]$: 
\begin{equation}
\frac{4(n-1)}{n-2}F_{h^{\prime }}(v_{i}^{t})(P)>S_{\mathbf{g}}(P)-\frac{n-4}{2}\frac{%
\bigtriangleup _{\mathbf{g}}f(P)}{f(P)}\,\,\,\,\,\,\, \forall P \in Max f
\end{equation}
as soon as $i$ is large enough.

We therefore fix $i$ and $s$ large enough to have (24) et (25) and we consider
\[
s_{0}=\inf \{t>1\,/\,F_{h^{\prime }}(v_{i}^{t})\text{ is weakly critical\}} 
\]

We then apply theorem 1 to the path $t\in [1,s_{0}]\mapsto
F_{h^{\prime }}(v_{i}^{t})$ to obtain that $F_{h^{\prime
}}(v_{i}^{s_{0}})$ is critical for $f$ and $\mathbf{g}$, with minimizing solutions. 
Therefore $h^{\prime }$ is critical for $f$ and $\mathbf{g}^{\prime }=(v_{i}^{s_{0}})^{\frac{4}{n-2}}\mathbf{g}$ 
with minimizing solutions.

This ends the proof.

\section{Critical triple 3}
Let $(M,\mathbf{g})$ be a compact riemannian manifold of dimension 
$n\geqslant 3$. Let $h$ be a fixed $C^{\infty }$ function such that $\triangle _{\mathbf{g}}+h$ 
is coercive. The problem we want to study is the following: \textit{can we find a function $f$ such that 
$(h,f,\g)$ is a critical triple ?}

We first make a remark. If $h\geqslant B_{0}(\mathbf{g})K(n,2){%
{}^{-2}}$, then $h$ is weakly critical for any function $f$, and there cannot exist a function $f$ such that 
($h,f,\mathbf{g})$ is subcritical. But more important is the next observation:

\textit{If there exist a non constant function $f$ such that $(h,f,\g)$ is critical with a minimizing solution $u$, 
then $(h,1,\g)$ is sub-critical. }

Indeed, as we saw in section 1, we can suppose that $Sup\,f=1$. Then, as $u>0$%
\[
J_{h,1}(u)<J_{h,f}(u)=\frac{1}{K(n,2)^{2}(Sup\,f)^{\frac{2}{2^{*}}}}=\frac{1%
}{K(n,2)^{2}} 
\]
and therefore $h$ is subcritical for $1$.

We want to prove that, at least if $dim M\geq5$, this necessary condition is sufficient, i.e we want to prove theorem 4. 
We thus suppose now that $(h,1,\mathbf{g})$\textit{\ is sub-critical.}

The proof will proceed in two steps:

First step: we prove that there exist a function $f\in C^{\infty }(M)$
such that $\stackunder{M}{Sup}f=1$, with
$\triangle _{\mathbf{g}}f$ being as large as we want in its maximum points, 
and such that $(h,f,\mathbf{g})$ is weakly critical.

Second step: being given this function $f$, we prove that there exists on the path 
\[
t\rightarrow f_{t}=t.1+(1-t)f 
\]
a function for which $h$ is critical.

\textbf{First step:}

We proceed by contradiction. We suppose that for any $f\in
C^{\infty }(M)$ such that $\stackunder{M}{Sup}f>0$, $(h,f,\mathbf{g})$
is sub-critical. Then, for all such function, there exit positive solution $u$ to the equation 
\[
\triangle _{\mathbf{g}}u+h.u=\lambda .f.u^{\frac{n+2}{n-2}} 
\]
where 
\[
\lambda =\stackunder{w\in \mathcal{H}_{f}}{\inf }I_{h,\mathbf{g}%
}(w)\,\,\,and\,\,\,\int_{M}f.u^{\frac{2n}{n-2}}dv_{\mathbf{g}}=1\,. 
\]

The metric $\mathbf{g}$ being fixed, we will not write $dv_{%
\mathbf{g}}$ in the integrals.

The idea is to build a familly of functions $f_{t}$ whose laplacians tend to infinity at the maximum points. 
One of these function will then give a weakly critical triple $(h,f_{t},\mathbf{g})$. 
Furthermore, our proof holding for any subsequence of this familly. this function will have a laplacian as large 
as we want in its point of maximum.

In $\Bbb{R}^{n}$, we build for $t\rightarrow 0$ a familly ($P_{t})$
of $C^{\infty }$ functions, similar to a regularizing sequence, such that
\begin{eqnarray*}
0 &\leqslant &P_{t}\leqslant 1 \\
\,\,P_{t}(x) &=&P_{t}(\left| x\right| ) \\
\,\,P_{t}(0) &=&1 \\
\left\| \nabla P_{t}\right\| &\sim &\frac{c_{1}}{t}\,\,\,on\,\,\,B(0,t) \\
\left| \bigtriangleup P_{t}(0)\right| &\sim &\frac{c_{2}}{t^{2}} \\
SuppP_{t} &=&B(0,t).
\end{eqnarray*}
Let now $x_{0}$ be a point of $M$ such that $h(x_{0})>0$; this point exists because 
$\triangle _{\mathbf{g}}+h$ is coercive. We define 
\[
f_{t}=P_{t}\circ \exp _{x_{0}}^{-1} 
\]
We are therefore supposing that, for all $t$, $(h,f_{t},\mathbf{g})$ is sub-critical and we are looking for a 
contradiction. For all $t$ we have a solution $u_{t}>0$ of 
\[
(E_{t}):\,\,\bigtriangleup _{\mathbf{g}}u_{t}+h.u_{t}=\lambda
_{t}.f_{t}.u_{t}^{\frac{n+2}{n-2}} 
\]
with $\int f_{t}u_{t}^{2^{*}}dv_{\mathbf{g}}=1$ and
\[
\lambda _{t}<K^{-2}(\stackunder{M}{Sup}f_{t})^{-\frac{n-2}{n}}=K^{-2}. 
\]
Then, ($u_{t})$ is bounded in $H_{1}^{2}(M)$ when $t\rightarrow 0$.
So ($u_{t})$ is bounded in $L^{2^{*}}$ and ($u_{t}^{2^{*}-1})$ is bounded in 
$L^{\frac{2^{*}}{2^{*}-1}}$. After extraction of a subsequence, if 
$f_{t}\stackrel{L^{2}}{\rightarrow }f$ and $u_{t}\stackrel{L^{2}}{\rightarrow }u$,
then
\[
f_{t}u_{t}^{2^{*}-1}\rightharpoondown fu^{2^{*}-1}. 
\]
But here, $f_{t}\stackrel{L^{p}}{\rightarrow }0$, therefore the equation ($E_{t})$
``converge'' to 
\[
\bigtriangleup _{\mathbf{g}}u+h.u=0 
\]
in the sense that $u$ is solution of this equation. But $\bigtriangleup
_{\mathbf{g}}+h$ is coercive, therefore $u=0$, i.e. $u_{t}\rightarrow 0$ in $L^{p}$ for $p<2^{*}$.

The sequence ($u_{t})$ therefore concentrates in the sense we saw in subsection 4.2. 
But in subsection 4.2, the function $f$ on the right handside of the equation was constant and it was on the 
left handside that we had a sequence $(h_{t})$. However the results we saw there remain true, only the blow-up 
necessary for the weak estimates requires a new treatment. We will go over these results, only detailing the new 
difficulties.

\textit{a/: There exists, up to a subsequence of $(u_{t})$, exactly one concentration point 
and it is the point $x_{0}$ where the $f_{t}$ are maximum on $M$. Moreover}

\[
\forall \delta >0\mathit{,\ }\overline{\stackunder{t\rightarrow 1}{\lim }}%
\int_{B(x_{0},\delta )}f_{t}u_{t}^{2^{*}}=1\,. 
\]

The method of subsection 4.2 works here. More precisely, as $Supp\,f_{t}=B(x_{0},t)$, we have for all 
$\delta >0$ and as soon as $t<\delta $:
\[
\int_{B(x_{0},\delta )}f_{t}u_{t}^{2^{*}}=1. 
\]
We can also suppose that
\[
\lambda _{t}\rightarrow \lambda =K^{-2}(\stackunder{M}{Sup}f_{t})^{-\frac{n-2%
}{n}}=K^{-2}. 
\]

\textit{b/: }$u_{t}\rightarrow 0$\textit{\ in }$C_{loc}^{0}(M-\{x_{0}\})$

Same proof as in subsection 4.2.

\textit{c/: weak estimates}

We consider a sequence of points $(x_{t})$ such that $$m_{t}=\stackunder{M}{Max}\,u_{t}=u_{t}(x_{t}):=\mu _{t}^{-\frac{n-2}{2}}.$$

From the previous point, $x_{t}\rightarrow x_{0}$ and $\mu
_{t}\rightarrow 0$. Remember that $\overline{u}_{t},\overline{f}_{t},%
\overline{h}_{t},\mathbf{g}_{t}$ are the functions and the metric seen in the chart $\exp _{x_{t}}^{-1}$, 
and $\,\,\widetilde{u}_{t}\,,\,\widetilde{h}_{t}\,,\,\widetilde{f}_{t},\widetilde{\mathbf{g}}%
_{t} $ are the functions after
blow-up of center $x_{t}$ and coefficient $k_{t}=\mu _{t}^{-1}$.

Reviewing the proof of the weak estimates in section 4.2, we see that it will work here if we obtain :

\[
\forall R>0\,\,\,:\stackunder{t\rightarrow 0}{\lim }\int_{B(x_{t},R\mu
_{t})}f_{t}u_{t}^{2^{*}}dv_{\mathbf{g}}=1-\varepsilon _{R}\,\,\,where%
\,\,\,\varepsilon _{R}\stackunder{R\rightarrow +\infty }{\rightarrow }0\,. 
\]
This relation is itself proved using blow-up theory once it is proved that $\,\widetilde{u}%
_{t}\stackrel{C_{loc}^{2}(\Bbb{R}^{n})}{\rightarrow }\,\widetilde{u}$ where
$\widetilde{u}$ is solution of :
\[
\bigtriangleup _{e}\widetilde{u}=K^{-2}\widetilde{u}^{\frac{n+2}{n-2}}. 
\]
This is where we have the main difficulty due to the presence of a familly ($f_{t}$). 
Indeed, after blow-up, the equation
\[
(E_{t}):\,\,\bigtriangleup _{\mathbf{g}}u_{t}+h.u_{t}=\lambda
_{t}.f_{t}.u_{t}^{\frac{n+2}{n-2}} 
\]
becomes
\[
(\widetilde{E}_{t})\,:\,\triangle _{\widetilde{\mathbf{g}}_{t}}\widetilde{u}%
_{t}+\mu _{t}^{2}.\widetilde{h}_{t}.\widetilde{u}_{t}=\lambda _{t}\widetilde{%
f}_{t}.\widetilde{u}_{t}^{\frac{n+2}{n-2}} 
\]
and to obtain that this equation "converges" to
\[
\bigtriangleup _{e}\widetilde{u}=K^{-2}\widetilde{u}^{\frac{n+2}{n-2}} 
\]
we need to show that $(\,\widetilde{f}_{t})$ is simply convergent to 1
(which is obvious when we have a constant function $f$ on the right handside of 
($E_{h,f,\mathbf{g}}$)). As the sequence ($\widetilde{f}_{t}$) is uniformly bounded by 1 
on $\Bbb{R}^{n}$ (considering we have extended
$\widetilde{f}_{t}$ by 0 on $\Bbb{R}^{n}\backslash B(0,\delta \mu
_{t}^{-1})$), we have, using e.g. theorem 8.25 of
Gilbard-Trudinger \cite{G-T} and Ascoli's theorem,
the existence of a function $\widetilde{u}\in C^{0}(\Bbb{R}^{n})$ such that, 
after extraction, $\widetilde{u}_{t}\stackrel{C_{loc}^{0}(\Bbb{R}%
^{n})}{\rightarrow }\,\widetilde{u}$, with $\widetilde{u}(0)=1$.

We are going to prove that $\widetilde{f}_{t}\stackrel{a.e.}{\rightarrow }1$ on 
$\Bbb{R}^{n}$ in two steps (we will prove a little bit more):

1/: There exists $\widetilde{f}\in L_{loc}^{2}(\Bbb{R}^{n})$ such that $%
\widetilde{f}_{t}\stackrel{a.e.}{\rightarrow }\widetilde{f}$ on $\Bbb{R}^{n}$

2/: $\widetilde{f}=1$ a.e. on $\Bbb{R}^{n}$\\

First step:

We have $\widetilde{f}_{t}(x)=\overline{f}_{t}(\mu _{t}x)$ and $\left| \nabla 
\overline{f}_{t}\right| \leq \frac{c}{t}$. Therefore
\[
\left| \nabla \widetilde{f}_{t}\right| \leq c.\frac{\mu _{t}}{t}\,. 
\]

We consider two cases:

a/: If ($\frac{\mu _{t}}{t}$) is bounded: Then for any compact set $%
K\subset \subset \Bbb{R}^{n}$, ($\widetilde{f}_{t}$) is bounded in $%
H_{1}^{n+1}(K)$ (where $n=\dim M)$. Thus, by compacity of the inclusion $%
H_{1}^{n+1}(K)\subset C^{0,\alpha }(K)$ for some $\alpha >0$, up to a subsequence, 
there exists $\widetilde{f}_{K}\in C^{0,\alpha }(K)$ such that
\[
\widetilde{f}_{t}\stackrel{C^{0,\alpha }(K)}{\rightarrow }\widetilde{f}_{K} 
\]
By diagonal extraction, we constuct $\widetilde{f}\in C^{0,\alpha }(%
\Bbb{R}^{n})$ such that
\[
\widetilde{f}_{t}\stackrel{C^{0,\alpha }(K^{\prime })}{\rightarrow }%
\widetilde{f} 
\]
for any compact set $K^{\prime }$ of $\Bbb{R}^{n}$, and moreover $\widetilde{f}%
\in H_{1,loc}^{n+1}(\Bbb{R}^{n}).$ So $\widetilde{f}_{t}%
\stackrel{a.e.}{\rightarrow }\widetilde{f}\,\,\,on\,\,\,\Bbb{R}^{n}\,.$

b/: If $\frac{\mu _{t}}{t}\rightarrow +\infty $ : the support of $\widetilde{f%
}_{t}$ is
\[
Supp\widetilde{f}_{t}=B(\frac{x_{0}(t)}{\mu _{t}},\frac{t}{\mu _{t}}), 
\]
where $x_{0}(t)=\exp _{x_{t}}^{-1}(x_{0})$.

If ($\frac{\left| x_{0}(t)\right| }{\mu _{t}}$) is bounded, there is after extraction a subsequence
\[
\frac{x_{0}(t)}{\mu _{t}}\rightarrow P\in \Bbb{R}^{n}; 
\]
and therefore
\[
\widetilde{f}_{t}\stackrel{C_{loc}^{0}(\Bbb{R}^{n}-\{P\})}{\rightarrow }0 
\]

If $\frac{\left| x_{0}(t)\right| }{\mu _{t}}\rightarrow \infty $, then
\[
\widetilde{f}_{t}\stackrel{C_{loc}^{0}(\Bbb{R}^{n})}{\rightarrow }0 
\]
In both cases, $\widetilde{f}_{t}\stackrel{p.p}{\rightarrow }%
0\,\,\,on\,\,\,\Bbb{R}^{n}.$

In case a/, $\widetilde{u}$ is a weak solution of 
\[
\bigtriangleup _{e}\widetilde{u}=K^{-2}\widetilde{f}\widetilde{u}^{\frac{n+2%
}{n-2}}\,\, 
\]
with $\widetilde{f}\geq 0$ as $\widetilde{f}_{t}\geq 0$, and $\widetilde{%
f}\in H_{1,loc}^{n+1}(\Bbb{R}^{n})\subset C^{0,\alpha }(\Bbb{R}^{n}).\,$

In case b/, $\,\widetilde{u}$ is a weak solution of
\[
\bigtriangleup _{e}\widetilde{u}=0.\,\, 
\]
In both cases, elliptic thory and standard regularity thorems gives the $C^{2}$ regularity of $\widetilde{u}$
, and therefore $\bigtriangleup _{e}\widetilde{u}\geq 0$. The maximum principle then shows that 
either $\widetilde{u}\equiv 0$ or $\widetilde{u}%
>0 $. But $\widetilde{u}(0)=1$ thus $\widetilde{u}>0$.

Second step:

We start using the iteration process : for some cut-off function 
$\eta $ equal to 1 near $x_{0}$, we multiply ($E_{t}$) by $\eta
^{2}u_{t}$, integrate and use the Sobolev inequality to obtain, remembering that $\lambda _{t}<K^{-2}(\stackunder{M}{Sup}%
f_{t})^{-\frac{n-2}{n}}$ and that $Sup\,f_{t}=1$:
\[
(\int_{M}(\eta u_{t})^{2^{*}})^{\frac{2}{2^{*}}}\leq \lambda
_{t}K^{2}\int_{M}\eta ^{2}f_{t}u_{t}^{2^{*}}+c\int_{Supp\,\eta
}u_{t}^{2}\,\,. 
\]
We take $\eta =1$ on $B(x_{0},\frac{3}{2}\delta )$ and $\eta =0$ on $%
M\backslash B(x_{0},2\delta )$. Then for $t$ close to 0 
\[
Supp\,f_{t}\subset B(x_{0},t)\subset B(x_{t},\delta )\subset B(x_{0},\frac{3%
}{2}\delta ) 
\]
So
\[
(\int_{B(x_{t},\delta )}u_{t}{}^{2^{*}})^{\frac{2}{2^{*}}}\leq
\int_{B(x_{t},\delta )}f_{t}u_{t}^{2^{*}}+c\int_{M}u_{t}^{2} 
\]
and after blow-up
\[
(\int_{B(0,\delta \mu _{t}^{-1})}\widetilde{u}_{t}^{2^{*}})^{\frac{2}{2^{*}}%
}\leq \int_{B(0,\delta \mu _{t}^{-1})}\widetilde{f}_{t}\widetilde{u}%
_{t}^{2^{*}}+c\int_{M}u_{t}^{2}=1+c\int_{M}u_{t}^{2}\,\,. 
\]
But $\int_{M}u_{t}^{2}\rightarrow 0$ therefore
\[
\stackunder{t\rightarrow 0}{\overline{\lim }}\int_{B(0,\delta \mu _{t}^{-1})}%
\widetilde{u}_{t}^{2^{*}}\leq 1\,\,. 
\]
Beside, we know that $\widetilde{f}_{t}\stackrel{a.e.}{\rightarrow }%
\widetilde{f}$ with $\widetilde{f}\leq 1$ and $\widetilde{u}_{t}(0)=1$.
Let suppose that there exists a set $A\subset \Bbb{R}^{n}$ with $mes(A)>0$ such that $\widetilde{f}<1$ on 
$A$ and write $\Bbb{R}^{n}=A\cup B$ with $%
\widetilde{f}=1$ a.e. on $B$. Then, as $\widetilde{f}_{t}\geq 0$ and as 
$\widetilde{u}_{t}\stackrel{C^{2}}{\rightarrow }\widetilde{u}>0$ :
\begin{eqnarray*}
1=\int_{B(0,\delta \mu _{t}^{-1})}\widetilde{f}_{t}\widetilde{u}_{t}^{2^{*}}
& &\leqslant \stackunder{t\rightarrow 0}{\overline{\lim }}\int_{B(0,\delta
\mu _{t}^{-1})\cap A}\widetilde{f}_{t}\widetilde{u}_{t}^{2^{*}}+\stackunder{%
t\rightarrow 0}{\overline{\lim }}\int_{B(0,\delta \mu _{t}^{-1})\cap B}%
\widetilde{f}_{t}\widetilde{u}_{t}^{2^{*}} \\ 
&& <\stackunder{t\rightarrow 0}{\overline{\lim }}\int_{B(0,\delta \mu
_{t}^{-1})\cap A}\widetilde{u}_{t}^{2^{*}}+\stackunder{t\rightarrow 0}{%
\overline{\lim }}\int_{B(0,\delta \mu _{t}^{-1})\cap B}\widetilde{u}%
_{t}^{2^{*}} \\ 
&& =\stackunder{t\rightarrow 0}{\overline{\lim }}\int_{B(0,\delta \mu
_{t}^{-1})}\widetilde{u}_{t}^{2^{*}}
\end{eqnarray*}
so 
\[
1<\stackunder{t\rightarrow 0}{\overline{\lim }}\int_{B(0,\delta \mu
_{t}^{-1})}\widetilde{u}_{t}^{2^{*}} 
\]
which is a contradiction, and therefore $\widetilde{f}_{t}\stackrel{a.e.}{%
\rightarrow }1$ on $\Bbb{R}^{n}$.

Thus, as we said
\[
(\widetilde{E}_{t})\,:\,\triangle _{\widetilde{\mathbf{g}}_{t}}\widetilde{u}%
_{t}+\mu _{t}^{2}.\widetilde{h}_{t}.\widetilde{u}_{t}=\lambda _{t}\widetilde{%
f}_{t}.\widetilde{u}_{t}^{\frac{n+2}{n-2}} 
\]
``converges'' to 
\[
\bigtriangleup _{e}\widetilde{u}=K^{-2}\widetilde{u}^{\frac{n+2}{n-2}} 
\]
in the sense that
\[
\,\widetilde{u}_{t}\stackrel{C_{loc}^{2}(\Bbb{R}^{n})}{\rightarrow }\,%
\widetilde{u} 
\]
where $\widetilde{u}$ is a solution of $\bigtriangleup _{e}\widetilde{u}%
=K^{-2}\widetilde{u}^{\frac{n+2}{n-2}}$. As $\widetilde{u}(0)=1$, 
\[
\widetilde{u}(x)=(1+\frac{K^{-2}}{n(n-2)}\left| x\right| ^{2})^{-\frac{n-2}{2%
}}\,\,. 
\]

Now, we can proceed exactly as in subsection 4.2.. We have:

\[
\forall R>0:\stackunder{t\rightarrow 0}{\lim }\int_{B(x_{t},R\mu
_{t})}f_{t}u_{t}^{2^{*}}dv_{\mathbf{g}}=1-\varepsilon _{R}\,\,\,where%
\,\,\,\varepsilon _{R}\stackunder{R\rightarrow +\infty }{\rightarrow }0 
\]
then

\[
\exists C>0\,\,\,such\,\,\,that\,\,\,\forall x\in M:\,d_{\mathbf{g}}(x,x_{t})^{%
\frac{n-2}{2}}u_{t}(x)\leq C. 
\]
and

\[
\forall \varepsilon >0,\exists R>0\,\,\,such\,\,\,that\,\,\,\forall
t,\,\forall x\in M:\,d_{\mathbf{g}}(x,x_{t})\geq R\mu _{t}\,\Rightarrow
\,\,d_{\mathbf{g}}(x,x_{t})^{\frac{n-2}{2}}u_{t}(x)\leq \varepsilon . 
\]

d/: We have here again the $L^{2}$-concentration:

If $\dim M\geq 4,$ 
\[
\forall \delta >0\,:\,\stackunder{t\rightarrow 0}{\lim }\frac{%
\int_{B(x_{0},\delta )}u_{t}^{2}dv_{\mathbf{g}}}{\int_{M}u_{t}^{2}dv_{%
\mathbf{g}}}=1 
\]

e/: We also have the strong estimates:
For $0<\nu <\frac{n-2}{2}$
\[
\exists C(\nu)>0\,\,\,such\,\,\,that\,\,\,\forall x\in M:\,d_{\mathbf{g}%
}(x,x_{t})^{n-2-\nu}\mu _{t}^{-\frac{n-2}{2}+\nu}u_{t}(x)\leq C, 
\]
and therefore the strong $L^{p}$-concentration:

\textbf{\ }$\forall R>0$ , $\forall \delta >0$ and $\forall p>\frac{n}{n-2}$
where $n=\dim M$:
\[
\stackunder{t\rightarrow 0}{\lim }\frac{\int_{B(x_{t},R\mu
_{t})}u_{t}^{p}dv_{\mathbf{g}}}{\int_{B(x_{t},\delta )}u_{t}^{p}dv_{\mathbf{g%
}}}=1-\varepsilon _{R}\,\,\,where\,\,\,\varepsilon _{R}\stackunder{%
R\rightarrow +\infty }{\rightarrow }0 
\]

We can now proceed with the central part of the proof of theorem 4:

We consider the euclidean Sobolev inequality and equation 
($E_{t})$ viewed in the chart $\exp _{x_{t}}^{-1}$. Using the same computations as in subsection 4.3, we get:
\begin{eqnarray*}
\int_{B(0,\delta )}\overline{h}_{t}(\eta \overline{u}_{t})^{2}dx\leq && \frac{%
1}{K(n,2){{}^{2}}(\stackunder{M}{Sup}f)^{\frac{n-2}{n}}}\int_{B(0,\delta )}%
\overline{f}_{t}\eta ^{2}\overline{u}_{t}^{2^{*}}dx\\
&&-\frac{1}{K(n,2){{}^{2}}}%
(\int_{B(0,\delta )}(\eta \overline{u}_{t})^{2^{*}}dx)^{\frac{2}{2^{*}}} \\ 
&&+C.\delta ^{-2}\int_{B(0,\delta )\backslash B(0,\delta /2)}\overline{u}%
_{t}^{2}dx+B_{t}+C_{t}
\end{eqnarray*}
with

$B_{t}=\frac{1}{2}\int_{B(0,\delta )}(\partial _{k}(\,\mathbf{g}%
\,_{t}^{ij}\Gamma (\,\mathbf{g}\,_{t})_{ij}^{k}+\partial _{ij}\,\mathbf{g}%
\,_{t}^{ij})(\eta \overline{u}_{t}^{2})dx$

$C_{t}=\left| \int_{B(0,\delta )}\eta ^{2}(\,\mathbf{g}\,_{t}^{ij}-\delta
^{ij})\partial _{i}\overline{u}_{t}\partial _{j}\overline{u}_{t}dx\right| $

$A_{t}=\frac{1}{K(n,2){{}^{2}}(\stackunder{M}{Sup}f)^{\frac{n-2}{n}}}%
\int_{B(0,\delta )}\overline{f}_{t}\eta ^{2}\overline{u}_{t}^{2^{*}}dx-\frac{%
1}{K(n,2){{}^{2}}}(\int_{B(0,\delta )}(\eta \overline{u}_{t})^{2^{*}}dx)^{%
\frac{2}{2^{*}}}$

We can write
\[
A_{t}\leq \frac{1}{K(n,2){{}^{2}}(\stackunder{M}{Sup}f_{t})^{\frac{n-2}{n}}}%
(A_{t}^{1}+A_{t}^{2}) 
\]
where $\,A_{t}^{1}=(\int_{B(0,\delta )}\overline{f}_{t}(\eta \overline{u}%
_{t})^{2^{*}}dx)^{\frac{n-2}{n}}\,-(Supf_{t}.\int_{B(0,\delta )}(\eta 
\overline{u}_{t})^{2^{*}}dx)^{\frac{n-2}{n}}\,.$

from the computation of subsection 4.3, 
\[
\stackunder{t\rightarrow 0}{\overline{\lim }}\frac{K(n,2){{}^{-2}}(%
\stackunder{M}{Sup}f_{t})^{-\frac{n-2}{n}}A_{t}^{2}+C.\delta
^{-2}\int_{B(0,\delta )\backslash B(0,\delta /2)}\overline{u}%
_{t}^{2}dx+B_{t}+C_{t}}{\int_{B(0,\delta )}\overline{u}_{t}^{2}dx}\leq \frac{%
n-2}{4(n-1)}S_{\mathbf{g}}(x_{0})+\varepsilon _{\delta } 
\]
where $\varepsilon _{\delta }\rightarrow 0$ when $\delta \rightarrow 0$.

We now consider : 
\[
\stackunder{t\rightarrow 0}{\overline{\lim }}\frac{A_{t}^{1}}{%
\int_{B(0,\delta )}\overline{u}_{t}^{2}dx} 
\]
We remark that from its definition, $f_{t}$ is decreasing when $t\rightarrow
0 $ in the sense that: 
\[
if\,\,\,t\leq t^{\prime }\,\,\,then\,\,\,f_{t}\leq f_{t^{\prime }}\,\,. 
\]
We fix a $t_{0}.$ Then, for any $t\leq t_{0}$%
\begin{eqnarray*}
\int_{B(0,\delta )}\overline{f}_{t}(\eta \overline{u}_{t})^{2^{*}}dx& & 
=\int_{B(x_{t},\delta )}f_{t}.(\eta \circ \exp
_{x_{t}}^{-1})^{2^{*}}.u_{t}^{2^{*}}.(\exp _{x_{t}}^{-1})^{*}dx \\ 
&& \leq \int_{B(x_{t},\delta )}f_{t_{0}}.(\eta \circ \exp
_{x_{t}}^{-1})^{2^{*}}.u_{t}^{2^{*}}.(\exp _{x_{t}}^{-1})^{*}dx \\ 
& &=\int_{B(0,\delta )}(f_{t_{0}}\circ \exp _{x_{t}})(\eta \overline{u}%
_{t})^{2^{*}}dx\,.
\end{eqnarray*}
We note: 
\[
\overline{f}_{t_{0},t}=f_{t_{0}}\circ \exp _{x_{t}}\, 
\]
and
\[
\widetilde{f}_{t_{0},t}=\overline{f}_{t_{0},t}\circ \psi _{\mu _{t}^{-1}}^{-1}\,. 
\]
Then : 
\begin{eqnarray*}
A_{t}^{1} & \leq (\int_{B(0,\delta )}\overline{f}_{t_{0},t}(\eta \overline{u}%
_{t})^{2^{*}}dx)^{\frac{n-2}{n}}\,-(Supf_{t}.\int_{B(0,\delta )}(\eta 
\overline{u}_{t})^{2^{*}}dx)^{\frac{n-2}{n}} \\ 
& \leq (\int_{B(0,\delta )}\overline{f}_{t_{0},t}(\eta \overline{u}%
_{t})^{2^{*}}dx)^{\frac{n-2}{n}}\,-(Supf_{t_{0}}.\int_{B(0,\delta )}(\eta 
\overline{u}_{t})^{2^{*}}dx)^{\frac{n-2}{n}}
\end{eqnarray*}
as $Supf_{t}=Supf_{t_{0}}=1=f_{t_{0}}(x_{0})\,$for all $t$.

We therefore obtain by the same method than that of section 4.3: 
\[
\stackunder{t\rightarrow 0}{\overline{\lim }}\frac{A_{t}^{1}}{%
\int_{B(0,\delta )}\overline{u}_{t}^{2}dv_{\mathbf{g}_{t}}}\leq -\frac{%
(n-2)(n-4)}{8(n-1)}\frac{\bigtriangleup _{\mathbf{g}}f_{t_{0}}(x_{0})\,}{%
f_{t_{0}}(x_{0})\,}+\varepsilon _{\delta } 
\]
and thus, after letting $\delta $ tend to 0, we obtain: 
\[
h(x_{0})\leq \frac{n-2}{4(n-1)}S_{\mathbf{g}}(x_{0})-\frac{(n-2)(n-4)}{8(n-1)%
}\frac{\bigtriangleup _{\mathbf{g}}f_{t_{0}}(x_{0})\,}{f_{t_{0}}(x_{0})\,} 
\]
But
\[
\bigtriangleup _{\mathbf{g}}f_{t}(x_{0})\,\sim +\frac{c}{t^{2}}\stackunder{%
t\rightarrow 0}{\rightarrow }+\infty 
\]
so taking $t_{0}$ close to 0 we obtain a contradiction.

This proves that we can find in the sequence ($f_{t}$) functions
whith laplacian in $x_{0}$, $\bigtriangleup _{\mathbf{g}}f_{t}(x_{0})$, 
as large as we want such that the equations: $\bigtriangleup _{%
\mathbf{g}}u+h.u=f_{t}.u^{\frac{n+2}{n-2}}$ do \textit{not} have minimizing solutions 
and therefore such that $h$ is weakly critical for $f_{t}$ and $\mathbf{g}$.
\\

Remark 1: We also have in this setting the analog of theorem 6 on the speed of convergence of $(x_{t})$ 
to $x_{0}$.

Remark 2: this can be apply to $h=cste<B_{0}K^{-2}$ or to $h=S_{\mathbf{g}}$ if $M$ is not the sphere.
\\

\textit{Second step:}

For our function $h$ such that $(h,1,\mathbf{g})$ is subcritical,
we know now that there exists a function $f$, with a laplacian as large as we want at its maximum points, such that $(h,f,\mathbf{g})$
is weakly critical. More precisely, we found a function $f$ such that:

1/: $(h,f,\mathbf{g})$ is weakly critical,

2/: $h(x_{0})>\frac{n-2}{4(n-1)}S_{\mathbf{g}}(x_{0})-\frac{(n-2)(n-4)}{%
8(n-1)}\frac{\bigtriangleup _{\mathbf{g}}f(x_{0})\,}{f(x_{0})\,}$ where

a/: $h(x_{0})>0$

b/: \{$x_{0}\}=\{x\,/\,f(x)=\stackunder{M}{Sup}f\}$ and $f(x_{0})=1$, $0\leq
f\leq 1$, $Supp\,f=B(x_{0},r)$

c/: $\nabla ^{2}f(x_{0})<0$ .

We now consider the path
\[
t\rightarrow f_{t}=(1-t).1+t.f. 
\]
Remark that for all $t$: $\bigtriangleup _{\mathbf{g}}f_{t}=t\bigtriangleup
_{\mathbf{g}}f$ and $f_{t}(x_{0})=1=\stackunder{M}{Sup}\,f_{t}$. We set
\[
\lambda _{t}=Inf\,J_{h,f_{t},\mathbf{g}}. 
\]
Then

\[
\lambda _{0}<K(n,2){{}^{-2}}(\stackunder{M}{Sup}f_{0})^{-\frac{n-2}{n}} 
\]
because ($h,1,\mathbf{g})$ is sub-critical and

\[
\lambda _{1}=K(n,2){{}^{-2}}(\stackunder{M}{Sup}f_{1})^{-\frac{n-2}{n}} 
\]
as ($h,f,\mathbf{g})$ is weakly critical. Remark that $\stackunder{M%
}{Sup}\,f_{t}\,$is always equal to 1.

Let
\[
t_{0}=Sup\{t\,/\,\lambda _{t}<K(n,2){{}^{-2}}(\stackunder{M}{Sup}f_{t})^{-%
\frac{n-2}{n}}\} 
\]
Then $0<t_{0}\leq 1$ and
\[
\lambda _{t_{0}}=K(n,2){{}^{-2}}(\stackunder{M}{Sup}f_{t_{0}})^{-\frac{n-2}{n%
}} 
\]
Before applying the method of section 4.3, we need to prove one more thing : as $h$ is weakly critical for $f_{t_{0}}$, 
we know that at the maximum point $x_{0}$ we have
\[
h(x_{0})\geq \frac{n-2}{4(n-1)}S_{\mathbf{g}}(x_{0})-\frac{(n-2)(n-4)}{8(n-1)%
}\frac{\bigtriangleup _{\mathbf{g}}f_{t_{0}}(x_{0})\,}{f_{t_{0}}(x_{0})\,} 
\]
because $\frac{\bigtriangleup _{\mathbf{g}}f_{t_{0}}(x_{0})\,}{f_{t_{0}}(x_{0})\,%
}=t_{0}\frac{\bigtriangleup _{\mathbf{g}}f(x_{0})\,}{f(x_{0})\,}$ with $%
t_{0}\leq 1$, but we need a strict inequality.

We consider the sequence ($f_{i})$, that we can construct using the first step: $f_{i}$ is such that 
($h,f_{i},\mathbf{g}$) is weakly critical with
\[
f_{i}(x_{0})=1=Supf_{i}\,\,\,et\,\,\,\bigtriangleup _{\mathbf{g}%
}f_{i}(x_{0})\rightarrow +\infty . 
\]
For each $f_{i}$, we note $t_{i}$ the ''$t_{0}$'' built above.
Therefore for any $i$ :\thinspace 
\[
h\text{ is weakly critical for }(1-t_{i}).1+t_{i}.f_{i}\text{ and }%
\mathbf{g}. 
\]
Suppose that liminf $t_{i}=0$, or, after extracting, that $%
t_{i}\rightarrow 0.$ Then, 
\[
(1-t_{i}).1+t_{i}.f_{i}\rightarrow 1 
\]
uniformly on $M$ as $0\leq f_{i}\leq 1$. But ($h,1,\mathbf{g})$ is sub-critical, thus there exists $u\in H_{1}^{2}(M)$ such that
\[
\frac{\int \left| \nabla u\right| ^{2}+\int hu^{2}}{(\int u^{2^{*}})^{\frac{2%
}{2^{*}}}}<K(n,2){{}^{-2}\,.} 
\]
But then
\[
\frac{\int \left| \nabla u\right| ^{2}+\int hu^{2}}{(\int
((1-t_{i}).1+t_{i}.f_{i})u^{2^{*}})^{\frac{2}{2^{*}}}}\rightarrow \frac{\int
\left| \nabla u\right| ^{2}+\int hu^{2}}{(\int u^{2^{*}})^{\frac{2}{2^{*}}}}%
<K(n,2){{}^{-2}} 
\]
whereas
\[
K(n,2){{}^{-2}=}K(n,2){{}^{-2}}(\stackunder{M}{Sup}%
((1-t_{i}).1+t_{i}.f_{i}))^{-\frac{n-2}{n}} 
\]
which contradict the fact that $(h,(1-t_{i}).1+t_{i}.f_{i},\mathbf{g})$ is weakly critical.

Therefore, up to extraction, $t_{i}\rightarrow t_{1}>0$

As $\bigtriangleup _{\mathbf{g}}f_{i}(x_{0})\rightarrow +\infty $, we can find $i$ 
large enough so that
\[
\frac{(n-2)(n-4)}{8(n-1)}t_{i}\frac{\bigtriangleup _{\mathbf{g}%
}f_{i}(x_{0})\,}{f_{i}(x_{0})\,}>\frac{n-2}{4(n-1)}S_{g}(x_{0})-h(x_{0})\,. 
\]
If we now denote $f$ this last function $f_{i}$ and $t_{0}$ this
\thinspace $t_{i}$, we get a path
\[
t\rightarrow f_{t}=(1-t).1+t.f 
\]
such that :

a/: $\forall t<t_{0}$ : $(h,f_{t},\mathbf{g})$ is sub-critical,

b/: ($h,f_{t_{0}},\mathbf{g})$ is weakly critical with :

b1/: \{$x_{0}\}=\{x\,/\,f_{t}(x)=\stackunder{M}{Sup}f_{t}\}$ and $%
f_{t}(x_{0})=1$ for all $t$

b2/: $h(x_{0})>\frac{n-2}{4(n-1)}S_{\mathbf{g}}(x_{0})-\frac{(n-2)(n-4)}{%
8(n-1)}\frac{\bigtriangleup _{\mathbf{g}}f_{t_{0}}(x_{0})\,}{%
f_{t_{0}}(x_{0})\,}$

b3/: $\nabla ^{2}f_{t_{0}}(x_{0})<0$

For any $t<t_{0}$ there exists a minimizing solution $u_{t}$ of the equation
\[
\bigtriangleup _{\mathbf{g}}u_{t}+h.u_{t}=\lambda _{t}.f_{t}.u_{t}^{\frac{n+2%
}{n-2}} 
\]
with $\int f_{t}u_{t}^{2^{*}}=1$. The sequence ($u_{t}$) is bounded in $%
H_{1}^{2}$ therefore 
\[
u_{t}\stackrel{H_{1}^{2}}{\stackunder{t\rightarrow t_{0}}{\rightharpoondown }%
}u 
\]
and we are once again in the situation where :

- either $u>0$ and then $u$ is a minimizing solution of $\bigtriangleup _{\mathbf{g}%
}u+h.u=\lambda _{t_{0}}f_{t_{0}}.u^{\frac{n+2}{n-2}}$, and therefore 
($h,f_{t_{0}},\mathbf{g})$ is critical.

- either $u\equiv 0$ $\,$ and once again the sequence ($u_{t}$) concentrates.
In this case, the sudy of the concentration phenomenom is easier than in the first step as the family $(f_{t})$ tend uniformly to $f$
when $t\rightarrow t_{0}$ with $Supp\,f_{t}=B(x_{0},r)$. We can find $%
\delta <r $ such that $f>0$ on $B(x_{0},\delta )$. Then there exists $c>0$ such that for any $t$ we have: 
\[
0<c\leq f_{t}\leq 1\text{ on }B(x_{0},\delta ), 
\]
Furthermore, the $f_{t}$ all reach their maximum at $x_{0}$, this maximum being always 1. 
We can then go over all the results and methods of section 4.3, 
the functions $f_{t}$ bringing this time no changes. We finally obtain
\[
h(x_{0})\leq \frac{n-2}{4(n-1)}S_{g}(x_{0})-\frac{(n-2)(n-4)}{8(n-1)}\frac{%
\bigtriangleup _{\mathbf{g}}f_{t_{0}}(x_{0})\,}{f_{t_{0}}(x_{0})\,} 
\]
thus a contradiction. Therefore ($h,f_{t_{0}},\mathbf{g})$ is critical with a minimizing solution.

This proof in fact shows the following result:

\textbf{Theorem 4':}

\textit{If }$h$\textit{\ is weakly critical for a function }$f$%
\textit{\ and a metric }$\mathbf{g},$ \textit{these datas satisfying:}

\textit{1/: }$h(x)>\frac{n-2}{4(n-1)}S_{\mathbf{g}}(x)-\frac{(n-2)(n-4)}{%
8(n-1)}\frac{\bigtriangleup _{\mathbf{g}}f(x)\,}{f(x)\,}$\textit{\ at the maximum points of }$f$

\textit{2/:}$\nabla ^{2}f(x)<0$\textit{\ at the maximum points of }$f$

\textit{3/: there exists a sequence }$f_{t}\stackrel{C^{2}}{\stackunder{%
t\rightarrow t_{0}}{\rightarrow }}f$\textit{\ with }$\stackunder{M}{Sup}%
f_{t}=\stackunder{M}{Sup}f$\textit{\ such that (}$h,f_{t},\mathbf{g})$%
\textit{\ is subcritical for }$t<t_{0}$

\textit{then (}$h,f,\mathbf{g})$\textit{\ is critical and has minimizing solutions.}

As we said in the introduction, this leads to another, dual, definition of critical functions, that is definition 3.
The natural question is then

\begin{center}
\textit{Is }$f$ \textit{critical for }$h$\textit{\ if and only if }$h$ 
\textit{is critical for }$f$ ?
\end{center}

Remark that in both cases, if $P$ is a point where $f$ is maximum on $M$\ :\textit{\ }$\frac{4(n-1)}{n-2}%
h(P)\geqslant S_{\mathbf{g}}(P)-\frac{n-4}{2}\frac{\bigtriangleup _{\mathbf{g%
}}f(P)}{f(P)}$\textit{\ }
\\

This problem seeems difficult. We prove here the result we obtain, theorem 5.

The proof starts with the following remark:
We have seen that if $h$ is weakly critical for $f$ and $\mathbf{g}$ and that $%
\bigtriangleup _{\mathbf{g}}u+h.u=f.u^{\frac{n+2}{n-2}}$ has a minimizing solution, then $h$ is 
critical for $f$ and $\mathbf{g}$. In the same way,
if $f$ is weakly critical for $h$ (in the sense that $\lambda _{h,f,%
\mathbf{g}}=K(n,2){{}^{-2}}(\stackunder{M}{Sup}f)^{-\frac{n-2}{n}}$ ) and if 
$\bigtriangleup _{\mathbf{g}}u+h.u=f.u^{\frac{n+2}{n-2}}$ has a minimizing solution $u>0$, 
then $f$ is critical for $h$. Indeed, if $f^{\prime
} $ is a function such that $Supf=Supf^{\prime }$ and $f^{\prime
}\gneqq f$, we have
\[
\int f^{\prime }u^{2^{*}}>\int fu^{2^{*}} 
\]
because $u>0$. Therefore
\[
J_{h,f^{\prime },\mathbf{g}}(u)<J_{h,f,\mathbf{g}}(u)=K(n,2){{}^{-2}}(%
\stackunder{M}{Sup}f)^{-\frac{n-2}{n}}=K(n,2){{}^{-2}}(\stackunder{M}{Sup}%
f^{\prime })^{-\frac{n-2}{n}}\,. 
\]
Using our work of section 4.3 and of this section, the proof is now short:

-If $h$ is critical for $f,$ we apply theorem 1: $\bigtriangleup _{\mathbf{g}}u+h.u=f.u^{\frac{n+2}{n-2}}$ 
has a minimizing solution, and therefore $f$ is critical for $h$.

-If $f$ is critical for $h$, these two functions (and the metric)
satisfying the hypothesis of the theorem, we have $\lambda _{h,f,%
\mathbf{g}}=K(n,2){{}^{-2}}(\stackunder{M}{Sup}f)^{-\frac{n-2}{n}},$ so $h$
is weakly critical for $f$. We then consider, for $t\stackrel{<}{%
\rightarrow }1$, the sequence
\[
t\rightarrow f_{t}=(1-t).Supf+t.f\,. 
\]
For allt $t:\,$we have $Supf_{t}=Supf$ and if $t<1$ then $f_{t}\gneqq f$. 
Therefore as $f$ is critical for $h$, by definition:
\[
\lambda _{h,f_{t},\mathbf{g}}<K(n,2){{}^{-2}}(\stackunder{M}{Sup}f_{t})^{-%
\frac{n-2}{n}}\,. 
\]
We then apply theorem 4' above to obtain that $h$ is critical for $f$ with minimizing solutions.

\section{The case of the dimension 3; ending remarks}
\subsection{The case of the dimension 3.}
We just state the results in the case of dimension 3, as they are immediate generalisations of results of O. Druet 
proved in the case where $f$ is a constant; we refer to his article for the proofs \cite{D}.
The dimension 3 requires fondamentaly the use of the Green function. We refer to the proof of proposition 8 in 
section 4.2 for the definition and the property of the Green function.
In dimension 3, for any point $x\in M$, and for $y$ close to $x$, $G_{h}$
can be writen in the following way: 
\[
G_{h}(x,y)=\frac{1}{\omega _{2}d_{\mathbf{g}}(x,y)}+M_{h}(x)+o(1)
\]
where $o(1)$ is to be taken for $y\rightarrow x$. We call $M_{h}(x)$
the mass of the Green function at $x$.

The generalisation of the results of O. Druet to the case of an arbitrary function $f$ in $\E$ gives the following:

Let $(M,\g)$ be a compact manifold of dimension 3, and let $f\in C^{\infty }(M)$ be such that $Sup f >0$.
We have the following results:
\begin{itemize}
\item For any function $h$ weakly critical for $f$ and $\g$, and for any $x\in Max f$, we have $M_{h}(x)\leq 0$.
\item For any $h\in C^{\infty }(M)$, let $B(h)=\inf \{B/$ $h+B\,\,is\,\,weakly\,\,critical\,\,for\,\,f\}$. 
Then $h+B(h)$is a critical function for $f$.
\item Let $h$ be a critical function for $f$ and $\g$. Then one of the following condition is true:
\begin{enumerate}
\item There exists $x\in Max f$ such that $M_{h}(x)=0$.
\item ($\E$) has minimizing solutions.
\end{enumerate}
\end{itemize}

Remarks:

-The condition
\[
M_{h}(x)\leqslant 0 
\]
appears as the analog of the condition
\[
\frac{4(n-1)}{n-2}h(P)\geqslant S_{\mathbf{g}}(P)-\frac{n-4}{2}\frac{%
\bigtriangleup _{\mathbf{g}}f(P)}{f(P)} 
\]
we had in dimension $\geqslant 4$. In the case $f=cst$, this condition must be satisfied on all of $M$.

-The particularity of dimension 3 is to offer critical functions of any shape, that is the meaning of the second point.

-The main difference with the case $f=cst$ studied by O. Druet is that the conditions on the mass of the Green 
function are to be considered only at the point of maximum of $f$.

\subsection{Degenerate hessian at the point of maximum and fundamental estimate.}
In theorem 6, we made the hypothesis that the hessian of $f$ is non degenerate at each of its points of maximum. 
We give here a conterexample to show that this hypothesis is necessary.
Consider the $n$-dimensional sphere $S^{n}$ with its standard metric $\mathbf{s}$.
Rewriting known results (c.f. for example \cite{H2}), there exists a unique critical function for 1 et $\mathbf{s}$,
which is
\[
h=\frac{n-2}{4(n-1)}S_{\mathbf{s}}=\frac{n-2}{4(n-1)} 
\]
and this critical function has only two type of extremal functions, the constants and the functions of the form
\[
u=a(b-\cos r)^{-\frac{n-2}{2}} 
\]
where $a\neq 0$, $b>1$, and $r$ is the geodesic distance to some fixed point of $S^{n}$. 
Consider now on $S^{n}$ a sequence of points $x_{t}$ converging to a poit $x_{0}$, and let 
\[
u_{t}=\mu _{t}^{\frac{n-2}{2}}(\mu _{t}^{2}+1-\cos r_{t})^{-\frac{n-2}{2}} 
\]
where $r_{t}(x)=d_{\mathbf{s}}(x,x_{t})$ and $\mu _{t}$ is a sequence of real converging to 0. Then 
\[
\int_{M}u_{t}^{2^{*}}dv_{\mathbf{s}}=1 
\]
and we obtain in this way a sequence of solutions of the equation 
\[
\triangle _{\mathbf{s}}u_{t}+\frac{n-2}{4(n-1)}.u_{t}=K(n,2)^{-2}u_{t}^{%
\frac{n+2}{n-2}} 
\]
where obviously the function $f=K(n,2)^{-2}$ has degenerate hessian at its maximum points ! Furthermore
\[
Sup_{M}u_{t}=u_{t}(x_{t})=\mu _{t}^{-\frac{n-2}{2}}. 
\]
This sequence concentrates and satisfies all the propositions 2 to 9 seen in section 4.2, whatever the choice of 
the sequence $x_{t}\rightarrow x_{0}$ and of the sequence $\mu _{t}\rightarrow 0$. By spherical symetry, we can 
easily find two sequences ($x_{t}$) and ($\mu _{t}$) such that
\[
\frac{d_{\mathbf{s}}(x_{t},x_{0})}{\mu _{t}}\rightarrow +\infty 
\]
by taking for example $\mu _{t}=d_{\mathbf{s}}(x_{t},x_{0})^{2}.$

Once again, it seems that the hypothesis on the hessian of $f$ ''fixes'' the position of the concentration point, 
and so imposes a speed of convergence of the sequence ($x_{t}$).
\subsection{Further questions.}
First a remark concerning the requirement of a strict inequality at the point of maximum of $f$ in theorem 1.
An easy but somewhat artificial extension of a result of hebey and Vaugon is the following:

\textit{Suppose that the manifold (}$M,\mathbf{g)}$\textit{\ is of dimension }$\geq 7$\textit{, 
and let (}$h,f,\mathbf{g)}$\textit{\ be a critical triple. Let }$T_{f}=\{x\in M/\,f(x)=Maxf\,\,and\,\,h(x)=\frac{n-2}{4(n-1)%
}S_{\mathbf{g}}(x)-\frac{(n-2)(n-4)}{8(n-1)}\frac{\bigtriangleup _{\mathbf{g}%
}f(P)}{f(P)}\}$\textit{. We suppose that }$T_{f}$\textit{\ is not dense in }$M$\textit{\ 
and that for any point }$x$\textit{\ of }$T_{f}$\textit{:}

\textit{1:The Weyl tensor vanishes on a neihbourhood of }$x$\textit{, }

\textit{2: }$\nabla ^{2}(h-\frac{n-2}{4(n-1)}S_{\mathbf{g}})$\textit{\ is not degenerate in }$x$\textit{,}

\textit{3: }$\bigtriangleup _{\mathbf{g}}f(x)=0$ if $x\in T_{f}$, \textit{and we suppose that  }$f$
\textit{\ is non degeneraye at the points of maximum which are not in }$T_{f}$\textit{.}

\textit{Then (}$h,f,\mathbf{g)}$\textit{\ has minimizing solutions.}

The main interest of this result is that we can expect existence of solutions in this case. Looking to our method, 
it seems that one need to find some other intrinsec parameters, i.e. invariant by the exponential charts $exp_{x_{t}}$.
See our thesis for more precision.
\\

Another question is the following: We saw that the study of equations $\bigtriangleup _{\mathbf{g}%
}u+hu=fu^{2^{*}-1}$ is linked to the study of the best constants in the Sobolev inclusions of $H_{1}^{2}$ in 
$L^{\frac{2n}{n-2}}$. In the same way, the study of the Sobolev inclusions of $H_{1}^{p}$
in $L^{\frac{pn}{n-p}}$, where $\frac{pn}{n-p}$ is the critical exponent,
and of the associated best constants, goes through the study of equations of the form
\[
\bigtriangleup _{p}u+hu=fu^{\frac{pn}{n-p}-1}
\]
where $\bigtriangleup _{p}u=-\nabla (\left| \nabla u\right| _{\mathbf{g}%
}^{p-2}\nabla u)$ is the p-laplacian; see for example O. Druet, E. Hebey and Z. Faget [F2]. Here also variational methods are used :
 the functional used is :
\[
I(u)=\int \left| \nabla u\right| _{\mathbf{g}}^{p}+\int hu^{p}
\]
from where we see the link with the Sobolev inclusion
\[
(\int u^{\frac{pn}{n-p}})^{\frac{n-p}{n}}\leq K(n,p)\int \left| \nabla
u\right| _{\mathbf{g}}^{p}+B\int u^{p}
\]
where $K(n,p)$ is the associated best constant. The starting point is again the following : If
\[
\stackunder{\int u^{\frac{pn}{n-p}}=1}{Inf\,}\,I(u)<K(n,p)^{-1}(Supf)^{-%
\frac{n-p}{n}}
\]
then the equation has a minimizing solution $u>0$ (knowing that the large inequality is always true). 
We therefore see that it is easy to extend the definition of critical functions to this case. It would therefore be 
interesting to know if our results can be extended to this setting.

Another question that can be asked after our work is the following :

\begin{center}
$f$\textit{ being given, is there constant critical functions ?}
\end{center}

This would give some kind of  ''best second constant $B_{0}(\mathbf{g},f)$'' linked to $f$.

At last, there is a question which emerges from our work:

\begin{center}
\textit{For a given arbitrary function }$h$\textit{\ on }$M$\textit{,
does there exist solutions (not minimizing) to the equation }$\bigtriangleup _{%
\mathbf{g}}u+hu=fu^{2^{*}-1}$\textit{\ ?}
\end{center}

Indeed, we saw that this equation has (minimizing) solutions when $h$ is sub-critical and when $h$ is 
critical with some hypothesis. However, variational methods do not give any answer is $h$ larger and different than 
some critical function, or if $\bigtriangleup _{\mathbf{g}}+h$ is not coercive. 
In this cases, if solutions exist, they cannot be minimizing. One therefore needs other methods for these cases. 
See A. Bahri \cite{B} who study the case $f=cst$ and $3\leq dimM\leq6$.
\\

\textit{The author wants to express here its deepest thanks to his thesis advisor, Michel Vaugon, for the generosity 
of his mathematical teaching, and, above all, for his friendship.}

\end{document}